\documentclass[12pt]{article}

\usepackage{latexsym,amsmath,amscd,amssymb,framed,graphics,mathrsfs}
\usepackage{enumerate}
\usepackage{tikz-cd}
\usepackage{graphicx}
\usepackage{diagrams}
\usepackage{amsthm}
\usepackage[square,numbers]{natbib}
\usepackage[colorlinks]{hyperref}
\usepackage{url}
\usepackage{colordvi}
\usepackage{color}
\newcommand*\diff{\mathop{}\!\mathrm{d}}
\usepackage[all]{xy}
\usepackage{scalerel,stackengine}
\stackMath
\newcommand\reallywidehat[1]{%
\savestack{\tmpbox}{\stretchto{%
  \scaleto{%
    \scalerel*[\widthof{\ensuremath{#1}}]{\kern-.6pt\bigwedge\kern-.6pt}%
    {\rule[-\textheight/2]{1ex}{\textheight}}
  }{\textheight}%
}{0.5ex}}%
\stackon[1pt]{#1}{\tmpbox}%
}

\makeatletter

\@addtoreset{figure}{section}
\def\thefigure{\thesection.\@arabic\c@figure}
\def\fps@figure{h, t}
\@addtoreset{table}{bsection}
\def\thetable{\thesection.\@arabic\c@table}
\def\fps@table{h, t}
\@addtoreset{equation}{section}

\makeatother

\newtheorem{theorem}{Theorem}[section]
\newtheorem{definition}[theorem]{Definition}
\newtheorem{remark}[theorem]{Remark}
\newtheorem{lemma}[theorem]{Lemma}
\newtheorem{proposition}[theorem]{Proposition}

\begin{document}

\title{Geometric, Variational, and Bracket Descriptions of Fluid Motion with Open Boundaries}

\author{Christopher Eldred\\
\texttt{celdred@sandia.gov}\\
Computer Science Research Institute\\
Sandia National Laboratories  \and
Fran\c{c}ois Gay-Balmaz\\
\texttt{francois.gb@ntu.edu.sg 
}\\
Division of Mathematical Sciences\\
Nanyang Technological University
\and
Meng Wu\\
\texttt{meng.wu@lmd.ipsl.fr 
}\\
Laboratoire de M\'et\'eorologie Dynamique\\
Ecole Normale Sup\'erieure de Paris\\
\&  Division of Mathematical Sciences\\
Nanyang Technological University
}

\date{}
\maketitle

\begin{abstract} We develop a Lie group geometric framework for the motion of fluids with permeable boundaries that extends Arnold's geometric description of fluid in closed domains. Our setting is based on the classical Hamilton principle applied to fluid trajectories, appropriately
amended to incorporate bulk and boundary forces via a Lagrange-d’Alembert approach,
and to take into account only the fluid particles present in the fluid domain at a given time. By applying a reduction by relabelling symmetries, we deduce a variational formulation in the Eulerian description that extends the Euler-Poincar\'e framework to open fluids. On the Hamiltonian side, our approach yields a bracket formulation that consistently extends the Lie-Poisson bracket of fluid dynamics and contains as a particular case a bulk+boundary bracket formulation proposed earlier. We illustrate the geometric framework with several examples, and also consider the extension of this setting to include multicomponent fluids, the consideration of general advected quantities, the analysis of higher-order fluids, and the incorporation of boundary stresses. Open fluids are found in a wide range of physical systems, including geophysical fluid dynamics, porous media, incompressible flow and magnetohydrodynamics. This new formulation permits a description based on geometric mechanics that can be applied to a broad class of models.
\end{abstract}

\tableofcontents

\section{Introduction}
Geometric mechanics formulations (ex. variational \cite{HoMaRa1998,GBYo2017a,GBYo2017b}, Hamiltonian, \cite{MoGr1980,HoKu1983,DzVo1980,MaWe1983,MaRaWe1984}, metriplectic/GENERIC \cite{Gr1984,Ka1984,Mo1984a,Mo1984b,GrOtt1997,
OtGr1997}) have proven to be a powerful tool in the development of physically based models for fluid systems, with deep connections to fundamental features such as conservation laws \cite{HoMaRaWe1985,webb2014noether,webb2014locala,webb2014localb} and involution constraints. However, existing formulations are largely restricted to the case of fixed\footnote{No exchange of momentum, since there is no movement of the boundary, although see \cite{GBMaRa2012}.}, material\footnote{No exchange of mass, also known as no-flux.} and adiabatic\footnote{No exchange of heat, and therefore energy.} boundaries. Real physical systems often have open boundaries where mass, momentum, energy, entropy and other quantities are exchanged with an environment; this exchange can extend even to the bulk of the fluid. For example, in geophysical fluid dynamics the atmosphere has surface exchanges of moisture and heat with the ocean or land surface through precipitation and evaporation, radiative interactions in the bulk, surfaces stresses due to wave-atmosphere interaction and land-atmosphere interaction; with similar processes occuring in the ocean. Plasma systems \cite{Makaremi2022,Meier2012} might have external magnetic fields such as those found in tokamaks and stellarators, and leakage or injection of charged material across (computational) boundaries. The general setting we are interested in is a fluid in a fixed\footnote{In this work we will restrict ourselves to the case of fixed boundaries, with an extension to moving boundaries \cite{GBMaRa2012} the subject of future work.} domain $\Omega$, but which can escape the boundary. The particular nature of the boundary exchange is left open, and might include no-flux or material, open and/or prescribed inflow or outflow boundary conditions.

Existing geometric mechanics formulations are incapable of directly handling such situations, especially the general case of arbitrary boundary conditions, both Lagrangian and Eulerian coordinates, and a full treatment of both variational and bracket-based approaches. Some work has been done in the context of Eulerian bracket-based formulations using GENERIC \cite{Ot2006,BaMaBeMe2018}, but it is limited to open boundary conditions, only certain fluid models, does not treat the variational case and does not handle prescribed boundary conditions (such as those needed for inflow). There is also significant work in this area that has occurred in the port-Hamiltonian literature \cite{rashad2021portA,rashad2021portB,lohmayer2022exergetic,mora2023irreversible,cardoso2021dissipative,haine2021incompressible,rashad2021exterior,califano2021geometric}, mainly in the context of single component incompressible and compressible fluids. This work is based on Stokes-Dirac structures using Eulerian coordinates, and is disconnected from first principles such as Hamilton's or Lagrange-d'Alembert. Related to this, the variational case is not treated and it is limited to certain fluid models.

This paper addresses these concerns by developing a geometric theory of fluid dynamics in open domains, applicable to both Lagrangian and Eulerian coordinates, starting from first principles and closely following the existing theory for closed domains. Furthermore, it is shown that several important extensions can be naturally incorporated within this geometric framework, such as the treatment of multicomponent fluids, general advected quantities, higher-order fluids, and boundary stresses. This is done through an application of the classical Hamilton's principle applied to Lagrangian fluid trajectories, appropriately amended to incorporate bulk and boundary exchanges via a Lagrange-d’Alembert approach, and to take into account only the fluid particles present in the fluid domain. By applying a reduction by relabelling symmetries, we deduce a variational formulation in the Eulerian description that extends the Euler-Poincar\'e framework to open fluids. On the Hamiltonian side, our approach yields a bracket formulation that consistently extends the Lie-Poisson bracket of fluid dynamics and contains as a particular case the GENERIC-based bulk+boundary bracket formulation proposed earlier \cite{Ot2006,BaMaBeMe2018}. This new theory is applicable to a wide range of fluid models with arbitrary advected quantities, and can handle arbitrary boundary conditions. 

Beyond providing a unifying geometric framework, our approach offers significant insights into the treatment of boundary conditions for open fluids, which we will elaborate on later. Specifically, it identifies non-trivial terms that arise when the relationship between (absolute) fluid momentum $\frac{\partial \mathfrak{l} }{\partial u}$ and the fluid velocity $u$ is no longer the usual  $\frac{\partial \mathfrak{l} }{\partial u} = \rho u$. Furthermore, this approach sheds light on the handling of fluids with tensorial advected quantities and higher-order fluids, where additional nontrivial contributions to the boundary conditions emerge in the context of open fluid systems.

The two diagrams below provide a schematic overview of our geometric framework, highlighting how the setting for open fluids naturally extends the existing framework for closed fluids, from both the Lagrangian and Hamiltonian perspectives, in both the material and Eulerian descriptions. All the notations used will be explained in details later.

\newpage

\textbf{(i) Overview of the case of \textit{closed} fluids:}

\begin{figure}[h!]
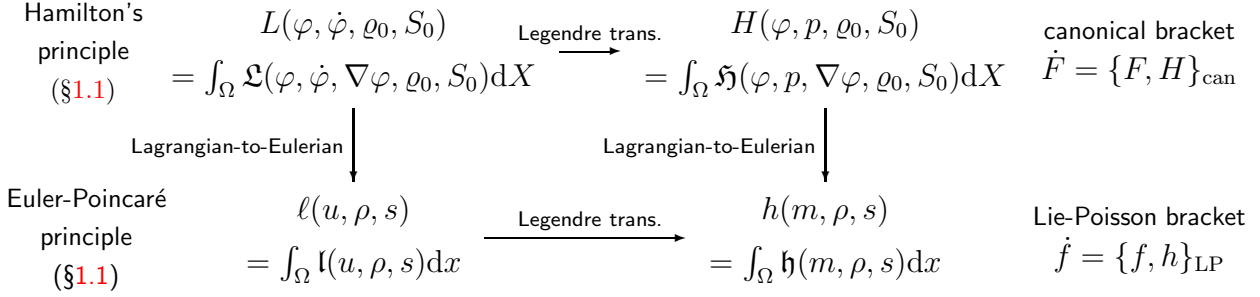

\begin{diagram}
\hspace{-0.2cm}
\begin{array}{c}
\text{\fontsize{10pt}{10pt}\selectfont \textsf{Hamilton's}}\\
\text{\fontsize{10pt}{10pt}\selectfont \textsf{principle}}\\
\text{\fontsize{10pt}{10pt}\selectfont (\S \ref{subsec_1_1})}
\end{array}
& \begin{array}{c}
\vspace{0.2cm}L(\varphi, \dot\varphi,\varrho_0,S_0)\\
 =\int_{\Omega}\mathfrak{L}(\varphi, \dot{\varphi},\nabla\varphi, \varrho_0,S_0){\rm d}X
 \end{array}& \rTo^{\tiny \hspace{-0.6cm}\textsf{Legendre trans.}\hspace{-0.6cm}} &\begin{array}{c}
\vspace{0.2cm}H(\varphi, p,\varrho_0,S_0)\\
 =\int_{\Omega}\mathfrak{H}(\varphi, p,\nabla\varphi, \varrho_0,S_0){\rm d}X
 \end{array}& \begin{array}{c} 
\textsf{\fontsize{10pt}{10pt}\selectfont canonical bracket}\\
\dot{F}=\{F,H\}_{\rm can}
\end{array}
\\
 & \dTo^{\tiny \textsf{Lagrangian-to-Eulerian}}& & \dTo^{\tiny \textsf{Lagrangian-to-Eulerian}} & \\
 \hspace{-0.4cm}
\begin{array}{c}
\textsf{\fontsize{10pt}{10pt}\selectfont Euler-Poincar\'e}\\
\textsf{\fontsize{10pt}{10pt}\selectfont principle}\\
\textsf{\fontsize{10pt}{10pt}\selectfont (\S\ref{subsec_1_1})}
\end{array}\!\!\!\! &  \begin{array}{c}
\vspace{0.2cm}\ell(u,\rho,s)\\
=\int_{\Omega}\mathfrak{l}(u,\rho,s){\rm d}x
\end{array}& \rTo^{\tiny\hspace{-0.6cm}\textsf{Legendre trans.}\hspace{-0.6cm}}& \begin{array}{c}
\vspace{0.2cm}h(m, \rho,s)\\
=\int_{\Omega}\mathfrak{h}(m,\rho,s){\rm d}x
\end{array} & \begin{array}{c} 
\textsf{\fontsize{10pt}{10pt}\selectfont Lie-Poisson bracket} \\
\dot{f}=\{f,h\}_{\rm LP}
\end{array}
\end{diagram}
\label{reversible-overview-fig}
\caption{\small Closed fluids: On the left side of the diagram, one starts with a formulation based on a \textit{Lagrangian} $L:T \operatorname{Diff}(\Omega)\rightarrow\mathbb{R} = L(\varphi, \dot\varphi,\varrho_0,S_0)$ using the classical \textit{Hamilton principle} in \textit{material} (configuration-velocity, unfortunately also known as \textit{Lagrangian}) variables, from which the variational principle in \textit{Eulerian} variables is derived through the use of \textit{relabelling symmetries} (\textit{Euler-Poincar\'e reduction}, left arrow pointing down, also known as the \textit{Lagrangian-to-Eulerian} map). Alternatively, one starts on the right side of the diagram with a formulation based on a \textit{Hamiltonian} $H:T ^*\operatorname{Diff}(\Omega)\rightarrow\mathbb{R} = H(\varphi, p,\varrho_0,S_0)$ in material (configuration-momentum) variables. These dynamics are governed by the \textit{canonical symplectic form}, with the associated \textit{canonical Poisson structure} $\{\cdot,\cdot\}_{\rm can}$. The \textit{(non-canonical) Lie-Poisson structure} $\{f,h\}_{\rm LP}$ in the \textit{Eulerian} formulation is then derived through relabelling symmetries (\textit{Lie-Poisson reduction}, right arrow pointing down, again also known as the \textit{Lagrangian-to-Eulerian} map). Under the appropriate regularity conditions, the Lagrangian and Hamiltonian approaches can be related through the \textit{Legendre transform}, at the level of either \textit{material} (top arrow) or \textit{Eulerian} (bottom arrow) variables. For example, if the Lagrangian $L$ comes from a regular Lagrangian density $\mathfrak{L}$, we can define the associated Hamiltonian function $H$.}
\end{figure}

\newpage

\textbf{(ii) Overview of the case of \textit{open} fluids:}

\begin{figure}[h!]
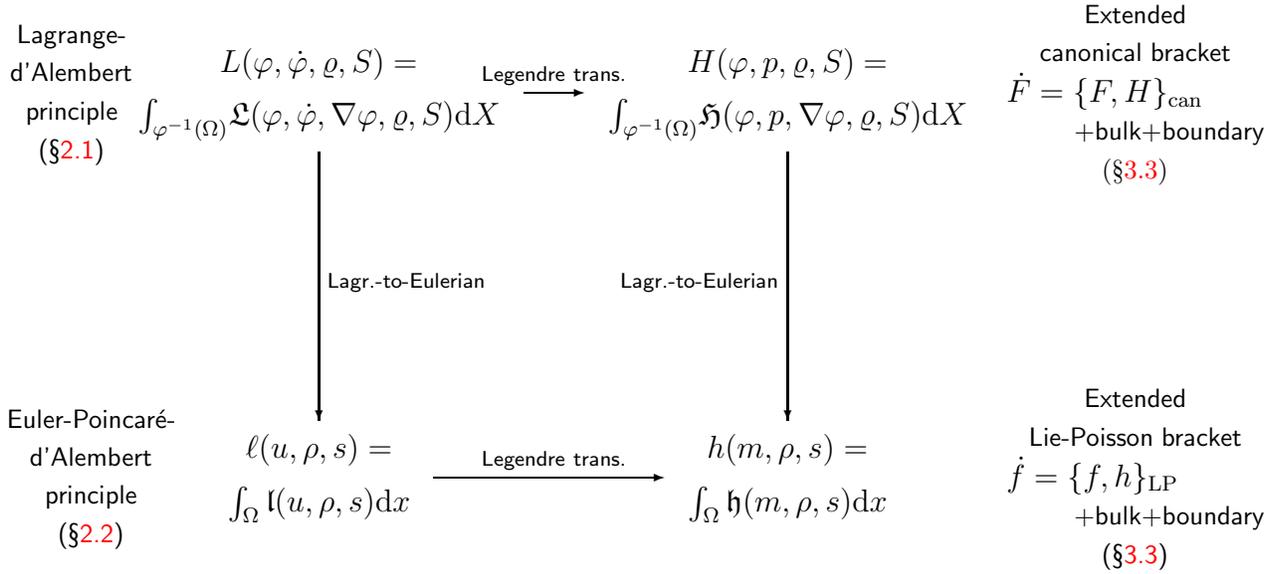

\begin{diagram}
\hspace{-0.2cm}
\begin{array}{c}
\textsf{\fontsize{10pt}{10pt}\selectfont Lagrange-}\\
\textsf{\fontsize{10pt}{10pt}\selectfont d'Alembert}\\
\textsf{\fontsize{10pt}{10pt}\selectfont principle}\\
\textsf{\fontsize{10pt}{10pt}\selectfont (\S \ref{sec_2_1})}
\end{array}\!\!\!\!
& \begin{array}{c}
\vspace{0.2cm}L(\varphi, \dot\varphi,\varrho,S)=\\
\int_{\varphi^{-1}(\Omega)}\!\mathfrak{L}(\varphi, \dot{\varphi},\nabla\varphi, \varrho,S){\rm d}X
\end{array}& \rTo^{\tiny \hspace{-0.6cm}\textsf{Legendre trans.}\hspace{-0.6cm}} &\begin{array}{c}
\vspace{0.2cm}H(\varphi, p,\varrho,S)=\\
\int_{\varphi^{-1}(\Omega)}\!\mathfrak{H}(\varphi, p,\nabla\varphi, \varrho,S){\rm d}X
\end{array}& \begin{array}{c} 
\textsf{\fontsize{10pt}{10pt}\selectfont Extended} \\
\textsf{\fontsize{10pt}{10pt}\selectfont canonical bracket}\\
\begin{array}{l}
\dot{F}=\{F,H\}_{\rm can}\\
\textsf{\fontsize{10pt}{10pt}\selectfont\qquad\;\;+\text{bulk}+\text{boundary}}
\end{array}\\
\text{\fontsize{10pt}{10pt}\selectfont (\S\ref{ELP})}
\end{array}
\\
& & &  & \\
& \dTo_{\tiny \textsf{Lagr.-to-Eulerian}}& & \dTo^{\tiny \textsf{Lagr.-to-Eulerian}} & \\
& & &  & \\
\hspace{-0.4cm}
\begin{array}{c}
\textsf{\fontsize{10pt}{10pt}\selectfont Euler-Poincar\'e-}\\
\textsf{\fontsize{10pt}{10pt}\selectfont d'Alembert}\\
\textsf{\fontsize{10pt}{10pt}\selectfont principle}\\
\textsf{\fontsize{10pt}{10pt}\selectfont (\S \ref{EPdA})}
\end{array}\hspace{-1cm}&  \begin{array}{c}
\vspace{0.2cm}\ell(u,\rho,s)=\\
\int_{\Omega}\mathfrak{l}(u,\rho,s){\rm d}x
\end{array}& \rTo^{\tiny\hspace{-0.6cm}\textsf{Legendre trans.}\hspace{-0.6cm}}& \begin{array}{c}
\vspace{0.2cm}h(m, \rho,s)=\\
\int_{\Omega}\mathfrak{h}(m,\rho,s){\rm d}x
\end{array} & \begin{array}{c} 
\textsf{\fontsize{10pt}{10pt}\selectfont Extended} \\
\textsf{\fontsize{10pt}{10pt}\selectfont Lie-Poisson bracket} \\
\begin{array}{l}
\dot{f}=\{f,h\}_{\rm LP}\\
\textsf{\fontsize{10pt}{10pt}\selectfont\qquad\;\;+\text{bulk}+\text{boundary}}
\end{array}\\
\textsf{\fontsize{10pt}{10pt}\selectfont (\S\ref{ELP})}
\end{array}
\end{diagram}
\label{open-reversible-overview-fig}
\caption{\small Open fluid: On the left side of the diagram, one again starts with a formulation based on a \textit{Lagrangian} $L:T \operatorname{Diff}(\Omega)\rightarrow\mathbb{R} = L(\varphi, \dot\varphi,\varrho_0,S_0)$ using now the \textit{Lagrange-d’Alembert principle} in material variables, incorporating \textit{bulk and boundary sources} of momentum, mass, and entropy. From this, the corresponding variational principle in Eulerian variables is derived through relabelling symmetries (Euler-Poincar\'e-d’Alembert reduction for the group $\operatorname{Diff}(\mathbb{R}^n)$ (not $\operatorname{Diff}(\Omega)$), left arrow pointing down). Alternatively, one again starts on the right side of the diagram with a formulation based on a \textit{Hamiltonian} $H:T ^*\operatorname{Diff}(\Omega)\rightarrow\mathbb{R} = H(\varphi, p,\varrho_0,S_0)$ in material variables. These dynamics are now governed by an \textit{extension} of the canonical Poisson structure $\{\cdot,\cdot\}_{\rm can}$ which incorporates the \textit{bulk and boundary contributions}. From this, an \textit{extended Lie-Poisson formulation} in the Eulerian picture is derived by relabelling symmetries, mimicking the Lie-Poisson reduction process (right arrow pointing down). Again, under the appropriate regularity conditions, the Lagrangian and Hamiltonian approaches can be related through the \textit{Legendre transform}, at the level of either \textit{material} (top arrow) or \textit{Eulerian} (bottom arrow) variables.}
\end{figure}

In summary, our geometric approach for open fluids meets the following criteria:
\begin{itemize}
\item[\rm (i)] The fluid's evolution equations can be derived equivalently from either a variational perspective or a bracket (extended Poisson) perspective.
\item[\rm (ii)] It consistently extends the geometric Lie group formulation of closed fluids on both the Lagrangian and Hamiltonian sides, in both the material and Eulerian descriptions.
\item[\rm (iii)] It incorporates the process of reduction by symmetries by the relabelling group of diffeomorphisms. 
\item[\rm (iv)] It coherently extends the fundamental variational and geometric formulations of finite dimensional mechanics, including the Hamilton and Lagrange-d’Alembert principles, as well as the canonical symplectic and Poisson structures.
\end{itemize}

\color{black}

\newpage

\paragraph{Plan of the paper.} The remainder of this paper is structured as follows. This section closes with a review of the existing general theory for closed fluids.  Section \ref{sec_L} discusses the Lagrangian variational principle for open fluids, in both Lagrangian and Eulerian coordinates. Section \ref{sec_H} discusses the Hamiltonian variational principle and associated brackets for open fluids, in both Lagrangian and Eulerian coordinates. Section \ref{extensions} introduces extensions of the theory to handle multicomponent fluids, general advected quantities, higher-order fluids and boundary stresses. Finally, conclusions are discussed in Section \ref{conclusions}. Appendix \ref{equivalance_lag_eul} demonstrates the equivalent between the Lagrangian and Eulerian coordinate versions of the theory.

\paragraph{Hamilton and Lagrange-d'Alembert principles.} Let us recall that the critical action principle of classical mechanics for a system with configuration manifold $Q$ and Lagrangian $L:TQ \rightarrow \mathbb{R} $ is given by the \textit{Hamilton principle}
\begin{equation}\label{HP} 
\delta \int_0^TL(q, \dot  q) {\rm d} t =0.
\end{equation}
It seeks for critical curves of the Lagrangian action functional, among curves with fixed endpoints, leading to $ \delta q(0)= \delta q(T)=0$. Here $TQ$ denotes the tangent bundle, or velocity space, of the manifold $Q$.
When the system is subject to an external force $F$, then the \textit{Lagrange-d'Alembert} principle must be used, in which the critical condition is amended by the virtual work of the force, giving
\begin{equation}\label{LdA} 
\delta \int_0^TL(q, \dot  q) {\rm d} t + \int_0^T \left\langle F(q, \dot  q), \delta q \right\rangle {\rm d} t=0.
\end{equation}
In order to derive the equations and boundary conditions for open fluids, we shall use a continuum version of this principle. 
As we shall see, in this continuum extension, the configuration manifold $Q$ of the system is a group of diffeomorphisms and the virtual force term has both distributed and boundary contributions.

\medskip 

We review below the geometric  variational description for closed fluids. It takes its origin in the geometric formulation due to \cite{Ar1966} of the ideal fluid equations as geodesics on the group of volume preserving diffeomorphisms of the fluid domain. For compressible fluids, one can more generally formulate the equations of motion as Euler-Lagrange equations on the group of all diffeomorphisms of the fluid domain. In the material picture, the equations thus follows from the standard Hamilton variational principle associated to the fluid Lagrangian given by its kintic minus potential energy. The induced variational principle in the Eulerian picture is best understood by using the theory of Euler-Poincar\'e reduction for systems on Lie groups \cite{HoMaRa1998}. This geometric setting is also useful to understand the Hamiltonian side and to systematically explain the occurrence of Lie-Poisson bracket structures in the Eulerian description, emerging as reduction by symmetry of canonical Poisson bracket in the material description, see \cite{MaWe1983,MaRaWe1984}.
These geometric developments only hold for closed fluids, and it is our purpose to extend them to fluids with open boundaries.

\subsection{Review of the case of closed fluids}\label{subsec_1_1}

\paragraph{The material description and Lagrangian densities.} We consider the motion of a fluid in a bounded domain $ \Omega \subset \mathbb{R} ^n$, $n=2,3$. If we assume that the motion is smooth and the fluid cannot escape the domain, then its evolution is fully described by a time dependent diffeomorphism of this domain, namely $ \varphi (t, \cdot ) \in \operatorname{Diff}( \Omega )$. Here $\operatorname{Diff}( \Omega )$ is the infinite dimensional Lie group of smooth diffeomorphisms of $ \Omega $, and it plays the role of the configuration manifold for the fluid system. Note in particular that the domain boundary is kept invariant during such motion, i.e. $ \varphi (t, \partial \Omega ) \subset  \partial  \Omega $. We shall use the standard continuum mechanics notation $x= \varphi (t,X)$ to describe the position $x \in \Omega $ of a fluid particle with label $ X \in \Omega $, which is referred to as the material or Lagrangian description of the fluid, see Remark \ref{rmk_double_domain}.

In addition to the fluid motion, for a full compressible fluid, one also needs to describe the evolution of its mass and entropy densities. For closed fluids and in absence of irreversible processes, these are described by time independent fields $ \varrho _0(X)$ and $S_0(X)$ in the material description.

As recalled earlier, following the classical Lagrangian mechanics setting, the Lagrangian function is defined on the tangent bundle $TQ$ of the configuration manifold $ Q$. For compressible fluids $Q= \operatorname{Diff}( \Omega )$ and for each fixed $ \varrho _0$, $S_0$, we have a Lagrangian function
\[
L( \cdot , \cdot , \varrho _0, S_0): T\operatorname{Diff}( \Omega ) \rightarrow \mathbb{R}, \qquad ( \varphi , \dot  \varphi ) \mapsto L( \varphi , \dot  \varphi,\varrho _0, S_0).
\]
Here $ \dot  \varphi \in T_ \varphi \operatorname{Diff}( \Omega )$ denotes an arbitrary tangent vector to the manifold $\operatorname{Diff}( \Omega )$ at $ \varphi $. It is given as a vector field on $ \Omega $ along $ \varphi $ parallel to the boundary, i.e., $\dot\varphi : \Omega \rightarrow T \Omega $ is such that $\dot  \varphi (X) \in T_{ \varphi (X)} \Omega $, for all $X \in \Omega $ and $ \dot  \varphi (X)\| \partial \Omega $ for all $X \in \partial \Omega $. In particular for $ \varphi =id$ the identity map this tangent space becomes the Lie algebra of $ \operatorname{Diff}( \Omega )$ given by the space of vector fields on $ \Omega $ tangent to the boundary: $T_ {id} \operatorname{Diff}( \Omega )= \mathfrak{X}_{\|}(\Omega )$.

In the main part of the  paper, we shall assume that the Lagrangian function is given in terms of a Lagrangian density $\mathfrak{L}$ as
\begin{equation}\label{mathfrakL} 
L( \varphi , \dot  \varphi , \varrho _0, S_0)=\int_ \Omega \mathfrak{L}( \varphi (X), \dot  \varphi(X) ,\nabla _X\varphi (X), \varrho _0(X),S_0(X)) {\rm d} X,
\end{equation} 
which is the case for the main fluid models.
We shall briefly consider a more general case in \S\ref{HOF}.

From now on we will not explicitly write the $X$-dependence of the fields in expressions such as \eqref{mathfrakL}. 

\paragraph{Hamilton's principle.} For fluid motion in a domain $ \Omega $, Hamilton's principle \eqref{HP} reads, for each fixed $ \varrho _0$ and $S_0$,
\begin{equation}\label{HP_fluid} 
\delta \int_0^T L( \varphi , \dot  \varphi , \varrho _0,S_0) {\rm d} t=0
\end{equation} 
for arbitrary variations $ \delta \varphi $ with $ \delta \varphi |_{t=0,T}=0$.

For $L$ given in terms of a density $ \mathfrak{L}$ as in \eqref{mathfrakL}, Hamilton's principle \eqref{HP_fluid} gives the Euler-Lagrange equations and boundary conditions
\begin{equation}\label{critical_cond} 
\frac{d}{dt} \frac{\partial \mathfrak{L}}{\partial \dot  \varphi} = \frac{\partial \mathfrak{L}}{\partial \varphi} - \operatorname{DIV} \frac{\partial \mathfrak{L}}{\partial \nabla \varphi } \quad\text{and}\quad \frac{\partial \mathfrak{L}}{\partial \nabla \varphi }(N,  \delta \varphi)=0, \;\forall \delta \varphi\, \| \,\partial \Omega,
\end{equation} 
where $N$ is the outward pointing unit normal vector field to $ \Omega $, seen as the domain of labels, see Remark \ref{rmk_double_domain}, and $ \operatorname{DIV}$ denotes the divergence operator with respect to the material variable $X$.

Using the local coordinate notations $X^A$ for the labels $X$ and $x^a$ for the spatial points $x$, the equations above have the local form 
\[
\frac{d}{dt} \frac{\partial \mathfrak{L}}{\partial \dot  \varphi^a } = \frac{\partial \mathfrak{L}}{\partial \varphi ^a} - \partial _A \frac{\partial \mathfrak{L}}{\partial \varphi ^a_{,A}} \quad\text{and}\quad \frac{\partial \mathfrak{L}}{\partial \varphi ^a_{,A}}N_A \delta \varphi ^a=0, \;\forall \delta \varphi\, \| \,\partial \Omega,
\]
where $ \varphi ^a_{,A}= \partial \varphi ^a / \partial X^A$. The local notation makes clear the meaning of $\frac{\partial \mathfrak{L}}{\partial \nabla \varphi }(N,  \delta \varphi)$ in \eqref{critical_cond}, namely, the evaluation of the 2-point tensor $\frac{\partial \mathfrak{L}}{\partial \nabla \varphi }$ along the material and spatial vectors $N$ and $ \delta \varphi $.

The boundary condition emerging from Hamilton's principle indicates that the traction on the boundary can only be normal to the boundary, i.e. $ \frac{\partial \mathfrak{L} }{\partial \varphi ^a_{,A}}N_A = \lambda n_a$ for some $ \lambda $ and with $n$ the outward pointing unit normal vector field to $ \Omega $ seen as the spatial domain, see Remark \ref{rmk_double_domain}. As we will see later, when $\mathfrak{L}$ satisfies the relabelling symmetry, this boundary condition becomes empty. For more general classes of Lagrangian than \eqref{mathfrakL}, such as those used in elasticity, such a boundary condition is not empty even under relabelling symmetry, see \cite{GBMaRa2012}.

Note that the writing of these equations needs the introduction of a metric, chosen here as the Euclidean metric on $ \Omega $ induced from the one on $ \mathbb{R} ^n $, see also Remark \ref{rmk_manifolds}.

\begin{remark}[Double role of $ \Omega$]\label{rmk_double_domain}\rm In the case of a closed fluid with fixed boundary, it is possible to identify the domain of labels with the spatial domain of the fluid, which we have done here by choosing them both equal to $ \Omega $. In general however they have to be thought of as being distinct \cite{MaHu1983}. For instance the metrics on the space of labels and the spatial domain do not necessarily coincide and have totally different meanings. Even though in our case we choose in both cases the Euclidean metric, we introduce two different notations for the outward pointing unit normal vector field, $N$ and $n$, according if we look at it on the domain of labels or the spatial domain. Such confusion can not arise for open fluids as these domains are distinct.
\end{remark}

\paragraph{Relabelling symmetries.}
By the relabelling symmetry (or material covariance) of fluid dynamics, the Lagrangian density $\mathfrak{L}$ must satisfy
\begin{equation}\label{relabelling_density} 
\mathfrak{L}\big( \varphi \circ \psi , \dot  \varphi \circ \psi , \nabla_X ( \varphi \circ \psi ) ,(\varrho_0 \circ \psi )J \psi,(S_0 \circ \psi )J \psi \big)= \mathfrak{L}\big( \varphi , \dot  \varphi , \nabla_X \varphi , \varrho_0 , S_0\big) \circ \psi J \psi ,
\end{equation}
for all $\psi \in \operatorname{Diff}( \Omega)$ with $ J\psi $ the Jacobian of $ \psi$, see \cite{GB2024}. 
This property turns out to be equivalent to the existence of a Lagrangian density $\mathfrak{l}(u, \rho, s )$ in the Eulerian description, such that
\begin{equation}\label{2_Lagr_densities} 
\mathfrak{L}\big( \varphi , \dot  \varphi , \nabla_X \varphi , \varrho_0  ,S_0 \big) \circ \varphi  ^{-1} J \varphi ^{-1} = \mathfrak{l} ( u  , \rho  , s )
\end{equation} 
for the Eulerian velocity, mass density, and entropy density
\[
u = \dot  \varphi \circ \varphi ^{-1} , \quad \rho  = (\varrho_0 \circ \varphi ^{-1})J \varphi ^{-1} , \quad  s  = (S_0 \circ \varphi ^{-1} ) J \varphi ^{-1},
\]
see \cite{GB2024} for details.

In terms of the associated Lagrangian function $L$ in \eqref{mathfrakL}, property \eqref{relabelling_density}  corresponds to its right $ \operatorname{Diff}( \Omega )$-invariance, namely 
\[
L( \varphi \circ \psi , \dot  \varphi \circ \psi , \varrho _0 \circ \psi J \psi , S_0 \circ \psi J \psi )=L( \varphi   , \dot  \varphi  , \varrho _0   , S_0  ) 
\]
for all $\psi \in \operatorname{Diff}( \Omega)$, and \eqref{2_Lagr_densities} yields
\[
L( \varphi , \dot  \varphi , S_0, \varrho _0)=\ell(u, \rho  , s),
\]
for the reduced Lagrangian function
\[
\ell(u, \rho  , s) = \int_ \Omega \mathfrak{l}(u, \rho  , s){\rm d} x.
\]

\paragraph{The Eulerian variational principle.}
As shown in \cite{HoMaRa1998} via Lagrange reduction by symmetry, the Hamilton principle \eqref{HP_fluid} induces the following variational principle in the Eulerian picture
\begin{equation}\label{EP_fluid}
\delta \int_0^T \ell(u, \rho  , s) {\rm d} t=0
\end{equation} 
for variations
\begin{equation}\label{EP_fluid_var} 
\delta u = \partial _t \zeta + \pounds _u \zeta , \qquad \delta \rho  = - \operatorname{div}( \rho  \zeta ), \qquad \delta s= - \operatorname{div}(s \zeta ),
\end{equation} 
where $ \zeta $ is an arbitrary time dependent vector field parallel to the boundary, vanishing at $t=0,T$.
The associated equations of motion read
\begin{equation}\label{Closed_fluid_Eulerian} 
\left\{
\begin{array}{l}
\vspace{0.2cm}\displaystyle\partial _t  \frac{\partial \mathfrak{l} }{\partial u} + \pounds _u \frac{\partial \mathfrak{l} }{\partial u} = \rho  \nabla  \frac{\partial \mathfrak{l} }{\partial \rho  }+ s\nabla  \frac{\partial \mathfrak{l} }{\partial s  } \\
\vspace{0.2cm}\displaystyle\partial _t \rho  + \operatorname{div}( \rho  u) = 0  , \qquad \partial _t s  + \operatorname{div}( s  u) = 0\\
\vspace{0.2cm}\displaystyle u \cdot n=0\quad\text{on}\quad \partial \Omega,
\end{array}
\right. 
\end{equation} 
with
\[
\pounds _u \frac{\partial \mathfrak{l} }{\partial u} = u \cdot \nabla \frac{\partial \mathfrak{l} }{\partial u} + \nabla u^\mathsf{T} \frac{\partial \mathfrak{l} }{\partial u} + \frac{\partial \mathfrak{l} }{\partial u}  \operatorname{div}u 
\]
the Lie derivative of the fluid momentum density along $u$.
We note that \eqref{EP_fluid}--\eqref{EP_fluid_var} only give the first equation. The advection equations follow from the definition $ \rho  = (\varrho _0 \circ \varphi ^{-1} )J \varphi ^{-1} $ and $ s  = (S _0 \circ \varphi ^{-1} )J \varphi ^{-1} $, while the boundary condition $u \cdot n=0$ follow since $u= \dot  \varphi \circ \varphi ^{-1} $ and $ \varphi \in \operatorname{Diff}( \Omega )$.

\begin{remark}[Equations on manifolds]\label{rmk_manifolds}\rm
In this paper we assume for simplicity that $\Omega $ is a bounded domain in the Euclidean space $ \mathbb{R} ^n$, however our approach can be extended to domains given by compact manifolds with boundary by following \cite{MaHu1983,GBMaRa2012}. In particular, an intrinsic formulation which makes explicit the dependence on Riemannian metrics on $ \Omega $ is crucial to discuss the notion of covariances. These aspects will be pursued somewhere else.
\end{remark}

\subsection{Statement of an equivalent principle for closed fluids}\label{equiv_HP}

While in the previous approach we fix the fields $ \varrho _0$ and $ S_0$ a priori and compute the critical condition by varying $ \varphi $ only, for our subsequent developments it is useful to reformulate the Hamilton principle \eqref{HP_fluid} as a critical action principle in which all the fields can be varied independently. We will utilize an extension of Hamilton's action functional recently developed in nonequilibrium thermodynamic, \cite{GBYo2017b,GBYo2019a}, which is reminiscent of a time integration by parts of a Clebsch-type variational principle, \cite{SeWh1968}.

\paragraph{Material description.} We consider the augmented Hamilton principle
\begin{equation}\label{extendedHP} 
\delta \int_0^T \left[ L( \varphi , \dot  \varphi , \varrho , S) + \varrho \dot  W + S \dot  \Gamma  \right] {\rm d} t=0
\end{equation} 
for arbitrary variations $ \delta \varphi $, $ \delta W$, $ \delta \Gamma $, $ \delta \varrho $, $ \delta S$, vanishing at $t=0,T$.

While the criticality condition with respect to $ \delta \varphi $ yields the same equations \eqref{critical_cond} as above, the other variations give, respectively,
\begin{equation}\label{additional_conditions} 
\frac{d}{dt} \varrho =0, \quad \frac{d}{dt} S=0, \quad \frac{d}{dt} W=- \frac{\partial \mathfrak{L}}{\partial \varrho } , \quad \frac{d}{dt} \Gamma =- \frac{\partial \mathfrak{L}}{\partial S}.
\end{equation} 
The first two equations recover the fact that the mass and entropy densities are constant in time, while the last two equations govern the dynamics of the additional fields $W$ and $ \Gamma $ which do not play a role in the dynamics. In this approach, $ \varrho $ and $S$ are seen as independent variables, not predetermined. The actual mass and entropy densities of the fluid $ \varrho _0(X)$ and $S_0(X)$ are inserted as initial conditions for $ \varrho (t,X)$ and $S(t,X)$. From the last two equations in \eqref{additional_conditions}, the fields $W$ and $\Gamma $ acquire the meaning of thermodynamic displacements, as defined in \cite{GBYo2017b,GBYo2019b}, which are relevant for the variational formulation of thermodynamic systems experiencing irreversible processes associated with heat and matter exchanges.

\paragraph{Eulerian description.} The Eulerian version of \eqref{extendedHP} will be derived in details later in the more general case of open fluids. For now we shall just notice that by defining the Eulerian thermodynamic displacments $w= W \circ \varphi ^{-1} $ and $ \gamma = \Gamma \circ \varphi ^{-1} $, the variational principle \eqref{extendedHP} yields
\begin{equation}\label{extendedHP_Euler} 
\delta \int_0^T \big[ \ell( u, \rho  , s) + \rho  D_t w + s D_t \gamma    \big] {\rm d} t=0
\end{equation} 
for variations $\delta u = \partial _t \zeta + \pounds _u \zeta $ and for free variations $ \delta \rho  $, $ \delta s$, $ \delta w$, and $ \delta \gamma $, vanishing at $t=0,T$ and where, as earlier, $ \zeta $ is an arbitrary time dependent vector field parallel to the boundary, vanishing at $t=0,T$. We have denoted by $D_t$ the Eulerian time derivative of scalar fields, i.e., 
\[
D_t w= \partial _t w+ u \cdot \nabla w, \qquad D_t \gamma = \partial _t \gamma +u \cdot \nabla \gamma .
\]
Now the following Eulerian version of \eqref{additional_conditions} emerges
\[
\bar D_t \rho  =0, \qquad \bar D_t s=0, \qquad D_tw= - \frac{\partial \mathfrak{l}}{\partial \rho  } , \qquad D_t \gamma = - \frac{\partial \mathfrak{l}}{\partial s } 
\]
from the variations  $ \delta w$, $ \delta \gamma $, $ \delta \rho  $, $ \delta s$, while the variations proportional to $ \zeta $ yield the fluid momentum equation, i.e., the first equation in \eqref{Closed_fluid_Eulerian}. Above we have introduced the notation $\bar D_t$ for the Eulerian time derivative of density fields i.e., 
\[
\bar D_t \rho  = \partial _t \rho  + \operatorname{div}( \rho  u) , \qquad \bar D_t s  = \partial _t s  + \operatorname{div}( s  u) .
\]
In conclusion, in a similar way with its Lagrangian version \eqref{extendedHP}  earlier, \eqref{extendedHP_Euler} is an alternative Lagrangian variational derivation of the system \eqref{Closed_fluid_Eulerian} for closed fluids, which does not a priori assume the advection equations for the mass and entropy density, but derives them as critical conditions.

The principles \eqref{extendedHP} and \eqref{extendedHP_Euler} will be appropriately extended to handle the case of an open fluid below.

\section{Lagrangian variational formulation for open fluids}\label{sec_L}

\subsection{Material description and the Lagrange-d'Alembert principle}\label{sec_2_1}

\paragraph{Lagrangian and configuration group for open fluids.} We now consider as before a domain $ \Omega $, but we allow the fluid to enter or leave the domain. In this context, while we remain focused on solving for the fluid motion within $ \Omega $, the fluid configuration map is no longer a diffeomorphism of $ \Omega $. Not only is the boundary 
$ \partial \Omega $ not preserved by the fluid motion, but the set of fluid labels corresponding to particles currently within the domain $\Omega$ also depends on the fluid configuration.

This fact is made precise by  considering as configuration group the group $ \operatorname{Diff}( \mathbb{R} ^n)$ of diffeomorphisms of the whole Euclidean space
and making the following definition. As before, we choose a Lagrangian density $ \mathfrak{L}$, and $ \varrho $, $S$ denote the mass and entropy densities in the material description.

\begin{definition}\label{def_Lopen} The Lagrangian function for an open fluid with domain $ \Omega $ and Lagrangian density $\mathfrak{L}$ is
\[
L( \cdot , \cdot , \varrho , S): T \operatorname{Diff}( \mathbb{R} ^n) \rightarrow \mathbb{R}  
\]
defined by:
\[
L( \varphi , \dot  \varphi , \varrho ,S)= \int_{ \varphi ^{-1} ( \Omega )} \mathfrak{L}( \varphi  , \dot  \varphi  , \nabla \varphi , \varrho ,S ) {\rm d} X.
\]
\end{definition}

We note the nontrivial dependence of $L$ on $ \varphi $ which arises from the fact that we integrate only on the domain $ \varphi ^{-1} ( \Omega ) \subset \mathbb{R} ^n$ which contains the labels $X $ such that $ x= \varphi (X)$ belongs to $ \Omega $. A similar Lagrangian appears in the treatment of fluid conveying tubes, see \cite{GBPu2019}, which is playing a key role in modeling the effect of boundary forces.

\begin{remark}[Case with open and closed boundaries]\rm In practical situations, the boundary can be made of several pieces some with closed and other with open boundary conditions. We shall see later how such cases can be treated, while still using $ \operatorname{Diff}( \mathbb{R} ^n)$ as the configuration Lie group.
\end{remark} 

\begin{remark}[Case of manifolds]\rm This setting can be extended to the case when $ \Omega $ is a compact manifold with smooth boundary. For this extension we use the fact that any such manifold $ \Omega $ can be properly embedded into a smooth manifold without boundary $\widetilde{ \Omega }$, which is given by a copy of its interior $ \operatorname{Int} \Omega$, see \cite{Le2013}. In this case one uses $ \operatorname{Diff}( \widetilde{ \Omega })$ instead of $ \operatorname{Diff}( \mathbb{R} ^n )$ in Definition \ref{def_Lopen}.
\end{remark}

\paragraph{Lagrange-d'Alembert principle for open fluids.} We shall now insert bulk and boundary external flows by using a continuum version of the Lagrange-d'Alembert principle \eqref{LdA}. 
We shall describe these bulk and boundary flows as being associated to the fluid momentum, the mass density, and the entropy density. Following the Lagrange-d'Alembert approach, these effects can be  naturally included by considering virtual works associated to the variation $ \delta \varphi $, $ \delta W$, and $ \delta \Gamma $, both at the interior and at the boundary. In order to be consistently paired with these variations, the distributed/bulk and boundary sources must be given, for each $ \varphi \in \operatorname{Diff}( \Omega )$, as the following maps:
\begin{itemize}
\item A distributed source of momentum $ \mathfrak{B} {\rm d} X=\mathfrak{B}( \varphi , \dot \varphi , \nabla \varphi , \varrho , S) {\rm d} X $, given as a map
\begin{equation}\label{B_def} 
\mathfrak{B}{\rm d} X:  \varphi ^{-1} ( \Omega )  \rightarrow T^* \mathbb{R} ^n \otimes \Lambda ^n (\varphi ^{-1} ( \Omega ));
\end{equation} 
\item A boundary source of momentum $ \mathfrak{J} {\rm d} X= \mathfrak{J}( \varphi , \dot \varphi , \nabla \varphi , \varrho , S) {\rm d} X$, given as a map
\begin{equation}\label{J_def} 
\mathfrak{J}{\rm d} A:  \varphi ^{-1} ( \partial \Omega ) \rightarrow T^* \mathbb{R} ^n \otimes \Lambda ^{n-1} (\varphi ^{-1} ( \partial \Omega ));
\end{equation} 
\item Distributed sources of mass $\Theta _ \varrho  {\rm d} X= \Theta _ \varrho ( \varphi , \dot \varphi , \nabla \varphi , \varrho , S) {\rm d} X $ and of entropy $\Theta _ S  {\rm d} X= \Theta _ S ( \varphi , \dot \varphi , \nabla \varphi , \varrho , S) {\rm d} X $, given as maps
\begin{equation}\label{Theta_def} 
\Theta_ \varrho {\rm d} X, \;\;\Theta_ S {\rm d} X:\varphi ^{-1} ( \Omega ) \rightarrow \Lambda ^n ( \varphi ^{-1} ( \Omega ));
\end{equation} 
\item Boundary sources of mass $\mathfrak{j}_ \varrho  ( \varphi , \dot \varphi , \nabla \varphi , \varrho , S){\rm d} A=\mathfrak{j}_ \varrho  {\rm d} A$ and of entropy $\mathfrak{j}_ S  {\rm d} A=\mathfrak{j}_ S ( \varphi , \dot \varphi , \nabla \varphi , \varrho , S) {\rm d} A$, given as maps
\begin{equation}\label{j_def} 
\mathfrak{j}_ \varrho  {\rm d} A, \;\;\mathfrak{j}_ S  {\rm d} A: \varphi ^{-1} ( \partial \Omega ) \rightarrow \Lambda ^{n-1} ( \varphi ^{-1} ( \partial \Omega )).
\end{equation}
\end{itemize}
Here $ {\rm d} A$ denotes the area element induced by $ {\rm d} X$ on $ \varphi ^{-1} ( \Omega )$ and $ \Lambda ^k(\varphi ^{-1} ( \Omega  ))$, resp., $ \Lambda ^k(\varphi ^{-1} ( \partial  \Omega  ))$, denotes the bundle of $k$-forms on $\varphi ^{-1} ( \partial \Omega  )$, resp., $\varphi ^{-1} ( \Omega  )$.

More precisely, for each $X \in \varphi ^{-1} (\Omega) $, we have 
\begin{align*}
& \mathfrak{B}(X) {\rm d} X \in T^*_{ \varphi (X)} \mathbb{R} ^n \otimes \Lambda ^n_X ( \varphi ^{-1} (\Omega ))\\
& \Theta_ \varrho (X) {\rm d} X \in\Lambda ^{n-1}_X ( \varphi ^{-1} ( \Omega ))\\
& \Theta_ S(X)  {\rm d} X \in\Lambda ^{n-1}_X ( \varphi ^{-1} ( \Omega ))
\end{align*} 
and, for each $X \in \varphi ^{-1} (\partial \Omega) $, we have 
\begin{align*}
&\mathfrak{J}(X) {\rm d} A \in T^*_{ \varphi (X)} \mathbb{R} ^n \otimes \Lambda ^{n-1}_X ( \varphi ^{-1} ( \partial \Omega ))\\
&\mathfrak{j}_ \varrho (X) {\rm d} A\in \Lambda ^{n-1}_X ( \varphi ^{-1} ( \partial \Omega ))\\
&\mathfrak{j}_ S (X) {\rm d} A \in \Lambda ^{n-1}_X ( \varphi ^{-1} ( \partial \Omega )),
\end{align*} 
with $ \Lambda ^k_X ( \varphi ^{-1} ( \Omega ))$ and $ \Lambda ^k_X ( \varphi ^{-1} ( \partial \Omega ))$ the vector fibers of the bundles at $X$.

In each case, the above clarification of the spaces in which these objects live will be useful to understand their transformation to the Eulerian frame. For this step it is important to have the dependence on $ \varphi \in \operatorname{Diff}( \mathbb{R} ^n  )$ explicitly written. Note that we have allowed these objects to have the same dependences as the Lagrangian density, i.e., on $( \varphi , \dot \varphi , \nabla \varphi , \varrho , S)$.

\medskip

Based on the preceding considerations regarding
\begin{itemize}
\item[(i)] The general statement of the Lagrange-d'Alembert principle \eqref{LdA} in classical mechanics;
\item[(ii)] The equivalent formulation of Hamilton's principle given in \S\ref{equiv_HP};
\item[(iii)] The definition of the Lagrangian for open fluids in Definition \ref{def_Lopen};
\end{itemize}
we propose the following \textit{Lagrange-d'Alembert principle for open fluids}: 
\begin{equation}\label{extended_HP}
\begin{aligned} 
&\left. \frac{d}{d\varepsilon}\right|_{\varepsilon=0}\int_0^T \!\!\int_{ \varphi _ \varepsilon ^{-1} ( \Omega )}\left[  \mathfrak{L}( \varphi _ \varepsilon , \dot  \varphi _ \varepsilon , \nabla \varphi _ \varepsilon, \varrho  _ \varepsilon , S_ \varepsilon ) + \varrho_ \varepsilon  \dot  W_ \varepsilon + S_ \varepsilon \dot  \Gamma _ \varepsilon \right]  {\rm d} X {\rm d} t\\
& \qquad  +  \int_0^T\!\!\int_{ \varphi ^{-1} ( \Omega)} \Big[ \mathfrak{B} \cdot \delta \varphi +\Theta_ \varrho   \delta W +\Theta_ S   \delta \Gamma  \Big]
{\rm d} X {\rm d} t\;\;\;\;\;\leftarrow \text{bulk contribution}\\
&  \qquad  + \int_0^T\!\!\int_{ \varphi ^{-1} (\partial\Omega)} \Big[ \mathfrak{J}  \cdot \delta \varphi +\mathfrak{j}_ \varrho  \delta W  +\mathfrak{j}_S  \delta \Gamma   \Big] {\rm d} A {\rm d} t=0\;\;\leftarrow \text{boundary contribution},
\end{aligned}
\end{equation} 
where $\delta ( \cdot )= \left. \frac{d}{d\varepsilon}\right|_{\varepsilon=0}$ and we assume $ \delta \varphi $, $ \delta W$, and $ \delta \Gamma $ vanish at $t=0,T$.
In particular, we have applied the Lagrange-d'Alembert principle \eqref{LdA} in which the action functional term $\int_0^TL(q, \dot  q) {\rm d} t$ is given by the one in \eqref{extendedHP} and the virtual force term $\int_0^T \left\langle F(q, \dot  q), \delta q \right\rangle {\rm d} t$ is given by the sum of the contribution of all the distributed and boundary forces described above. Note that in this treatment the configuration variable $q$ in \eqref{LdA} is now represented by the list of variables $( \varphi , \varrho , S, W, \Gamma )$. Note also that while the application of this type of principle implies that all variations vanish at $t=0,T$, we didn't assume that $ \delta S$ and $ \delta \varrho $ vanish. This is because, in our case, the action functional does not involve $ \dot  S$ and $ \dot  \varrho $.



\begin{proposition}\label{critical_condition} The variational principle \eqref{extended_HP} gives the following open Euler-Lagrange equations and boundary conditions:
\begin{equation}\label{EL_force} 
\left\{
\begin{array}{l}
\vspace{0.2cm}\displaystyle \frac{d}{dt}  \frac{\partial \mathfrak{L}}{\partial \dot  \varphi }  +\operatorname{DIV} \frac{\partial \mathfrak{L}}{\partial \nabla \varphi } - \frac{\partial \mathfrak{L}}{\partial \varphi } = \mathfrak{B} \\
\vspace{0.2cm}\displaystyle \frac{d}{dt}  \varrho = \Theta_ \varrho , \qquad \frac{d}{dt}  S = \Theta_ S\\
\vspace{0.2cm}\displaystyle N \cdot \Big( \frac{\partial \mathfrak{L}}{\partial \nabla \varphi }   + (\nabla_x  \varphi ^{-1} \circ \varphi )\cdot \dot \varphi \,\frac{\partial \mathfrak{L}}{\partial \dot  \varphi } + \Big( \varrho \frac{\partial \mathfrak{L}}{\partial \varrho} + S \frac{\partial \mathfrak{L}}{\partial S} - \mathfrak{L}\Big)  \nabla_x  \varphi ^{-1} \circ \varphi \Big) =-\mathfrak{J}\quad\text{on}\quad \partial \Omega \\
\displaystyle N \cdot \left((\nabla_x \varphi ^{-1} \circ \varphi )\cdot \dot \varphi\, \varrho \right)=- \mathfrak{j}_ \varrho  , \quad N \cdot \left((\nabla_x \varphi ^{-1} \circ \varphi )\cdot \dot \varphi\, S \right)=- \mathfrak{j}_ S\quad\text{on}\quad \partial \Omega ,
\end{array}
\right. 
\end{equation} 
where $N$ is the outward pointing unit normal vector field to the boundary $ \partial \varphi ^{-1} ( \Omega )$.
\end{proposition}
\noindent\textbf{Proof.} The treatment of this variational principle is not trivial since we take variations of a Lagrangian which is integrated on a domain that is $ \varepsilon $-dependent: $ \varphi _ \varepsilon ^{-1} ( \Omega )$.
The following formula will be useful
\begin{equation}\label{variat}
\left. \frac{d}{d\varepsilon}\right|_{\varepsilon=0}\int_{ \varphi _ \varepsilon ^{-1} ( \Omega )} f_ \varepsilon {\rm d} X = \int_{  \varphi ^{-1} (\partial  \Omega )} f(\delta \varphi ^{-1} \circ \varphi )\cdot N {\rm d} A + \int_{ \varphi ^{-1} ( \Omega )} \delta f {\rm d} X,
\end{equation}
for $f_ \varepsilon$ and $ \varphi _ \varepsilon $ families of functions and diffeomorphisms, both smoothly depending on $ \varepsilon$, with $f_{ \varepsilon =0}=f$ and $ \varphi _{ \varepsilon =0}= \varphi $.

In the following computation of \eqref{extended_HP}, we omit the dependence on $S$ since it is identical to the treatment of $ \varrho $ for simplicity. We compute
\begin{align*}
&\left.\frac{d}{d\varepsilon}\right|_{\varepsilon=0}\int_0^T \!\!\!\int_{ \varphi _ \varepsilon^{-1}  ( \Omega )} 
\mathfrak{L}( \varphi _ \varepsilon , \dot  \varphi _ \varepsilon , \nabla_X \varphi _ \varepsilon,  \varrho  _ \varepsilon ){\rm d} X{\rm d} t \\
& = \int_0^T\!\!\! \int_{ \varphi ^{-1} ( \Omega )} \left[ \frac{\partial \mathfrak{L}}{\partial \varphi } \cdot \delta \varphi + \frac{\partial \mathfrak{L}}{\partial \dot  \varphi } \cdot \delta \dot  \varphi + \frac{\partial \mathfrak{L}}{\partial \nabla \varphi } \cdot \nabla_X \delta \varphi + \frac{\partial \mathfrak{L}}{\partial \varrho  }  \delta \varrho \right] {\rm d} X {\rm d} t \\
& \;\;\;+ \int_0^T\!\!\!\int_{ \varphi ^{-1} ( \partial \Omega )} \mathfrak{L} \delta \varphi ^{-1} \circ \varphi \cdot N {\rm d} A {\rm d} t\\
&= \int_0^T\!\!\! \int_{ \varphi ^{-1} ( \Omega )} \left[ \Big[\frac{\partial \mathfrak{L}}{\partial \varphi }  - \frac{d}{dt}  \frac{\partial \mathfrak{L}}{\partial \dot  \varphi }  - \operatorname{DIV} \frac{\partial \mathfrak{L}}{\partial \nabla \varphi } \Big] \cdot \delta \varphi + \frac{\partial \mathfrak{L}}{\partial \varrho  } \delta \varrho \right] {\rm d} X {\rm d} t \\
&\;\;\;+ \int_0^T \frac{d}{dt}  \int_{ \varphi ^{-1} ( \Omega )} \frac{\partial \mathfrak{L}}{\partial \dot  \varphi } \cdot \delta \varphi {\rm d} X {\rm d} t\\
& \;\;\;+ \int_0^T \!\!\! \int_{  \varphi ^{-1} ( \partial \Omega )} \left[\frac{\partial \mathfrak{L}}{\partial \nabla \varphi } \cdot \delta \varphi + \mathfrak{L} \delta \varphi ^{-1} \circ \varphi - \left(\frac{\partial \mathfrak{L}}{\partial \dot  \varphi } \cdot \delta   \varphi \right) \frac{d}{dt} \varphi^{-1} \circ \varphi  \right] \cdot N {\rm d} A  {\rm d} t \\
&= \int_0^T\!\!\! \int_{ \varphi ^{-1} ( \Omega )} \left[ \Big[\frac{\partial \mathfrak{L}}{\partial \varphi }  - \frac{d}{dt}  \frac{\partial \mathfrak{L}}{\partial \dot  \varphi }  - \operatorname{DIV} \frac{\partial \mathfrak{L}}{\partial \nabla \varphi } \Big] \cdot \delta \varphi + \frac{\partial \mathfrak{L}}{\partial \varrho  } \delta \varrho \right] {\rm d} X {\rm d} t \\
& \;\;\;+ \int_0^T\!\!\! \int_{  \varphi ^{-1} ( \partial  \Omega )} \left[ \Big[ \frac{\partial \mathfrak{L}}{\partial \nabla \varphi }  - \mathfrak{L} \nabla_x  \varphi ^{-1} \circ \varphi + ( \nabla_x  \varphi ^{-1} \circ \varphi) \cdot \dot \varphi \,\frac{\partial \mathfrak{L}}{\partial \dot  \varphi } \Big]\cdot \delta \varphi  \right] \cdot N {\rm d} A  {\rm d} t,
\end{align*}
where we used $\delta \varphi ^{-1} = - \nabla_x  \varphi ^{-1} \cdot \delta \varphi \circ \varphi ^{-1}$ as well as $\delta \varphi_0 = \delta \varphi_T = 0$.
We also have
\begin{align*}
&\left.\frac{d}{d\varepsilon}\right|_{\varepsilon=0}\int_0^T\!\!\! \int_{ \varphi _ \varepsilon^{-1}  ( \Omega )} 
\varrho _ \varepsilon  \dot  W _ \varepsilon {\rm d} X {\rm d} t\\
&=\int_0^T\!\!\!\int_{ \varphi ^{-1} ( \Omega )} \left[  \dot  W \delta \varrho + \frac{d}{dt} (\varrho \delta W)- \dot \varrho \delta W \right] {\rm d} X {\rm d} t + \int_0^T\!\!\!\int_{  \varphi ^{-1} ( \partial \Omega )}  \varrho \dot W \delta \varphi ^{-1} \circ \varphi \cdot N {\rm d} A  {\rm d} t  \\
& = \int_0^T \!\!\!\int_{ \varphi ^{-1} ( \Omega )} \left[  \dot  W \delta \varrho - \dot \varrho \delta W \right] {\rm d} X {\rm d} t + \int_0^T\!\!\! \int_{  \varphi ^{-1} (\partial  \Omega )} \left[ \varrho \dot W \delta \varphi ^{-1} \circ \varphi - \varrho \delta W \frac{d}{dt}  \varphi ^{-1} \circ \varphi  \right]\cdot N {\rm d} A  {\rm d} t \\
& = \int_0^T\!\!\! \int_{ \varphi ^{-1} ( \Omega )} \left[  \dot  W \delta \varrho - \dot \varrho \delta W \right] {\rm d} X {\rm d} t \\
&\;\;+ \int_0^T\!\!\! \int_{  \varphi ^{-1} ( \partial \Omega )} \left[ -\varrho \dot W  (\nabla_x \varphi ^{-1} \circ \varphi) \cdot \delta \varphi + ( \nabla_x  \varphi ^{-1} \circ \varphi) \cdot \dot \varphi\, \varrho \delta W \right]\cdot N {\rm d} A  {\rm d} t.
\end{align*}

From the variations $ \delta \varrho $, we get the condition
\[
\dot W = - \frac{\partial \mathfrak{L}}{\partial \varrho}.
\]
Then, by using this condition, and collecting the terms proportional to the variations $ \delta \varphi $ and $ \delta W$ at the interior and at the boundary, we get the system \eqref{EL_force}. $\qquad\blacksquare$

\begin{remark}[Local expression] \rm For further understanding of the notations, it is useful to rewrite the equations \eqref{EL_force} locally.  Using as before $X^A$ and $x^a$ for the material and spatial coordinates, we have
\begin{equation}\label{EL_force_local} 
\left\{
\begin{array}{l}
\vspace{0.2cm}\displaystyle \frac{d}{dt}  \frac{\partial \mathfrak{L}}{\partial \dot  \varphi ^a}  +  \partial _A\frac{\partial \mathfrak{L}}{\partial  \varphi_{,A}^a } - \frac{\partial \mathfrak{L}}{\partial \varphi ^a} = \mathfrak{B}_a \\
\vspace{0.2cm}\displaystyle \frac{d}{dt}  \varrho = \Theta_ \varrho, \qquad \frac{d}{dt}  S = \Theta_ S \\
\vspace{0.2cm}\displaystyle N_A  \Big( \frac{\partial \mathfrak{L}}{\partial  \varphi_{,A}^a }   + \big((\varphi ^{-1})_{,b}^A \circ \varphi \big)  \dot \varphi^b \frac{\partial \mathfrak{L}}{\partial \dot  \varphi ^a} + \Big( \varrho \frac{\partial \mathfrak{L}}{\partial \varrho} +S \frac{\partial \mathfrak{L}}{\partial S} - \mathfrak{L} \Big) (\varphi ^{-1})^A_{,a} \circ \varphi \Big) =-\mathfrak{J}_a\;\text{on}\; \partial \Omega \\
\displaystyle N_A \big( (\varphi ^{-1} )^A_{,b}\circ \varphi \big)  \dot  \varphi^b \varrho =- \mathfrak{j}_ \varrho  \quad\text{on}\quad \partial \Omega , \quad N_A  \big((\varphi ^{-1} )^A_{,b}\circ \varphi  \big) \dot  \varphi^b S =- \mathfrak{j}_ S  \quad\text{on}\quad \partial \Omega .
\end{array}
\right. 
\end{equation} 
\end{remark}
 
\paragraph{Example: the compressible Euler fluid with open boundaries.} The Lagrangian density is given as
\begin{equation}\label{Euler_mathfrakL} 
\mathfrak{L}( \varphi , \dot  \varphi , \nabla \varphi , \varrho , S  )= \frac{1}{2} \varrho | \dot  \varphi | ^2 - \varepsilon \Big( \frac{ \varrho }{J \varphi }, \frac{ S }{J \varphi } \Big) J \varphi,
\end{equation} 
with $ \varepsilon$ the energy density of the fluid, expressed in terms of the Eulerian mass and entropy densities $ \rho  $ and $s$. We shall denote by $ \frac{\partial \varepsilon }{\partial \rho  } $ and $ \frac{\partial \varepsilon }{\partial s} $ the partial derivative, even though $ \rho  $ and $ s$ do not appear in the material description (recall that $ \rho  = (\varrho \circ  \varphi ) J \varphi $ and $ s  = (S \circ  \varphi ) J \varphi $).
The partial derivatives are: 
\begin{equation}\label{partial_der} \begin{aligned}
&\frac{\partial \mathfrak{L}}{\partial \varphi }=0, \quad \frac{\partial \mathfrak{L}}{\partial \dot  \varphi }= \varrho \dot  \varphi , \quad \frac{\partial \mathfrak{L}}{\partial \varrho  }=\frac{1}{2} | \dot  \varphi | ^2 - \frac{\partial \varepsilon }{\partial \rho } , \quad 
\frac{\partial \mathfrak{L}}{\partial S  }= - \frac{\partial \varepsilon }{\partial s } \\
&\frac{\partial \mathfrak{L}}{\partial \nabla \varphi }= \Big( \frac{\partial \varepsilon }{\partial \rho } \frac{\varrho}{J \varphi} + \frac{\partial \varepsilon }{\partial s } \frac{S}{J \varphi} - \varepsilon \Big) J \varphi \left( \nabla_x\varphi^{-1} \circ \varphi \right)=(p \circ \varphi )J \varphi \left( \nabla_x\varphi^{-1} \circ \varphi \right).
\end{aligned}
\end{equation}
In the last equality, we used the definition of the pressure in terms of the internal energy density as
\[
p= \frac{\partial \varepsilon }{\partial \rho  } ( \rho  , s)\rho  + \frac{\partial \varepsilon }{\partial s}( \rho  , s) s- \varepsilon ( \rho  , s)
\]
which, in material coordinates, results in
\[
p \circ  \varphi = \frac{\partial \varepsilon }{\partial \rho } \Big( \frac{ \varrho }{J \varphi }, \frac{ S }{J \varphi } \Big)\frac{\varrho}{J \varphi} + \frac{\partial \varepsilon }{\partial s } \Big( \frac{ \varrho }{J \varphi }, \frac{ S }{J \varphi } \Big)\frac{S}{J \varphi} - \varepsilon \Big( \frac{ \varrho }{J \varphi }, \frac{ S }{J \varphi } \Big).
\]
Then, we have 
\[
\operatorname{DIV} \left(p \circ \varphi  \left( \nabla\varphi^{-1} \circ \varphi \right) J \varphi\right) = \nabla_X (p \circ \varphi) \cdot  \left( \nabla\varphi^{-1} \circ \varphi \right)J \varphi =( \nabla_x p \circ \varphi) J \varphi,
\]
by using the Piola identity 
$\operatorname{DIV}(J \varphi \left( \nabla _x\varphi^{-1} \circ \varphi \right))=0$. With this, the first equation of \eqref{EL_force} reduces to 
\[
\frac{d}{dt} ( \varrho \dot  \varphi )+ (\nabla _x p \circ \varphi )J \varphi = \mathfrak{B},
\]
which is the fluid momentum equation in the material description.

By using \eqref{partial_der} remarkable cancellations occur in the third equation of \eqref{EL_force} yielding simply
\[
N \cdot \left[(\nabla_x \varphi^{-1} \circ \varphi) \cdot \dot \varphi\right] \varrho \dot \varphi = -\mathfrak{J}.
\]
We shall mention later the origin of these cancellations, see Remark \ref{cancellation}. We have thus obtained the following statement.

\begin{proposition} The equation for an open compressible Euler fluid in material coordinates are given by
\begin{equation}\label{EL_force_EULER} 
\left\{
\begin{array}{l}
\vspace{0.2cm}\displaystyle  \frac{d}{dt} ( \varrho \dot  \varphi )+ (\nabla _x p \circ \varphi )J \varphi = \mathfrak{B} \\
\vspace{0.2cm}\displaystyle \frac{d}{dt}  \varrho = \Theta_ \varrho, \qquad \frac{d}{dt} S = \Theta_ S \\
\vspace{0.2cm}\displaystyle N \cdot \left[(\nabla_x \varphi^{-1} \circ \varphi) \cdot \dot \varphi\right] \varrho \dot \varphi = -\mathfrak{J} \quad\text{on}\quad \partial \Omega\\
\displaystyle N \cdot \left[(\nabla_x \varphi ^{-1} \circ \varphi) \cdot \dot \varphi \right] \varrho =- \mathfrak{j}_ \varrho  , \quad N \cdot \left[(\nabla_x \varphi ^{-1} \circ \varphi) \cdot \dot \varphi \right] S =- \mathfrak{j}_ S\quad\text{on}\quad \partial \Omega.
\end{array}
\right. 
\end{equation}

\end{proposition} 

\paragraph{Balance laws in the material description.} The seemingly involved expressions of the boundary conditions in \eqref{EL_force} in the material description, see also \eqref{EL_force_EULER}, is due to the fact that the domain of labels is time dependent. These expressions naturally appear in the boundary integrals when computing the balance laws.

For instance, for the balance of mass, one notes
\begin{align*} 
\frac{d}{dt} \int_{ \varphi ^{-1} ( \Omega )} \varrho \,{\rm d} X &= \int_{\varphi ^{-1} ( \Omega )} \frac{d}{dt} \varrho \,{\rm d} X + \int_{ \varphi ^{-1}  ( \partial \Omega )} \frac{d}{dt} ( \varphi ^{-1} ) \circ \varphi  \cdot N {\rm d} A \\
&= \int_{\varphi ^{-1} ( \Omega )} \frac{d}{dt} \varrho \,{\rm d} X - \int_{ \varphi ^{-1}  ( \partial \Omega )} N \cdot \left[(\nabla_x \varphi ^{-1} \circ \varphi) \cdot \dot \varphi \right] \varrho \, {\rm d} A\\
&= \int_{\varphi ^{-1} ( \Omega )} \Theta _ \varrho \,{\rm d} X+  \int_{\varphi ^{-1} ( \partial \Omega )} \mathfrak{j}_ \varrho {\rm d} A,
\end{align*} 
which used the open boundary condition for the mass density in \eqref{EL_force} in the last equality. Similarly, for the entropy density, one has
\[
\frac{d}{dt} \int_{ \varphi ^{-1} ( \Omega )} S \,{\rm d} X  = \int_{\varphi ^{-1} ( \Omega )} \Theta _ S \,{\rm d} X+  \int_{\varphi ^{-1} ( \partial \Omega )} \mathfrak{j}_S {\rm d} A.
\]
This is consistent with the interpretation of $ \Theta _ \varrho $, $ \Theta _ S$ as the distributed/bulk source of mass and entropy, and $ \mathfrak{j}_ \varrho $ and $ \mathfrak{j} _S$ as the boundary sources of mass and entropy.

Regarding the energy balance, we recall that the total energy density in the material description, for a given Lagrangian density $\mathfrak{L}$ is found as 
\[
\mathfrak{E}= \frac{\partial \mathfrak{L} }{\partial \dot  \varphi } \cdot \dot  \varphi - \mathfrak{L} .
\]
From \eqref{EL_force}, we get the local energy balance
\begin{equation}\label{local_mathfrakE} 
\frac{d}{dt} \mathfrak{E}= - \operatorname{DIV} \Big( \frac{\partial \mathfrak{L} }{\partial \nabla \varphi } \cdot \dot  \varphi \Big) + \dot  \varphi \cdot \mathfrak{B}   - \frac{\partial \mathfrak{L} }{\partial \varrho }\Theta _ \rho  - \frac{\partial \mathfrak{L} }{\partial S} \Theta _S.
\end{equation} 
Hence, the total energy balance follows from the three boundary conditions in \eqref{EL_force} as
\[
\frac{d}{dt} \int_{ \varphi ^{-1} ( \Omega )} \mathfrak{E} {\rm d}X= \int_ { \varphi ^{-1} ( \Omega )} \left[ \dot  \varphi \cdot \mathfrak{B} - \Theta _ \varrho \frac{\partial \mathfrak{L} }{\partial \varrho } - \Theta _ S \frac{\partial \mathfrak{L} }{\partial S } \right] {\rm d} X -  \int_{ \varphi ^{-1} ( \partial \Omega )} \left[ \dot  \varphi \cdot \mathfrak{J} -\frac{\partial \mathfrak{L} }{\partial \varrho } \mathfrak{j}  _ \rho  - \frac{\partial \mathfrak{L} }{\partial S} \mathfrak{j} _S \right] {\rm d} A,
\]
which show the change of energy due to bulk (first integral) and boundary (second integral) contributions. In particular, the somehow involved expression of the first boundary condition in \eqref{EL_force} is used to express the boundary integral arising from the divergence term in \eqref{local_mathfrakE} thereby yielding the simple final form.

\subsection{Eulerian variational principle}\label{EPdA}

As recalled in \S\ref{subsec_1_1} for the case of closed fluids, the variational formulation in the Eulerian frame is best explained via the process of Lagrangian reduction by symmetry. 
In our case, the fluid Lagrangian is a function 
\[
L: T \operatorname{Diff}( \mathbb{R} ^n) \times \Omega  ^n( \mathbb{R} ^n) \times \Omega  ^n( \mathbb{R} ^n) \rightarrow \mathbb{R},
\]
\[
L( \varphi , \dot  \varphi , \varrho , S)= \int_{ \varphi ^{-1} ( \Omega )} \mathfrak{L} \big( \varphi (X), \dot  \varphi(X), \nabla_X  \varphi(X) , \varrho (X), S(X)\big) {\rm d} X
\]
and the relabelling symmetry is given as the $\operatorname{Diff}( \mathbb{R} ^n )$-invariance
\begin{equation}\label{relabelling_symm} 
L\big( \varphi \circ \psi , \dot  \varphi \circ \psi , (\varrho \circ \psi) J \psi,(S \circ \psi) J \psi \big)=L( \varphi , \dot  \varphi, \varrho , S ),\quad \text{for all $ \psi \in \operatorname{Diff}( \mathbb{R} ^n )$}.
\end{equation} 
From this, $L$ induces a unique reduced Lagrangian
\begin{equation}\label{red_ell} 
\ell: \mathfrak{X}( \mathbb{R} ^n ) \times  \Omega ^n( \mathbb{R} ^n) \times \Omega  ^n( \mathbb{R} ^n) \rightarrow \mathbb{R} 
\end{equation} 
defined by $\ell(u, \rho  ,s)=L( \varphi , \dot  \varphi , \varrho , S) $, where the variables on the right hand side are such that $u= \dot  \varphi \circ \varphi  ^{-1} $, $ \rho  = (\varrho \circ \varphi ^{-1} )J \varphi ^{-1} $, and $ S  = (S \circ \varphi ^{-1} )J \varphi ^{-1} $. Let us emphasize that even though the Lagrangian is not integrated over the whole space $ \mathbb{R} ^n $, it is still $ \operatorname{Diff}( \mathbb{R} ^n )$-invariant, so that the ideas of Lagrangian reduction still apply in the open case.
We refer to Remark \ref{remark_L_red} for the link with the theory of Lagrangian reduction theory by symmetry. See Remark \ref{remark_quotient_space_H} for similar comments on the Hamiltonian side.

\paragraph{Relabelling symmetries.} For our development it is useful to put the emphasis on the symmetries of the Lagrangian density rather than those of the Lagrangian function. We will say that the Lagrangian density $ \mathfrak{L} $ is materially covariant if it satisfies
\begin{equation}\label{mat_cov} 
\mathfrak{L}\big( \varphi \circ \psi , \dot  \varphi \circ \psi , \nabla_X ( \varphi \circ \psi ) ,(\varrho \circ \psi )J \psi,(S \circ \psi )J \psi \big)= \mathfrak{L}\big( \varphi , \dot  \varphi , \nabla_X \varphi , \varrho , S\big) \circ \psi J \psi ,
\end{equation} 
for all $\psi \in \operatorname{Diff}( \mathbb{R} ^n )$, see \cite{GB2024}. This property is the natural extension of the material covariance property of stored energy functions in elasticity \cite{MaHu1983,GBMaRa2012} to the case of general Lagrangian density. It is easy to see that it implies the $ \operatorname{Diff}( \mathbb{R} ^n )$-invariance \eqref{relabelling_symm} of the associated Lagrangian function $L$. Also, the symmetry \eqref{mat_cov} is equivalent to 
the existence of a Lagrangian density $\mathfrak{l}(u, \rho, s )$ in the Eulerian description, such that
\begin{equation}\label{mathfrak_ell} 
\mathfrak{L}\big( \varphi , \dot  \varphi , \nabla_X \varphi , \varrho  ,S \big) \circ \varphi  ^{-1} J \varphi ^{-1} = \mathfrak{l} \big(u, \rho  , s \big).
\end{equation} 
where $u = \dot  \varphi \circ \varphi ^{-1}$, $\rho  = (\varrho \circ \varphi ^{-1})J \varphi ^{-1}$, $s  = (S \circ \varphi ^{-1} ) J \varphi ^{-1}$. From this, the reduced Lagrangian $\ell$ in \eqref{red_ell} associated to $L$ is given in terms of $\mathfrak{l}$
\[
\ell: \mathfrak{X} ( \mathbb{R} ^n ) \times \Omega  ^n( \mathbb{R}   ^n ) \times \Omega  ^n( \mathbb{R}   ^n ) \rightarrow \mathbb{R}, \quad \ell(u, \rho  , s)= \int_ \Omega \mathfrak{l}\big(u(x), \rho (x) , s(x)\big) {\rm d} x .
\]
It will be important to distinguish between the Lagrangian function $\ell$ and the Lagrangian density $ \mathfrak{l}$. For instance, we would like to emphasize that the Lagrangian $\ell$ is defined on the whole Lie algebra $ \mathfrak{X} ( \mathbb{R} ^n)$, but its expression effectively depends only on the restriction $u|_{  \Omega }$ of $ u \in \mathfrak{X} ( \mathbb{R} ^n)$ to $ \Omega $. Similarly for the variables $\rho  $ and $s$.

\begin{remark}[Properties of $\mathfrak{L}$]\rm While the properties \eqref{mat_cov} and \eqref{mathfrak_ell} are written, for convenience, in terms of fields, such as $\varphi $, $ \varrho $, and $u$, these properties only depends on the point values of such fields, as a careful examination shows, so that they are really properties of the densities $\mathfrak{L}$ and $\mathfrak{l}$, see  \cite{GB2024} for details.
\end{remark}

\begin{remark}[Property of materially covariant Lagrangian densities]\label{cancellation}\rm We noticed earlier a drastic simplification of one of the boundary condition in the case of the compressible Euler fluid. It is not associated to this particular example, but follows from the material covariance \eqref{mat_cov} of the Lagrangian density. Indeed, the infinitesimal property associated to the material covariance \eqref{mat_cov} yields the following relation between the partial derivatives of $ \mathfrak{L} $: 
\[
\frac{\partial \mathfrak{L} }{\partial \nabla \varphi } + \Big(\frac{\partial \mathfrak{L} }{\partial \varrho   }  \varrho +\frac{\partial \mathfrak{L} }{\partial S  }S - \mathfrak{L} \Big) \nabla _x \varphi ^{-1} \circ \varphi =0,
\]
see, e.g., \cite{GB2024}. In this case, the first boundary condition in \eqref{EL_force} reduces to just
\[
\Big(N \cdot (\nabla_x  \varphi ^{-1} \circ \varphi )\cdot \dot \varphi \Big)\frac{\partial \mathfrak{L}}{\partial \dot  \varphi }   =-\mathfrak{J},\quad\text{i.e.,}\quad 
N_A    \big((\varphi ^{-1})_{,b}^A \circ \varphi \big)  \dot \varphi^b \frac{\partial \mathfrak{L}}{\partial \dot  \varphi ^a}  =-\mathfrak{J}_a.
\]
\end{remark} 

\medskip

In a similar way with the Lagrangian density, the bulk and boundary fluxes must also satisfy symmetries in order to admit a purely Eulerian formulation. They are given by
\begin{equation}\label{B_Theta}
\begin{aligned} 
&\mathfrak{B}\big( \varphi \circ \psi , \dot  \varphi \circ \psi , \nabla_X ( \varphi \circ \psi ) ,(\varrho \circ \psi )J \psi,(S \circ \psi )J \psi \big)= \mathfrak{B}(\varphi , \dot \varphi , \nabla \varphi , \varrho , S)\circ \psi J \psi\\
&\Theta_ \varrho  \big( \varphi \circ \psi , \dot  \varphi \circ \psi , \nabla_X ( \varphi \circ \psi ) ,(\varrho \circ \psi )J \psi,(S \circ \psi )J \psi \big)= \Theta _ \varrho (\varphi , \dot \varphi , \nabla \varphi , \varrho , S) \circ \psi J \psi,
\end{aligned}
\end{equation} 
similarly for $\Theta _S$, and by
\begin{equation}\label{J_j} 
\begin{aligned}
&\mathfrak{J}\big( \varphi \circ \psi , \dot  \varphi \circ \psi , \nabla_X ( \varphi \circ \psi ) ,(\varrho \circ \psi )J \psi,(S \circ \psi )J \psi \big)= \mathfrak{J}(\varphi , \dot \varphi , \nabla \varphi , \varrho , S)\circ \psi _ \partial J \psi_ \partial \\
&\mathfrak{j} _ \varrho  \big( \varphi \circ \psi , \dot  \varphi \circ \psi , \nabla_X ( \varphi \circ \psi ) ,(\varrho \circ \psi )J \psi,(S \circ \psi )J \psi \big)= \mathfrak{j}  _ \varrho (\varphi , \dot \varphi , \nabla \varphi , \varrho , S) \circ \psi _ \partial J \psi_ \partial ,
\end{aligned} 
\end{equation}
similarly for $ \mathfrak{j}_ S$, for all $ \psi \in \operatorname{Diff}( \mathbb{R} ^n )$. We used the notation $ \psi _ \partial $ for the restriction of the diffeomorhism $ \psi $ to the boundary, and $J \psi _ \partial $ for the Jacobian of $ \psi _ \partial $. These properties are equivalent to the existence of the following Eulerian objects:
\begin{itemize}
\item a distributed source of momentum
\[
b(u, \rho  , s) {\rm d} x : \Omega  \rightarrow T^* \mathbb{R} ^n \otimes \Lambda ^n ( \Omega );
\]
\item a boundary source of momentum
\[
\mathfrak{j} (u, \rho  , s) {\rm d} a:   \partial \Omega  \rightarrow T^* \mathbb{R} ^n \otimes \Lambda ^{n-1} ( \partial \Omega );
\]
\item distributed sources of mass and entropy
\[
\theta _ \rho ( u, \rho  , s){\rm d} x, \;\;\theta _s  ( u, \rho  , s){\rm d} x:  \Omega  \rightarrow \Lambda ^n ( \Omega );
\]
\item boundary sources of mass and entropy
\[
j_ \rho (u, \rho  , s) {\rm d} a, \;\;j_ s (u, \rho  , s) {\rm d} a:  \partial \Omega \rightarrow \Lambda ^{n-1} (  \partial \Omega );
\]
\end{itemize}
such that
\begin{equation}\label{B_red}
\begin{aligned}
&\mathfrak{B}(\varphi,...) \circ \varphi ^{-1} J \varphi ^{-1}  = b( u, \rho  , s) \\&\Theta _ \varrho (\varphi,...) \circ \varphi ^{-1} J \varphi ^{-1}= \theta _ \rho( u, \rho  , s) , \qquad \Theta _ S (\varphi,...) \circ \varphi ^{-1} J \varphi ^{-1} = \theta _ s( u, \rho  , s)
\end{aligned} 
\end{equation} 
and 
\begin{equation}\label{J_red}
\begin{aligned}
&\mathfrak{J}(\varphi,...) \circ \varphi_ \partial  ^{-1} J \varphi _ \partial ^{-1} = J(u, \rho  , s),\\&\mathfrak{j}_ \varrho (\varphi,...) \circ \varphi_ \partial  ^{-1} J \varphi _ \partial ^{-1} = j_ \rho (u, \rho  , s) , \qquad  \mathfrak{j}_ S (\varphi,...) \circ \varphi_ \partial  ^{-1} J \varphi _ \partial ^{-1} = j_ s(u, \rho  , s),
\end{aligned}
\end{equation} 
where $u= \dot  \varphi \circ \varphi  ^{-1} $, $ \rho  = (\varrho \circ \varphi ^{-1} )J \varphi ^{-1} $, and $ S  = (S \circ \varphi ^{-1} )J \varphi ^{-1} $. Above, we have used the abbreviation $(\varphi,...)= ( \varphi , \dot  \varphi , \varrho , S)$.

Note that the transformations appearing in \eqref{B_Theta}--\eqref{J_red} are the natural ones associated to the geometric type of these objects, see \eqref{B_def}--\eqref{j_def}.

\paragraph{Eulerian variational principle for open fluids.} Under the relabelling symmetries described above, we can state the following variational formulation of open fluids in the Eulerian frame:

\begin{proposition}[Eulerian variational principle for open fluids] \label{EP_open} Assume that the Lagrangian density and the distributed and boundary sources of momentum, mass and entropy satisfy the relabelling symmetries. Then the following hold:
\begin{itemize}
\item[\bf (i)]
The Lagrange-d'Alembert principle \eqref{extended_HP} yields the following variational formulation in the Eulerian frame (Euler-Poincar\'e-d'Alembert):
\begin{equation}\label{extended_reduced}
\begin{aligned} 
&\delta \int_0^T \!\!\int_ \Omega \Big[ \mathfrak{l}(u, \rho ,s ) + \rho D_t w + s D_t \gamma \Big]  {\rm d} x {\rm d} t \\
& \qquad \quad + \int_0^T \!\!\int_ \Omega \Big[b \cdot \zeta + \theta_ \rho   D_ \delta w  + \theta_s D_ \delta \gamma  \Big] {\rm d} x {\rm d} t\hspace{1.5cm}\leftarrow \text{bulk contribution}\\
& \qquad \quad  + \int_0^T \!\!\int_{ \partial \Omega } \Big[J \cdot \zeta + j _ \rho  D_ \delta w + j_ s D_ \delta \gamma \Big]   {\rm d} a {\rm d} t=0\hspace{0.5cm}\leftarrow \text{boundary contribution}
\end{aligned} 
\end{equation}
with respect to variations $\delta u=\partial _t \zeta+ [ \zeta, u]$ and free variations $ \delta \rho  $, $ \delta s$, $ \delta w$, $ \delta \gamma $, with $ \delta w$ and $ \delta \gamma $ vanishing at $t=0,T$ and $\zeta $ an arbitrary time dependent vector field vanishing at $t=0,T$.

The relation with the material variables used in \eqref{extended_HP} is given by
\[
u= \dot  \varphi \circ  \varphi ^{-1} , \quad \rho  = (\varrho \circ \varphi ^{-1} )J \varphi ^{-1} , \quad s  = (S \circ \varphi ^{-1} )J \varphi ^{-1} , \quad w= W \circ \varphi ^{-1} , \quad \gamma = \Gamma \circ \varphi ^{-1} .
\]

\item[\bf (ii)] This principle yields the following equations and boundary conditions for open fluids:
\begin{equation}\label{Open_fluid_Eulerian} 
\left\{
\begin{array}{l}
\vspace{0.2cm}\displaystyle\partial _t  \frac{\partial \mathfrak{l} }{\partial u} + \pounds _u \frac{\partial \mathfrak{l} }{\partial u} - \rho  \nabla  \frac{\partial \mathfrak{l} }{\partial \rho  }- s\nabla  \frac{\partial \mathfrak{l} }{\partial s  } = b\\
\vspace{0.2cm}\displaystyle\partial _t \rho  + \operatorname{div}( \rho  u) = \theta_ \rho  , \qquad \partial _t s  + \operatorname{div}( s  u) = \theta_s\\
\vspace{0.2cm}\displaystyle (u \cdot n) \frac{\partial \mathfrak{l} }{\partial u} =-J\quad\text{on}\quad \partial \Omega \\
\vspace{0.2cm}\displaystyle\rho   u   \cdot n=- j _ \rho , \qquad  s   u   \cdot n=- j _s\quad\text{on}\quad \partial \Omega .
\end{array}
\right. 
\end{equation}
\end{itemize}
\end{proposition}
\noindent\textbf{Proof.} {\bf (i)} By using the relations \eqref{B_red}, \eqref{J_red}, as well as $\delta W= D_ \delta w \circ \varphi $ and $\delta \Gamma = D_ \delta \gamma  \circ \varphi $, the bulk and boundary source terms in  \eqref{extended_HP} yield the ones in  \eqref{extended_reduced}.
From the relabelling symmetry of $\mathfrak{L}$, we have from \eqref{mathfrak_ell} 
\begin{align*} 
\left.\frac{d}{d\varepsilon}\right|_{\varepsilon=0}\int_{ \varphi _ \varepsilon^{-1}  ( \Omega )} 
\mathfrak{L}( \varphi _ \varepsilon , \dot  \varphi _ \varepsilon , \nabla_X \varphi _ \varepsilon, \varrho  _ \varepsilon , S_ \varepsilon  ){\rm d} X &= \left. \frac{d}{d\varepsilon}\right|_{\varepsilon=0}\int_{\Omega }  \mathfrak{L}( \varphi _ \varepsilon , \dot  \varphi _ \varepsilon, \nabla_X \varphi _ \varepsilon , \varrho  _ \varepsilon ,S_ \varepsilon) \circ \varphi^{-1} _ \varepsilon J \varphi^{-1} _ \varepsilon  {\rm d} x\\
&= \int_{\Omega } \left. \frac{d}{d\varepsilon}\right|_{\varepsilon=0} \mathfrak{L}( \varphi _ \varepsilon , \dot  \varphi _ \varepsilon, \nabla_X \varphi _ \varepsilon , \varrho  _ \varepsilon, S_ \varepsilon ) \circ \varphi^{-1} _ \varepsilon J \varphi^{-1} _ \varepsilon  {\rm d} x\\
&= \int_{\Omega }  \left. \frac{d}{d\varepsilon}\right|_{\varepsilon=0} \mathfrak{l} (u_ \varepsilon , \rho _ \varepsilon, s_ \varepsilon) {\rm d} x,
\end{align*} 
where $u_ \varepsilon = \dot  \varphi _ \varepsilon \circ \varphi ^{-1} $ and $ \rho  _ \varepsilon = \varrho _ \varepsilon \circ \varphi _ \varepsilon ^{-1} J \varphi _ \varepsilon ^{-1} $, $ s _ \varepsilon = S_ \varepsilon \circ \varphi _ \varepsilon ^{-1} J \varphi _ \varepsilon ^{-1} $. From these relations, we have
\[
\delta u = \partial _t \zeta + \pounds _u \zeta 
\]
while $\delta \rho  $ and $\delta s$ are free since $ \delta \varrho $ and $ \delta S $ are free. We also have, by using $w= W \circ \varphi ^{-1}$,
\[
\int_{ \varphi ^{-1} ( \Omega )} 
\varrho \dot W {\rm d} X = \int_{ \Omega }
\varrho \circ \varphi ^{-1} J \varphi ^{-1} \dot  W \circ \varphi ^{-1} {\rm d} x = \int_ \Omega \rho  D_ t w {\rm d} x.
\] 
Similarly we have  $\int_{ \varphi ^{-1} ( \Omega )} 
S \dot \Gamma  {\rm d} X=\int_ \Omega s  D_ t \gamma  {\rm d} x$. In conclusion, the principle \eqref{extended_HP} yields \eqref{extended_reduced}. 

{\bf (ii)} The treatment of $(s, \gamma )$ is similar to that of $( \rho  , w)$, so it is enough to present the computation for the variation of $\int_0^T \int_ \Omega \left[  \mathfrak{l} (u, \rho) + \rho D_tw \right]   {\rm d} x {\rm d} t $. The variation with respect to $ \rho  $ gives the condition 
\begin{equation}\label{w_condition} 
D_tw = - \frac{\partial \mathfrak{l} }{\partial \rho  }.
\end{equation} 
The variation of $\int_0^T \int_ \Omega \left[  \mathfrak{l} (u, \rho) + \rho D_tw \right]   {\rm d} x {\rm d} t $ with respect to $u$ and $w$ is computed as follows:
\begin{align*}
&= \int_0^T \int_ \Omega  \left[ \frac{\partial \mathfrak{l} }{\partial u} \cdot \left( \partial _t \zeta + \pounds_ u \zeta \right) + \rho  D_t \delta w + \left( \partial _t \zeta + \pounds_ u \zeta \right) \cdot \rho  \nabla w \right]{\rm d} x {\rm d} t \\
&= \int_0^T \int_ \Omega  \Bigg[ \left(- \partial _t \frac{\partial \mathfrak{l} }{\partial u} - \pounds_ u \frac{\partial \mathfrak{l} }{\partial u} \right) \cdot \zeta + \operatorname{div} \left( \Big( \frac{\partial \mathfrak{l} }{\partial u} \cdot \zeta \Big)  u \right) +\bar D_t( \rho   \delta w) - \bar D_t \rho   \delta w + \bar D_t( \rho  \zeta \cdot \nabla w)  \\
& \qquad \qquad \qquad  - \big( \partial _t ( \rho  \nabla w) + \pounds _u( \rho  \nabla w)\big) \cdot \zeta  \Bigg]{\rm d} x {\rm d} t  \\
&= \int_0^T \int_ \Omega  \Bigg[ \left(- \partial _t \frac{\partial \mathfrak{l} }{\partial u} - \pounds_ u \frac{\partial \mathfrak{l} }{\partial u} \right) \cdot \zeta + \operatorname{div} \left( \Big( \frac{\partial \mathfrak{l} }{\partial u} \cdot \zeta + \rho  D_ \delta w \Big)  u \right) -\bar D_t \rho   D_ \delta w -  \rho   \zeta \cdot  \nabla D_t w  \Bigg]{\rm d} x {\rm d} t  \\
&=  \int_0^T \int_ \Omega  \Bigg[ \left(- \partial _t \frac{\partial \mathfrak{l} }{\partial u} - \pounds_ u \frac{\partial \mathfrak{l} }{\partial u}+ \rho  \nabla \frac{\partial \mathfrak{l} }{\partial \rho  }  \right) \cdot \zeta  -\bar D_t \rho   D_ \delta w \Bigg]{\rm d} x {\rm d} t \\
& \qquad \qquad \qquad +  \int_0^T \int_ {\partial \Omega} \left[ \left(\frac{\partial \mathfrak{l} }{\partial u} \cdot \zeta \right) (u \cdot n) + \rho  u \cdot n D_ \delta w\right] {\rm d} a {\rm d} t,
\end{align*}
where we used \eqref{w_condition}, as well as $ \zeta $ and $ \delta w$ vanishing at $t=0,T$. Hence, when \eqref{extended_reduced} is considered, and by collecting the terms proportional to $ \zeta $, $ D_ \delta w$, and $D_ \delta \gamma $, both at the interior and at the boundary, and using the equalities $D_tw = - \frac{\partial \mathfrak{l} }{\partial \rho  } $, $D_t \gamma  = - \frac{\partial \mathfrak{l} }{\partial s } $ arising from the variations $ \delta \rho  $ and $\delta s $, we get the six equations \eqref{Open_fluid_Eulerian}. $\qquad\blacksquare$

\begin{remark}[Lagrangian reduction by symmetry]\label{remark_L_red}\rm Let us explain how the variational setting employed above fits in the Lagrangian reduction framework. In general, given a configuration manifold $Q$ and a Lagrangian $L:TQ \rightarrow \mathbb{R} $, we assume that the Lagrangian is invariant under the (tangent lifted) action of a Lie group action of $G$ on $Q$. In the case when the quotient $(TQ)/G$ is a smooth manifold, $L$ induces a reduced Lagrangian $\ell:(TQ)/G \rightarrow \mathbb{R} $ and the Hamilton principle for the action functional of $L$ induces a variational principle for the action functional of $\ell$.
The so called reduced principle is in general more involved than the Hamilton principle since it uses constrained variations, namely, it considers only the variations of the curves in $(TQ)/G$ that are induced by the free variations of the curves in $Q$. The same idea extends easily to the case with forcing, i.e., the Lagrange-d'Alembert principle, when the force is equivariant with respect to $G$. For the case of open fluid, the configuration manifold is originally $ \operatorname{Diff}( \mathbb{R} ^n  ) \times \Omega  ^n ( \mathbb{R} ^n ) \times \Omega ^n ( \mathbb{R} ^n)\ni ( \varphi , \varrho , S)$, but we extend it to
\[
Q= \operatorname{Diff}(\mathbb{R} ^n ) \times \Omega  ^n ( \mathbb{R} ^n) \times \Omega  ^n ( \mathbb{R} ^n ) \times \Omega  ^0 ( \mathbb{R} ^n ) \times  \Omega  ^0 ( \mathbb{R} ^n )\ni q=( \varphi , \varrho , S, W, \Gamma )
\]
in order to include the thermodynamic displacements $W$ and $ \Gamma $.
The Lagrangian $L(q, \dot  q)$ is 
\[
L(q, \dot  q)= \int_{ \varphi ^{-1} ( \Omega )}\left[  \mathfrak{L}( \varphi   , \dot  \varphi   , \nabla \varphi , \varrho   ) + \varrho   \dot  W  + S  \dot  \Gamma   \right]  {\rm d} X {\rm d} t
\]
and the force $F(q, \dot  q)$ is given by
\[
\left\langle F(q, \dot  q), \delta q \right\rangle =
\int_{ \varphi ^{-1} ( \Omega)} \Big[ \mathfrak{B} \cdot \delta \varphi +\Theta_ \varrho   \delta W +\Theta_ S   \delta \Gamma  \Big]
{\rm d} X + \int_{ \varphi ^{-1} (\partial\Omega)} \Big[ \mathfrak{J}  \cdot \delta \varphi +\mathfrak{j}_ \varrho  \delta W  +\mathfrak{j}_S  \delta \Gamma   \Big] {\rm d} A .
\]
This identification explicitly shows how \eqref{extended_HP} emerges from \eqref{LdA}. Now the symmetry Lie group is $G= \operatorname{Diff}( \mathbb{R} ^n )$ and it acts on $Q$ by the right action
\[
q= ( \varphi , \varrho , S, W, \Gamma ) \mapsto q \cdot \psi  = ( \varphi \circ  \psi  , (\varrho  \circ \psi ) J \psi , (S \circ \psi ) J \psi, W \circ  \psi , \Gamma  \circ \psi ),
\]
thereby leading to the Eulerian variables we have used above.
\end{remark}

\paragraph{Balance laws in the Eulerian description.} From the local balance of mass and entropy and the last two boundary conditions in \eqref{Open_fluid_Eulerian}, the corresponding global balance equations are
\begin{equation}\label{balance_mass&entropy}
\frac{d}{dt} \int_ \Omega \rho  {\rm d} x = \int_ \Omega \theta _ \rho  {\rm d} x+ \int_{ \partial \Omega }j _ \rho  {\rm d} a \qquad\text{and}\qquad   \frac{d}{dt} \int_ \Omega s  {\rm d} x = \int_ \Omega \theta_s {\rm d} x  +\int_{ \partial \Omega }j_ s {\rm d} a,
\end{equation}
showing the transfer of mass and entropy via bulk and boundary contributions.

We now consider the total energy density, defined for a general Lagrangian density $ \mathfrak{l} $ as 
\begin{equation}\label{e_density} 
\mathfrak{e}=\frac{\partial \mathfrak{l}}{\partial u}\cdot u  -\mathfrak{l}(u, \rho, s )
\end{equation} 
and we compute the local energy balance
\begin{equation}\label{loc_energy_Eulerian}
\begin{aligned} 
\partial _t \mathfrak{e}&= \partial_t \frac{\partial \mathfrak{l}}{\partial u} \cdot u   - \frac{\partial \mathfrak{l}}{\partial \rho} \cdot \partial_t \rho - \frac{\partial \mathfrak{l}}{\partial s} \cdot \partial_t s\\
&=- \pounds_ u \frac{\partial \mathfrak{l}}{\partial u}\cdot u  + \rho u \cdot \nabla \frac{\partial \mathfrak{l}}{\partial \rho} + s u \cdot \nabla \frac{\partial \mathfrak{l}}{\partial s} + b \cdot u + \frac{\partial \mathfrak{l}}{\partial \rho}  \operatorname{div} (\rho u) - \theta_{\rho} \frac{\partial \mathfrak{l}}{\partial \rho} + \frac{\partial \mathfrak{l}}{\partial s}  \operatorname{div} (s u) - \theta_s \frac{\partial \mathfrak{l}}{\partial s}\\
&= \operatorname{div}  \left(- \Big(\frac{\partial \mathfrak{l} }{\partial u} \cdot u \Big)u  + \rho \frac{\partial \mathfrak{l}}{\partial \rho} u + s \frac{\partial \mathfrak{l}}{\partial s} u \right) + b \cdot u  - \theta_{\rho} \frac{\partial \mathfrak{l}}{\partial \rho} - \theta_s \frac{\partial \mathfrak{l}}{\partial s}.
\end{aligned}
\end{equation} 
From this, using all three boundary conditions in \eqref{Open_fluid_Eulerian}, the global energy balance follows as
\begin{equation}\label{glob_energy_Eulerian} 
\frac{d}{dt} \int_ \Omega \mathfrak{e}\, {\rm d} x =\int_ \Omega \left[  b \cdot u  - \theta_{\rho} \frac{\partial \mathfrak{l}}{\partial \rho} - \theta_s \frac{\partial \mathfrak{l}}{\partial s} \right]  {\rm d} x + \int_{\partial \Omega} \left[ J\cdot u - j_ \rho   \frac{\partial \mathfrak{l}}{\partial \rho  }- j_ s   \frac{\partial \mathfrak{l}}{\partial s  }  \right] {\rm d} a,
\end{equation} 
showing the energy transfer between the fluid and its exterior via both bulk and boundary contributions.
We note that from the boundary conditions in \eqref{Open_fluid_Eulerian} we have the following relations at the boundary
\[
j_s= \frac{ s }{\rho} j_ \rho  \quad\text{and}\quad  J= j_ \rho\frac{1}{ \rho  } \frac{\partial \mathfrak{l} }{\partial u}=j_ s\frac{1}{ s  } \frac{\partial \mathfrak{l} }{\partial u} \quad\text{on}\quad  \partial \Omega,
\]
hence the boundary contribution to the energy balance in \eqref{glob_energy_Eulerian} can be written in either of the following forms
\begin{equation}\label{forms_boundary_flow} \begin{aligned} 
\int_{\partial \Omega} \left[ j_ \rho  \left( \frac{1}{ \rho  } \frac{\partial \mathfrak{l} }{\partial u} \cdot u -  \frac{\partial \mathfrak{l}}{\partial \rho  } \right) - j_ s   \frac{\partial \mathfrak{l}}{\partial s  }  \right] {\rm d} a&= \int_{\partial \Omega}  j_ \rho  \left( \frac{1}{ \rho  } \frac{\partial \mathfrak{l} }{\partial u} \cdot u -  \frac{\partial \mathfrak{l}}{\partial \rho  } - \frac{s}{ \rho  }   \frac{\partial \mathfrak{l}}{\partial s  }  \right) {\rm d} a\\
&= \int_{\partial \Omega}  \left(  \rho\frac{\partial \mathfrak{l}}{\partial \rho  }   +  s\frac{\partial \mathfrak{l}}{\partial s  }  - u \cdot \frac{\partial \mathfrak{l} }{\partial u}  \right)( u \cdot n ){\rm d} a.
\end{aligned}
\end{equation}

\paragraph{The compressible Euler fluid with open boundaries.} Let us consider the Lagrangian density of compressible Euler fluid as given in \eqref{Euler_mathfrakL}. One can readily checks that it satisfies the material covariance assumption \eqref{mat_cov}. In the Eulerian description the Lagrangian density takes the standard form 
\begin{equation}\label{l_Euler}
\mathfrak{l}(u, \rho ,s )= \frac{1}{2} \rho  | u | ^2 - \varepsilon ( \rho  , s),
\end{equation}
as it is seen from \eqref{mathfrak_ell}. In this case the system \eqref{Open_fluid_Eulerian} becomes
\begin{equation}\label{Euler_rho_u} 
\left\{
\begin{array}{l}
\vspace{0.2cm}\displaystyle\partial _t  (\rho u) + \operatorname{div}( \rho   u \otimes u)  = - \nabla p + b\\
\vspace{0.2cm}\displaystyle\partial _t \rho  + \operatorname{div}( \rho  u) = \theta_ \rho  , \qquad \partial _t s  + \operatorname{div}( s  u) = \theta_s\\
\vspace{0.2cm}\displaystyle (u \cdot n) \rho u =-J\quad\text{on}\quad \partial \Omega \\
\vspace{0.2cm}\displaystyle\rho   u   \cdot n=- j _ \rho , \qquad  s   u   \cdot n=- j _s\quad\text{on}\quad \partial \Omega .
\end{array}
\right. 
\end{equation} 
We note that when written in terms of the velocity, the momentum equation acquires an additional term due to the presence of the distributed source of mass $ \theta _ \rho  $, namely,
\begin{equation}\label{Euler_u} 
\rho  ( \partial _t u + u \cdot \nabla u)=- \nabla p +b- \theta _ \rho  u.
\end{equation}



By using  $\frac{\partial \mathfrak{l} }{\partial \rho  } =\frac{1}{2} |u| ^2 -g$ and $\frac{\partial \mathfrak{l} }{\partial s} =  -T$ with $g=\frac{\partial \varepsilon }{\partial \epsilon } = (\varepsilon+ p  - sT)/ \rho$ the Gibbs potential and $T=\frac{\partial \varepsilon }{\partial s}$ the temperature, the local energy balance \eqref{loc_energy_Eulerian} becomes
\[
\partial _t \mathfrak{e}+ \operatorname{div}(( \mathfrak{e}+p)u) = b \cdot u + \theta _ \rho  \Big( g - \frac{1}{2} |u| ^2 \Big) + \theta _s T,
\]
giving the total energy balance \eqref{glob_energy_Eulerian} in the form.
\begin{equation}\label{balance_energy} 
\frac{d}{dt} \int_ \Omega \mathfrak{e} \, {\rm d} x=\int_ \Omega \Big[  b \cdot u  +\theta_{\rho}\Big( g - \frac{1}{2} |u| ^2 \Big) + \theta_s T\Big]  {\rm d} x + \int_{\partial \Omega} \Big[ J\cdot u + j_ \rho   \Big( g - \frac{1}{2} |u| ^2 \Big) + j_ s   T  \Big] {\rm d} a.
\end{equation} 
Following \eqref{forms_boundary_flow}, other useful expressions for the boundary contribution to energy change are given by
\[
\int_{\partial \Omega} \Big[  j_ \rho   \Big( \frac{1}{2} |u| ^2  + g \Big) + j_ s   T  \Big] {\rm d} a =\int_{\partial \Omega}   j_ \rho   \Big( \frac{1}{2} |u| ^2  + h \Big) {\rm d} a=-\int_{\partial \Omega}  ( \mathfrak{e}+p) u \cdot n {\rm d} a,
\]
with $h=( \varepsilon +p)/ \rho  $ the enthalpy. One of these expressions may be useful depending on whether it is desirable to have the boundary sources
$j_\rho$ and $j_s$ explicitly appear or not.

Apart from the overall energy balance, in the case of Euler fluid, we can calculate the balance law for the kinetic energy density $ \mathfrak{k} = \frac{1}{2} \rho  |u| ^2 $ and the internal energy density $ \varepsilon $. Using the equations, we get
\[
\partial _t \mathfrak{k} + \operatorname{div}( \mathfrak{k} u) = - u \cdot \nabla p + b \cdot u - \theta _ \rho \frac{1}{2} |u| ^2   
\]
\[
\partial _t \varepsilon + \operatorname{div}( \varepsilon u) = - p \operatorname{div}u+  \theta _ \rho   g + \theta _s T.
\]
Then, from the boundary conditions the global balance equations become
\begin{align*} 
\frac{d}{dt} \int_ \Omega \mathfrak{k}\, {\rm d} x& = \int_{\Omega} \Big[ - u \cdot \nabla p + b \cdot u - \theta _ \rho \frac{1}{2} |u| ^2\Big] {\rm d}  x + \int_ {\partial \Omega} j_ \rho  \frac{1}{2} | u| ^2 {\rm d} a\\
\frac{d}{dt} \int_ \Omega \varepsilon \,{\rm d} x &= \int_ \Omega \Big[  u \cdot \nabla p +  \theta _ \rho   g + \theta _s T \Big] {\rm d} x+ \int_{\partial \Omega }  \Big[ j_ \rho   g + j_ s   T  \Big] {\rm d} a\\
&= \int_ \Omega \Big[  u \cdot \nabla p +  \theta _ \rho   g + \theta _s T \Big] {\rm d} x+ \int_{\partial \Omega }   j_ \rho   h {\rm d} a.
\end{align*}

Adding these two equations up, we recover exactly the total energy balance law \eqref{balance_energy}, however these detailed balance equations show that in the interior of the fluid domain, the source term $b$ contributes uniquely to the kinetic energy and the source term $\theta_s$ completely to the internal energy, while the source $\theta_{\rho}$ splits in two parts $-\theta_{\rho}\frac{1}{2}|u|^2 $ for the kinetic energy contribution and $\theta_{\rho}g$ for the internal energy contribution. A similar situation occurs on the boundary. The term $ \int_{\Omega} u \cdot \nabla p \diff x$ determines the exchange rate of energy between the kinetic and the internal energy. 


\paragraph{Effect of rotation and gravity.} In the case of a rotating fluid subject to gravity, the Lagrangian density is
\[
\mathfrak{l} (u, \rho  , s)= \frac{1}{2} \rho  |u| ^2 + \rho  R \cdot u - \varepsilon ( \rho  , s) - \rho  \phi.
\]
Here  $R(x)$ is the vector potential for the Coriolis parameter, i.e., $ \operatorname{curl}R= 2 \omega$ with $ \omega $ the angular velocity of the Earth, and $ \phi (x)$ is the gravitational potential. In this case, system \eqref{Open_fluid_Eulerian} becomes
\[
\left\{
\begin{array}{l}
\vspace{0.2cm}\displaystyle\partial _t  (\rho u) + \operatorname{div}( \rho   u \otimes u)  +   2 \omega \times \rho  u = - \nabla p - \rho  \nabla \phi  + b -  \theta _ \rho  R\\
\vspace{0.2cm}\displaystyle\partial _t \rho  + \operatorname{div}( \rho  u) = \theta_ \rho  , \qquad \partial _t s  + \operatorname{div}( s  u) = \theta_s\\
\vspace{0.2cm}\displaystyle (u \cdot n) \rho (u+R) =-J\quad\text{on}\quad \partial \Omega \\
\vspace{0.2cm}\displaystyle\rho   u   \cdot n=- j _ \rho , \qquad  s   u   \cdot n=- j _s\quad\text{on}\quad \partial \Omega .
\end{array}
\right. 
\]

We note the appearance of the vector potential $R$ in the fluid momentum equation, due to the bulk source of mass $ \theta _ \rho  $ and in the boundary condition for the fluid momentum, compare with the non rotating case \eqref{Euler_rho_u} . The fluid momentum equation can be rewritten in a similar way with \eqref{Euler_u} as
\begin{equation}\label{rot_Euler_u} 
\rho  ( \partial _t u + u \cdot \nabla u + 2 \omega \times u)=- \nabla p- \rho  \nabla \phi  +b- \theta _ \rho (u+R),
\end{equation} 
which also show the occurrence of $R$.

We note that the total energy density as defined in \eqref{e_density} becomes
\[
\mathfrak{e}= \frac{1}{2} \rho  |u| ^2 + \varepsilon ( \rho  , s) + \rho   \phi,
\]
in which the contribution of rotation as been eliminated.
In order to further show the impact of rotation and gravity we now list all the global and local balance equations for kinetic, internal, potential, and total energy densities:
\begin{align*} 
&\partial _t \mathfrak{k} + \operatorname{div}( \mathfrak{k} u) = - u \cdot \nabla p - \rho  u \cdot \nabla \phi + b \cdot u - \theta _ \rho \Big( \frac{1}{2} |u| ^2+ u \cdot R \Big) \\
&\partial _t \varepsilon + \operatorname{div}( \varepsilon u) = - p \operatorname{div}u+  \theta _ \rho   g + \theta _s T\phantom{\Big(}\\
&\partial _t ( \rho  \phi ) + \operatorname{div}( \rho  \phi  u) = \theta _ \rho  \phi + \rho  u \cdot \nabla  \phi\phantom{\Big(} \\
& \partial _t \mathfrak{e} + \operatorname{div}(( \mathfrak{e}+p)u) = b \cdot u + \theta _ \rho  \Big( g + \phi  - u \cdot R- \frac{1}{2} |u| ^2 \Big) + \theta _s T,
\end{align*} 
and
\begin{align*}
&\frac{d}{dt} \int_\Omega \mathfrak{k} \, {\rm d}  = \int_{\Omega} \Big[ - u \cdot \nabla p - \rho  u \cdot \nabla \phi  + b \cdot u - \theta _ \rho \Big(\frac{1}{2} |u| ^2 + u \cdot R\Big)\Big] {\rm d}  x + \int_ {\partial \Omega} j_ \rho  \frac{1}{2} | u| ^2 {\rm d} a\\
&\frac{d}{dt} \int_ \Omega \varepsilon \,{\rm d} x = \int_ \Omega \Big[  u \cdot \nabla p +  \theta _ \rho   g + \theta _s T \Big] {\rm d} x+ \int_{\partial \Omega }   j_ \rho   h {\rm d} a\\
&\frac{d}{dt} \int_ \Omega \rho  \phi  \,{\rm d} x = \int_ \Omega \Big[ \theta _ \rho  \phi + \rho  u \cdot \nabla \phi  \Big] {\rm d} x +\int_{ \partial \Omega }  j_ \rho \phi  {\rm d} a\\
&\frac{d}{dt} \int_ \Omega\mathfrak{e}   \,{\rm d} x = \int_ \Omega \Big[ b \cdot u - \theta _ \rho  \Big( \frac{1}{2} |u| ^2 + u \cdot R -g - \phi \Big)+ \theta _sT\Big] {\rm d} x +\int_{ \partial \Omega }  j_ \rho\Big( \frac{1}{2} |u| ^2 + h+ \phi \Big)  {\rm d} a.
\end{align*} 

\begin{remark}[Absolute Momentum vs Mass Flux $\rho u$]\rm A major insight of our variational approach is the clear distinction between the absolute fluid momentum $\frac{\partial \mathfrak{l}}{\partial u}$, and the mass flux $\rho u$. We see that both of them appear in the expression of the boundary conditions for the general open fluid equations in \eqref{Open_fluid_Eulerian}, namely:
\[
(u \cdot n) \frac{\partial\mathfrak{l}}{\partial u}=-J \quad\text{and}\quad \rho u \cdot n = - j_ \rho,
\]
(instead of $\frac{\partial\mathfrak{l}}{\partial u}\cdot n=-j_\rho$ for instance for the second one).
These quantities coincide if $\mathfrak{l}$ depends on $u$ only through the usual expression $\frac{1}{2}\rho|u|^2$ which is the case for nonrotating Euler, amongst other equations, but is not true in the general case. When they do coincide we have the relation $j_\rho u=J$ on the boundary, which is not true in general, for example the rotating case giving $j_\rho(u+R)=J$.
\end{remark}

\begin{remark}[Choice of the rotational vector potential]\rm For open fluids, the choice of the rotational vector potential $R$ such that $\operatorname{curl}R=2\omega$ impacts both the choice of the distributed momentum source $b$ and mass source $j_\rho$, since it changes the relationship between velocity $u$ and absolute momentum $\frac{\partial \mathfrak{l}}{\partial u}$ by a factor of $\rho R$. Therefore, when using $R'= R+ \nabla f$ for some function $f$, one needs to accordingly make the changes
\[
b\rightarrow b'=b-\theta_\rho \nabla f\quad\text{and}\quad J \rightarrow J'=J-j_\rho \nabla f
\]
in the distributed and boundary source of fluid momentum, to have the physics unchanged.
\end{remark}

\color{black}
\paragraph{The rotating shallow water equations.} The Lagrangian density for the rotating shallow water equations is
\[
\mathfrak{l}(u, h)= \frac{1}{2} h|u| ^2 + h R \cdot u- \frac{1}{2} g (h+Z) ^2 
\]
with $Z(x)$ the bottom topography and $h$ the water depth, so that the surface level is $h+Z$.
The associated open fluid equations found from \eqref{Open_fluid_Eulerian}  are
\[
\left\{
\begin{array}{l}
\vspace{0.2cm}\displaystyle\partial _t  (h u) + \operatorname{div}(hu \otimes u)   + 2 \omega \times hu = - gh \nabla (h+Z)  + b - \theta _h R\\
\vspace{0.2cm}\displaystyle\partial _t h  + \operatorname{div}( h  u) = \theta_ h \\
\vspace{0.2cm}\displaystyle h(u \cdot n)  (u+R) =-J\quad\text{on}\quad \partial \Omega \\
\vspace{0.2cm}\displaystyle h   u   \cdot n=- j _ h \quad\text{on}\quad \partial \Omega .
\end{array}
\right. 
\]
The fluid momentum equation can be written as
\[
h( \partial _t u + u \cdot \nabla u+ 2 \omega  \times u)= - gh \nabla (h+Z) - \theta _h (u+R) +b,
\]
showing again the role of the source term $ \theta _h$. Distributed source terms 
$\theta_h$ in the continuity equation of shallow water models appear, for instance, in the modeling of moist shallow water equations, as discussed in \cite{bouchut2009fronts}.

The total energy density as defined in \eqref{e_density} is now
\[
\mathfrak{e}= \frac{1}{2} h  |u| ^2 +  \frac{1}{2} g(h+Z) ^2 ,
\]
in which the contribution of rotation as been eliminated. The list of global and local balance equations for kinetic, potential $ \mathfrak{p} = \frac{1}{2} g(h+Z) ^2$, and total energy densities reads as follows:
\begin{align*} 
&\partial _t \mathfrak{k} + \operatorname{div}( \mathfrak{k} u) = - gh u \cdot \nabla (h+Z) + b \cdot u - \theta _ h \Big( \frac{1}{2} |u| ^2+ u \cdot R \Big) \\
&\partial _t \mathfrak{p}  + \operatorname{div}(\mathfrak{p}u)  =\operatorname{div} \Big(\frac{1}{2} g( Z ^2 - h ^2 )u\Big) + ghu \cdot \nabla (h+Z) + g(h+Z) \theta _h\\
& \partial _t \mathfrak{e} + \operatorname{div}\Big( \mathfrak{e}u+ \frac{1}{2} g( h ^2 - Z ^2 )u\Big) = b \cdot u + \theta _ h  \Big( g (h+Z)  - u \cdot R- \frac{1}{2} |u| ^2 \Big),
\end{align*} 
and
\begin{align*}
&\frac{d}{dt} \int_ \Omega \mathfrak{k} \, {\rm d} x = \int_ \Omega \Big[ - ghu \cdot \nabla (h+Z) + b \cdot u- \theta _h \Big( \frac{1}{2} |u| ^2 +R \cdot u \Big) \Big] {\rm d} x+\int_{ \partial \Omega } j_h \frac{1}{2} |u| ^2 {\rm d} a\\
&\frac{d}{dt} \int_ \Omega \mathfrak{p} \, {\rm d} x=\int_ \Omega \Big[  ghu \cdot \nabla (h+Z) + \theta _h  g(h+Z)\Big] {\rm d} x+\int_{ \partial \Omega } j_h g(h+Z){\rm d} a\\
&\frac{d}{dt} \int_ \Omega \mathfrak{e} \, {\rm d} x =\int_ \Omega \Big[b \cdot u - \theta _h \Big( \frac{1}{2} |u| ^2 +R \cdot u -  g(h+Z)\Big)\Big] {\rm d} x+\int_{ \partial \Omega } j_h \Big(\frac{1}{2} |u| ^2 + g(h+Z)\Big){\rm d} a.
\end{align*}

\paragraph{On the choice of the fluxes and inflow/outflow boundary conditions.} We now focus on the non rotating case and analyse how our setting allows the treatment of inflow and outflow boundary conditions.
For this discussion we consider the Lagrangian \eqref{l_Euler} and we assume $\rho>0$ everywhere. The boundary conditions take the following form:
\[
(u \cdot n) \rho  u = - J, \quad \rho  u \cdot n=- j_ \rho  , \quad su \cdot n = - j_s.
\]

(i) If we first assume $j_\rho=0$ then we have $u\cdot n=0$ on $\partial\Omega$ which corresponds to an impermeable boundary. In this case, the relations above impose the choice $j_s=0$ and $J=0$ and there are no boundary conditions for $\rho|_{\partial\Omega}$ and $s|_{\partial\Omega}$, consistently with the treatment of fluids with impermeable boundaries. It is interesting to note that, while our geometric setting is primarily designed for open fluid motion, it allows the treatment of closed fluids by appropriate choice of $J,j_\rho,j_s$.

(ii) If $j_\rho\neq 0$, then $u\cdot n\neq 0$ and $J$ must satisfy $J\cdot n\neq 0$. We get the following boundary \textit{relations}
\[
u|_{\partial\Omega}=\frac{J}{j_\rho}, \quad \rho|_{\partial\Omega}= - \frac{j_\rho^2}{J\cdot n}, \quad s|_{\partial\Omega}= - \frac{j_\rho j_s}{J\cdot n}.
\]
In particular, we note the relations $j_s=\frac{\rho}{s}j_\rho$ and $J\cdot n<0$ on the boundary.
Let us show that these relations do yield appropriate boundary \textit{conditions} for inflow and outflow.

(ii.a) The case $j_\rho>0$ corresponds to inflow: $u \cdot n<0$. In this situation $J$, $j_\rho$ and $j_s$ can be freely chosen, independently on the current fluid motion, so that $u|_{\partial\Omega}$, $\rho|_{\partial\Omega}$, $s|_{\partial\Omega}$ are all prescribed on the boundary, consistently with the treatment of inflow. Equivalently, $u|_{\partial\Omega}$, $\rho|_{\partial\Omega}$, $T|_{\partial\Omega}$ are all prescribed, with $T$ the temperature.

(ii.b) The case $j_\rho<0$ corresponds to outflow: $u \cdot n>0$. In this situation, we choose 
\[
j_\rho= -\rho u_0\cdot n\quad\text{and}  \quad J=- (u_0\cdot n)\rho u_0,
\]
where we recall that a dependence of $j_\rho$ on the current fluid variables (here $\rho$) is allowed in our approach, see \S\ref{EPdA}, and where $u_0$ is a prescribed velocity field on $\partial\Omega$ with $u_0\cdot n>0$. In this case the first two boundary relations just give one boundary condition
\[
u|_{\partial \Omega}=u_0,
\]
consistently with the treatment of outflow boundary conditions for which there is no prescription on $\rho|_{\partial\Omega}$. The temperature can be prescribed at the outflow boundary by choosing
\[
j_s= -s(\rho, T_0) u_0\cdot n,
\]
where the state equation is used to write $s$ in terms of $\rho$ and $T_0$.
With these choices, we see that the three boundary relations above yield the two boundary conditions $u|_{\partial\Omega}=u_0$ and $T|_{\partial \Omega}=T_0$, consistently with what is assumed for outflow in the viscous case, see, e.g., \cite{Se1959}. In the inviscid case, we choose
\[
j_\rho= -\rho \nu_0, \quad j_s= -s \nu_0,\quad\text{and}  \quad J=- \nu_0\rho u,
\]
for $\nu_0>0$. With these choices, we see that the three boundary relations above yield the single boundary condition
\[
u\cdot n|_{\partial\Omega}=\nu_0,
\]
consistently with what is assumed for outflow in the non-viscous case, see, e.g., \cite{Se1959}.

While the boundary conditions as detailed above are used for well-posedness questions, \cite{Se1959}, in practice the number of boundary conditions at the inflow and outflow depend on the (subsonic/supersonic) fluid regime. Appropriate choices of $J$, $j_\rho$, and $j_s$ can also handle these cases.

\color{black}

\section{Hamiltonian variational and bracket formulations.}
\label{sec_H}

In this section we present the Hamiltonian counterpart of the approach presented in \S\ref{sec_L}. This approach is useful for deriving the appropriate extension of the Lie-Poisson bracket in fluid dynamics that can accommodate open boundary conditions. As in the previous section, we start with the material description as it is in this description that the variational principle is simpler. The Eulerian version is then derived through reduction. We do not include detailed computations here, as they can be carried out in parallel with those on the Lagrangian side.

We are using the word Hamiltonian in this section, since the setting is based on using Hamiltonian functions and densities, rather than Lagrangian ones. However, the systems under consideration, being open, are \textit{not Hamiltonian systems} and, consequently, we do not aim to find a Poisson bracket for this system, but to describe a natural extension of the Lie-Poisson bracket that accounts for the bulk and boundary fluxes.

\subsection{Material formulation}\label{Mat_H}

Start by assuming the existing of a Hamiltonian function $H : T^* \mathrm{Diff} (\mathbb{R}^n) \times \Omega  ^n ( \mathbb{R} ^n ) \times \Omega  ^n ( \mathbb{R} ^n )\rightarrow \mathbb{R}$ and Hamiltonian density $\mathfrak{H}(\varphi, M, \nabla \varphi, \varrho, S)$, related by
\begin{equation}\label{typical_H} 
H(\varphi, M, \varrho, S) := \int_{\varphi^{-1}(\Omega)} \mathfrak{H}(\varphi, M, \nabla \varphi, \varrho, S) \diff X.
\end{equation} 

The next result is the Hamiltonian analog to \eqref{extended_HP} and  Proposition \ref{critical_condition}.

\begin{proposition}
The Hamilton-d'Alembert phase-space principle 
\begin{equation}\label{extended_PS}
\begin{aligned} 
& \left. \frac{d}{d\varepsilon}\right|_{\varepsilon=0}  \int_0^T \!\!\!\int_{ \varphi_ \varepsilon  ^{-1} ( \Omega )}\left[ M_ \varepsilon \cdot \dot  \varphi _ \varepsilon+ \varrho _ \varepsilon\dot  W_ \varepsilon + S_ \varepsilon \dot \Gamma _ \varepsilon- \mathfrak{H}( \varphi _ \varepsilon, M_ \varepsilon, \nabla \varphi _ \varepsilon, \varrho _ \varepsilon, S_ \varepsilon) \right]  {\rm d} X {\rm d} t\\
& \qquad + \int_0^T \!\!\! \int_{ \varphi ^{-1} ( \Omega)} \!\!\!\big( \mathfrak{B} \cdot \delta \varphi +\Theta_ \varrho   \delta W +\Theta_ S   \delta \Gamma  \big) 
{\rm d} X {\rm d} t\;\;\;\;\;\;\;\leftarrow \text{bulk contribution}\\
& \qquad + \int_0^T \!\!\!\int_{ \partial (\varphi ^{-1} (\Omega)) } \!\!\!\big( \mathfrak{J}  \cdot \delta \varphi +\mathfrak{j}_ \varrho  \delta W  +\mathfrak{j}_S  \delta \Gamma   \big) {\rm d} A{\rm d} t=0\;\;\leftarrow \text{boundary contribution}
\end{aligned}
\end{equation} 
with free variations $\delta M$, $\delta\varphi$, $\delta\varrho$, $\delta W$, $\delta S$, $\delta \Gamma$, such that $\delta\varphi$, $\delta W$, and $\delta \Gamma $ vanish at $t=0,T$, gives the following equations: 
\begin{equation}\label{PS_force} 
\left\{
\begin{array}{l}
\vspace{0.2cm}\displaystyle \frac{d}{dt} \varphi = \frac{\partial \mathfrak{H} }{\partial M}, \qquad  \frac{d}{dt}  M  -\operatorname{DIV} \frac{\partial \mathfrak{H}}{\partial \nabla \varphi } + \frac{\partial \mathfrak{H}}{\partial \varphi } = \mathfrak{B} \\
\vspace{0.2cm}\displaystyle \frac{d}{dt} \varrho = \Theta_ \varrho , \qquad \frac{d}{dt}  S = \Theta_ S\\
\vspace{0.2cm}\displaystyle N \cdot \Big( -\frac{\partial \mathfrak{H}}{\partial \nabla \varphi }   + (\nabla _x \varphi ^{-1} \circ \varphi )\cdot \frac{\partial \mathfrak{H}}{\partial M}\, M + \Big(\mathfrak{H} - M \cdot \frac{\partial \mathfrak{H}}{\partial M} -\varrho \frac{\partial \mathfrak{H}}{\partial \varrho} - S \frac{\partial \mathfrak{H}}{\partial S} \Big)  \nabla _x \varphi ^{-1} \circ \varphi \Big) =-\mathfrak{J}\quad\text{on}\quad \partial \Omega \\
\displaystyle N \cdot \Big(( \nabla _x \varphi ^{-1} \circ \varphi )\cdot \frac{\partial \mathfrak{H}}{\partial M} \varrho \Big)=- \mathfrak{j}_ \varrho  , \quad N \cdot \Big((\nabla _x \varphi ^{-1} \circ \varphi) \cdot \frac{\partial \mathfrak{H}}{\partial M} S \Big)=- \mathfrak{j}_ S\quad\text{on}\quad \partial \Omega ,
\end{array}
\right. 
\end{equation} 
\end{proposition}

\begin{remark}[Hamilton-d'Alembert phase space principle]\rm The variational principle \eqref{extended_PS} is an instance of the Hamilton-d'Alembert phase space principle, given in general as
\[
\delta \int_0^ T\Big[\left\langle p,\dot  q \right\rangle - H(q,p) \Big] {\rm d} t + \int_0^T \left\langle F(q, \dot  q), \delta q \right\rangle {\rm d} t=0,
\]
for a Hamiltonian $H: T^*Q \rightarrow \mathbb{R}$ and a force $F:TQ \rightarrow T^*Q$.
It is the Hamiltonian counterpart of the Lagrange-d'Alembert principle \eqref{LdA}.
\end{remark} 

\subsubsection{Relating Hamiltonian and Lagrangian formulations}
Given a Lagrangian density $\mathfrak{L}( \varphi , \dot  \varphi , \nabla\varphi  , S)$, the corresponding Hamiltonian density is obtained by the \textit{Legendre transform}, assuming that the map
\begin{equation}\label{LT} 
\dot  \varphi \mapsto M=\frac{\partial \mathfrak{L} }{\partial \dot  \varphi } ( \varphi , \dot  \varphi , \nabla \varphi, \varrho , S)
\end{equation} 
is a diffeomorphism for all $\varphi , \dot  \varphi , \varrho , S$ with $ \varrho >0$\footnote{In other words, the Lagrangian is \textit{non-degenerate} and therefore the \textit{Legendre transform} is invertible.}. In this case we can define the Hamiltonian density $\mathfrak{H}$ in the following way: 
\[
\mathfrak{H}(\varphi, M, \nabla \varphi, \varrho, S) := M \cdot \dot \varphi - \mathfrak{L}(\varphi, \dot \varphi, \nabla \varphi, \varrho, S),
\]
where $\dot \varphi$ in the formula above is the inverse image of $M$  by \eqref{LT}. Alternatively, one could start with a Hamiltonian and obtain a Lagrangian, assuming that the \textit{inverse Legendre transform} is a diffeomorphism. This is obviously the case for all the fluid models treated above, although there are fluid models where this is not the case and there is either a Lagrangian or a Hamiltonian formulation, but not both. For example, the hydrostatic primitive equations in Eulerian coordinates have only a Lagrangian representation \cite{dubos2014equations}. We note the relations 
\[
\frac{\partial \mathfrak{H}}{\partial M} = \dot \varphi, \quad \frac{\partial \mathfrak{H}}{\partial \cdot } = -\frac{\partial \mathfrak{L}}{\partial \cdot }, \; \text{for } \varphi, \nabla \varphi, \varrho, S.
\]

\paragraph{Example: the compressible Euler fluid with open boundaries.} By taking the Legendre transform of the Lagrangian density \eqref{Euler_mathfrakL} of the compressible Euler fluid, we get
\[
\mathfrak{H}(\varphi, M, \nabla \varphi, \varrho, S) = \frac{|M|^2}{2 \varrho } + \varepsilon \left( \frac{ \varrho }{J \varphi } , \frac{ S}{J \varphi }\right) J \varphi . 
\]
The partial derivatives are
\begin{equation}\label{partial_der_H} \begin{aligned}
&\frac{\partial \mathfrak{H}}{\partial \varphi }=0, \quad \frac{\partial \mathfrak{H}}{\partial M }= \frac{M}{\varrho } , \quad \frac{\partial \mathfrak{H}}{\partial \varrho  }=\frac{\partial \varepsilon }{\partial \rho } - \frac{| M | ^2}{2 \varrho^2} , \quad 
\frac{\partial \mathfrak{H}}{\partial S  }=  \frac{\partial \varepsilon }{\partial s } \\
&\frac{\partial \mathfrak{L}}{\partial \nabla \varphi }= \Big( \varepsilon - \frac{\partial \varepsilon }{\partial \rho } \frac{\varrho}{J \varphi} - \frac{\partial \varepsilon }{\partial s } \frac{S}{J \varphi} \Big) J \varphi \left( \nabla_x\varphi^{-1} \circ \varphi \right)=-(p \circ \varphi )J \varphi \left( \nabla_x\varphi^{-1} \circ \varphi \right)
\end{aligned}
\end{equation}
hence the system \eqref{PS_force} becomes 
\begin{equation}\label{PS_force_EULER} 
\left\{
\begin{array}{l}
\vspace{0.2cm}\displaystyle  \frac{d}{dt} \varphi = \frac{M}{ \varrho }, \qquad 
\frac{d}{dt} M+ (\nabla _x p \circ \varphi ) J \varphi  = \mathfrak{B} \\
\vspace{0.2cm}\displaystyle \frac{d}{dt}   \varrho = \Theta_ \varrho, \qquad \frac{d}{dt}  S = \Theta_ S \\
\vspace{0.2cm}\displaystyle N \cdot \Big[(\nabla_x \varphi^{-1} \circ \varphi) \cdot \frac{M}{\varrho}\Big] M = -\mathfrak{J} \quad\text{on}\quad \partial \Omega\\
\displaystyle N \cdot \Big[(\nabla_x \varphi ^{-1} \circ \varphi) \cdot \frac{M}{\varrho} \Big] \varrho =- \mathfrak{j}_ \varrho  , \quad N \cdot \Big[(\nabla_x \varphi ^{-1} \circ \varphi) \cdot \frac{M}{\varrho} \Big] S =- \mathfrak{j}_ S\quad\text{on}\quad \partial \Omega.
\end{array}
\right. 
\end{equation}
This system of equations and boundary conditions for the open compressible fluid in the material descriptions is readily seen to be equivalent to the one derived in the Lagrangian side in \eqref{EL_force_EULER} by writing $ M= \varrho \dot  \varphi $.

\subsection{Eulerian formulation}\label{Eul_H}

We now assume relabelling symmetries and deduce from it the variational principle in the Eulerian description, induced by the Hamilton-d'Alembert phase space principle \eqref{PS_force}.

\paragraph{Relabelling symmetries.} We say that a Hamiltonian density is materially covariant it it satisfies 
\begin{equation}\label{mat_cov_H} 
\mathfrak{H}\big( \varphi \circ \psi , (M \circ \psi )J \psi , \nabla_X ( \varphi \circ \psi ) ,(\varrho \circ \psi )J \psi,(S \circ \psi )J \psi \big)= \mathfrak{H}\big( \varphi , M, \nabla_X \varphi , \varrho , S\big) \circ \psi J \psi ,
\end{equation} 
for all $ \psi \in \operatorname{Diff}( \mathbb{R} ^n )$. It is clear that if a Lagrangian density $ \mathfrak{L} $ is materially covariant, see \eqref{mat_cov}, then the associated Hamiltonian density obtained via the Legendre transform is also.
The symmetry \eqref{mat_cov_H} is equivalent to 
the existence of a Hamiltonian density $\mathfrak{h}(m, \rho, s )$ in the Eulerian description, such that
\begin{equation}\label{mathfrak_h} 
\mathfrak{H}\big( \varphi ,M, \nabla_X \varphi , \varrho  ,S \big) \circ \varphi  ^{-1} J \varphi ^{-1} = \mathfrak{h} \big(m, \rho  , s \big).
\end{equation} 
where $m = (M\circ \varphi ^{-1})J \varphi ^{-1} $, $\rho  = (\varrho \circ \varphi ^{-1})J \varphi ^{-1}$, $s  = (S \circ \varphi ^{-1} ) J \varphi ^{-1}$. Under the assumption \eqref{mat_cov_H}, the associated Hamiltonian function $H:T^* \operatorname{Diff}( \mathbb{R} ^n ) \times \Omega  ^n ( \mathbb{R} ^n ) \times \Omega  ^n ( \mathbb{R} ^n ) \rightarrow \mathbb{R}$ is $\operatorname{Diff}( \mathbb{R} ^n )$-invariant, namely
\begin{equation}\label{relabelling_symm_H} 
H\big( \varphi \circ \psi , (M \circ \psi) J \psi  , (\varrho \circ \psi) J \psi,(S \circ \psi) J \psi \big)=H( \varphi , M, \varrho , S ),\quad \text{for all $ \psi \in \operatorname{Diff}( \mathbb{R} ^n )$}.
\end{equation} 
Therefore, it induces a unique reduced Hamiltonian $h$, which turns out to be associated to the Hamiltonian density in \eqref{mathfrak_h} 
\begin{equation}\label{red_h} 
h: \big(\Omega ^1 ( \mathbb{R} ^n ) \otimes \Omega ^n ( \mathbb{R} ^n ) \big) \times \Omega ^n( \mathbb{R} ^n) \times \Omega ^n( \mathbb{R} ^n) \rightarrow \mathbb{R} , \quad h(m, \rho  , s)= \int_ \Omega \mathfrak{h}( m(x), \rho  (x), s(x)) {\rm d} x.
\end{equation}
Again, we emphasize the distinction between the function $h$, which is defined on forms on $ \mathbb{R}^n $, and the density $ \mathfrak{h} $ depending on the point values of these forms on $ \Omega $ only.

If the Legendre transform is invertible, the reduced Hamiltonian density can be equivalently defined in the following way :
\[
\mathfrak{h}(m,\rho,s) = m \cdot u - \mathfrak{l}(u, \rho, s)
\]
where $\mathfrak{l}$ is the reduced Lagrangian density and $u$ is the inverse image of $m$ via the Legendre transform $u \mapsto \frac{\partial \ell}{\partial u}$.

\begin{remark}[Reduction of the phase space]\label{remark_quotient_space_H}\rm We note that the passing from the Hamiltonian function $H:T^* \operatorname{Diff}( \mathbb{R} ^n ) \times \Omega ^n ( \mathbb{R} ^n ) \times \Omega  ^n ( \mathbb{R} ^n ) \rightarrow \mathbb{R}$ to its reduced version $h: \big(\Omega  ^1 ( \mathbb{R} ^n ) \otimes \Omega  ^n ( \mathbb{R} ^n ) \big) \times \Omega  ^n( \mathbb{R} ^n) \times \Omega  ^n( \mathbb{R} ^n) \rightarrow \mathbb{R}$, is simply understood as taking the function $h$ induced by $H$ on the quotient by the symmetry group $ \operatorname{Diff}(\mathbb{R} ^n)$:
\[
\Big[T^* \operatorname{Diff}( \mathbb{R} ^n ) \times \Omega  ^n ( \mathbb{R} ^n ) \times \Omega ^n ( \mathbb{R} ^n )\Big]/ \operatorname{Diff}( \mathbb{R} ^n ) \simeq  \big(\Omega ^1 ( \mathbb{R} ^n ) \otimes \Omega ^n ( \mathbb{R} ^n ) \big) \times \Omega ^n( \mathbb{R} ^n) \times \Omega^n( \mathbb{R} ^n).
\]
While this reduction process is formally the same as that for a closed fluid, the key difference is that it is applied to Hamiltonians in which the domain of the fluid enter in a nontrivial way, see \eqref{typical_H}, being found by integrating the density on $ \varphi ^{-1} ( \Omega )$. We stress again that the relabelling symmetry group is not $ \operatorname{Diff}( \Omega )$, with $ \Omega $ the fluid domain, but is the group $ \operatorname{Diff}( \mathbb{R} ^n )$.
\end{remark}


Based on these developments, we can state the following result which is the Hamiltonian counterpart of Proposition \ref{EP_open}.

\begin{proposition}[Hamiltonian variational principle for open fluids] \label{EP_open_H} Assume that the Hamiltonian density and the distributed and boundary sources of momentum, mass and entropy satisfy the relabelling symmetries. Then the following hold:
\begin{itemize}
\item[\bf (i)]
The Hamilton-d'Alembert phase space principle \eqref{extended_PS} yields the following Hamltonian variational formulation in the Eulerian frame:
\begin{equation}\label{ham_extended_reduced}
\begin{aligned} 
&\delta \int_0^T \!\!\int_ \Omega \Big[ m \cdot u + \rho D_t w + s D_t \gamma - \mathfrak{h}(m, \rho ,s )\Big]  {\rm d} x {\rm d} t \\
& \qquad \quad + \int_0^T \!\!\int_ \Omega \Big[b \cdot \zeta + \theta_ \rho   D_ \delta w  + \theta_s D_ \delta \gamma  \Big] {\rm d} x {\rm d} t\hspace{1.5cm}\leftarrow \text{bulk contribution}\\
& \qquad \quad  + \int_0^T \!\!\int_{ \partial \Omega } \Big[J \cdot \zeta + j _ \rho  D_ \delta w + j_ s D_ \delta \gamma \Big]   {\rm d} a {\rm d} t=0\;\;\;\;\leftarrow \text{boundary contribution}
\end{aligned} 
\end{equation}
with respect to variations $\delta u=\partial _t \zeta+ [ \zeta, u]$ and free variations $ \delta m$, $ \delta \rho  $, $ \delta s$, $ \delta w$, $ \delta \gamma $, with $ \delta w$ and $ \delta \gamma $ vanishing at $t=0,T$ and $\zeta $ an arbitrary time dependent vector field vanishing at $t=0,T$.

The relation with the material variables used in \eqref{extended_PS} is given by
\[
m= (M \circ  \varphi ^{-1})J \varphi ^{-1}  , \quad u = \dot  \varphi \circ \varphi ^{-1} ,
\]
\[
\rho  = (\varrho \circ \varphi ^{-1} )J \varphi ^{-1} , \quad s  = (S \circ \varphi ^{-1} )J \varphi ^{-1} , \quad w= W \circ \varphi ^{-1} , \quad \gamma = \Gamma \circ \varphi ^{-1} .
\]

\item[\bf (ii)] This principle yields the following equations and boundary conditions for open fluids:
\begin{equation}\label{LP_open} 
\left\{
\begin{array}{l}
\vspace{0.2cm}\displaystyle \partial _t  m + \pounds _{ \frac{\partial \mathfrak{h}}{\partial  m} } m + \rho  \nabla  \frac{\partial \mathfrak{h}}{\partial \rho  } + s  \nabla  \frac{\partial \mathfrak{h}}{\partial s }= b \\
\vspace{0.2cm}\displaystyle\partial _t \rho  + \operatorname{div} \Big(  \rho  \frac{\partial \mathfrak{h}}{\partial m} \Big)  = \theta_{\rho},\quad \partial _t s  + \operatorname{div} \Big(  s  \frac{\partial \mathfrak{h}}{\partial m} \Big)  = \theta_{s}\\
\vspace{0.2cm}\displaystyle \Big( \frac{\partial \mathfrak{h}}{\partial m} \cdot n \Big)  m =-J\quad\text{on}\quad \partial \Omega \\
\displaystyle\rho \frac{\partial \mathfrak{h}}{\partial m} \cdot n = -j_{\rho}, \qquad s\frac{\partial \mathfrak{h}}{\partial m} \cdot n = -j_{s} \quad\text{on}\quad \partial \Omega .
\end{array}
\right. 
\end{equation} 
\end{itemize}
\end{proposition}

\paragraph{The compressible Euler fluid with open boundaries.} The  reduced Hamiltonian density has the form:
\[
\mathfrak{h}(m, \rho ,s )= \frac{| m | ^2}{2 \rho }  + \varepsilon ( \rho  , s),
\]
from which the system \eqref{LP_open} gives the equations and boundary condition for the compressible fluid as follows  
\begin{equation}\label{Euler_rho_m} 
\left\{
\begin{array}{l}
\vspace{0.2cm}\displaystyle\partial _t  m + \operatorname{div}( m \otimes \frac{m}{\rho}  )  = - \nabla p + b\\
\vspace{0.2cm}\displaystyle\partial _t \rho  + \operatorname{div}( m) = \theta_ \rho  , \qquad \partial _t s  + \operatorname{div}( \frac{sm}{ \rho  } ) = \theta_s\\
\vspace{0.2cm}\displaystyle ( \frac{m}{\rho  }  \cdot n) \rho u =-J\quad\text{on}\quad \partial \Omega \\
\displaystyle m   \cdot n=- j _ \rho , \qquad  \frac{sm }{ \rho  }   \cdot n=- j _s\quad\text{on}\quad \partial \Omega .
\end{array}
\right. 
\end{equation} 

\paragraph{Rotating compressible fluid and shallow water equations.} They can be treated in a similar way by noting the the Hamiltonian densities are given by
\begin{align*}
\mathfrak{h} (m, \rho  , s)&= \frac{1}{2 \rho  ^2 }|m- \rho  R| ^2 + \varepsilon ( \rho  , s) + \rho  \phi \\
\mathfrak{h} (m, h)&=\frac{1}{2 h  ^2 }|m- h R| ^2 + \frac{1}{2} g(h+Z) ^2 .
\end{align*}

\subsection{Extended Lie-Poisson bracket formulation}\label{ELP}

In this paragraph we derive a bracket formulation for open fluids based on the above description. This formulation extends the Lie-Poisson bracket formulation of fluid dynamics (\cite{MoGr1980,DzVo1980,HoKu1983,MaWe1983,MaRaWe1984}) to fluids with open boundaries.

In the context of Hamiltonian systems, the most efficient way to derive Lie-Poisson brackets is to obtain them via reduction by symmetry of a canonical Poisson bracket on the phase space of the system, given by the cotangent bundle of the configuration manifold. For fluid dynamics, this corresponds to passing from the material to the Eulerian description. This is the approach taken in \cite{MaRaWe1984} for compressible fluids, in which the canonical Poisson bracket on $T^* \operatorname{Diff}( \Omega )$ corresponding to the material description, is reduced by the group of relabelling symmetries to give, in the Eulerian description, a Lie-Poisson bracket associated to a semidirect product Lie algebra, see \cite{MaWe1983} for the incompressible fluid. We refer to \cite{GBYo2019b,eldred2020single} for the derivation of bracket formalisms from a variational perspective for irreversible processes.

We shall take the same approach here, by deriving the bracket formulation associated to the material and Eulerian descriptions  described in \S\ref{Mat_H} and \S\ref{Eul_H} via reduction by symmetries.
This results in a suitable extension of the Lie-Poisson bracket. As a particular case of our bracket formulation, we derive derives the ``bulk+boundary" bracket formulation proposed in \cite{Ot2006}, see also \cite{BaMaBeMe2018}.

\paragraph{Bracket formulation in the material description.} In order to get the bracket formulation in the material description, we shall focus on functionals having the same form as the Hamiltonian for open fluids, namely $F: T^* \operatorname{Diff}( \mathbb{R} ^n ) \times \Omega ^n ( \mathbb{R} ^n ) \times \Omega ^n ( \mathbb{R} ) \rightarrow \mathbb{R} $ given in terms of a density as follows
\begin{equation}\label{functional_F} 
F( \varphi , M, \varrho , S)= \int_{ \varphi ^{-1} ( \Omega )} \mathfrak{F}( \varphi , M, \nabla \varphi , \varrho , S) {\rm d} X,
\end{equation}
for some fixed domain $ \Omega  \subset \mathbb{R} ^n $.
In a similar way with the Poisson formulation $ \dot  F=\{F,H\}$ of Hamiltonian systems, our goal is to find an expression for the time derivative of any functional $F$ of the form \eqref{functional_F} along a curve in $T^* \operatorname{Diff}( \mathbb{R} ^n ) \times \Omega ^n ( \mathbb{R} ^n ) \times \Omega ^n ( \mathbb{R} )$, which extends $\{F,H\}$ and completely characterizes the evolution equations and boundary conditions of the system \eqref{PS_force}.  
In order to get this expression, we compute the evolution of such functionals along the solution of \eqref{PS_force}, thereby giving
\begin{align}
\frac{d}{dt} F&=\int_{ \varphi ^{-1} ( \Omega )} \Big[\Big( \frac{\partial \mathfrak{F} }{\partial \varphi }- \operatorname{DIV}\frac{\partial \mathfrak{F} }{\partial \nabla \varphi } \Big)\cdot \frac{\partial \mathfrak{H} }{\partial M} - \Big( \frac{\partial \mathfrak{H} }{\partial \varphi }- \operatorname{DIV}\frac{\partial \mathfrak{H} }{\partial \nabla \varphi } \Big)\cdot \frac{\partial \mathfrak{F} }{\partial M} \Big]  {\rm d} X\label{line_1}\\
& \qquad + \int_{ \varphi ^{-1} ( \Omega )}\Big[\frac{\partial \mathfrak{F} }{\partial M } \cdot \mathfrak{B} + \frac{\partial \mathfrak{F} }{\partial \varrho  } \cdot \Theta _ \varrho + \frac{\partial \mathfrak{F} }{\partial S  } \cdot \Theta _ S\Big] {\rm d} X\label{line_2}\\
& \qquad + \int_{ \varphi ^{-1} (\partial  \Omega )} \frac{\partial \mathfrak{H} }{\partial M}  \cdot \Big(  \frac{\partial \mathfrak{F} }{\partial \nabla \varphi }-  \mathfrak{F} \, ( \nabla _x \varphi ^{-1} \circ \varphi )\Big) \cdot N {\rm d} A.\label{line_3}
\end{align}
We now multiply the three boundary conditions in \eqref{PS_force} by, respectively, $ \frac{\partial \mathfrak{F} }{\partial M}$, $ \frac{\partial \mathfrak{F} }{\partial \varrho }$, and $ \frac{\partial \mathfrak{F} }{\partial S}$, integrate them on $ \varphi ^{-1} ( \partial \Omega  )$ and add them to the above expression. The boundary integral \eqref{line_3} then remarkably combines with these added terms to yield the expression
\begin{align} 
\frac{d}{dt} F&= \int_{ \varphi ^{-1} ( \Omega )} \Big[\Big( \frac{\partial \mathfrak{F} }{\partial \varphi }- \operatorname{DIV}\frac{\partial \mathfrak{F} }{\partial \nabla \varphi } \Big)\cdot \frac{\partial \mathfrak{H} }{\partial M} - \Big( \frac{\partial \mathfrak{H} }{\partial \varphi }- \operatorname{DIV}\frac{\partial \mathfrak{H} }{\partial \nabla \varphi } \Big)\cdot \frac{\partial \mathfrak{F} }{\partial M} \Big]  {\rm d} X\label{term_1}\\
& \qquad +\int_{ \varphi ^{-1} ( \partial \Omega )} \frac{\partial \mathfrak{H}}{\partial M} \cdot \Big( \frac{\partial \mathfrak{F}}{\partial \nabla \varphi } + \Big( \frac{\partial \mathfrak{F} }{\partial M } \cdot M + \frac{\partial \mathfrak{F} }{\partial \varrho } \varrho + \frac{\partial \mathfrak{F} }{\partial S }S - \mathfrak{F}\Big) \nabla _x \varphi ^{-1} \circ \varphi \Big) \cdot N{\rm d} A\label{term_2}\\
& \qquad -\int_{ \varphi ^{-1} ( \partial \Omega )} \frac{\partial \mathfrak{F}}{\partial M} \cdot \Big( \frac{\partial \mathfrak{H}}{\partial \nabla \varphi } + \Big( \frac{\partial \mathfrak{H} }{\partial M } \cdot M + \frac{\partial \mathfrak{H} }{\partial \varrho } \varrho + \frac{\partial \mathfrak{H} }{\partial S }S - \mathfrak{H}\Big) \nabla _x \varphi ^{-1} \circ \varphi \Big) \cdot N{\rm d} A\label{term_3}\\
& \qquad + \underbrace{\int_{ \varphi ^{-1} ( \Omega )}\Big[\frac{\partial \mathfrak{F} }{\partial M } \cdot \mathfrak{B} + \frac{\partial \mathfrak{F} }{\partial \varrho  } \cdot \Theta _ \varrho + \frac{\partial \mathfrak{F} }{\partial S  } \cdot \Theta _ S\Big]{\rm d} X}_{  \rm bulk\; contribution}\label{term_4}\\
& \qquad +\underbrace{\int_{ \varphi ^{-1} ( \partial \Omega )} \Big[\frac{\partial \mathfrak{F}}{\partial M} \cdot \mathfrak{J}+ \frac{\partial \mathfrak{F}}{\partial \varrho }  \mathfrak{j}_ \varrho + \frac{\partial \mathfrak{F}}{\partial S}  \mathfrak{j}_S\Big] {\rm d} A}_{  \rm boundary\; contribution}\label{term_5}.
\end{align}
A direct verification then shows that the previous equality holding for all $F$ is equivalent to the system \eqref{PS_force}. The two terms \eqref{term_4} and \eqref{term_5} contain the contribution of each of the bulk and boundary source terms, which are paired with the corresponding partial derivative of $\mathfrak{F}$ on the domain and on the boundary. The nontrivial appearance of the boundary terms \eqref{term_2} and \eqref{term_3} is due to the domain of integration being a variable. Remarkably these two terms vanish for materially covariant densities $\mathfrak{F}$ and $ \mathfrak{H}$, see Remark \ref{cancellation_H}, which is always the case for fluid dynamics, thereby yielding
\[
\frac{d}{dt} F= \eqref{term_1}+ \eqref{term_4}  + \eqref{term_5} 
\]
for materially covariant densities. 
While the bracket \eqref{term_1}--\eqref{term_3} takes a somehow involved form which deserves further investigations of its geometric properties, it takes the usual Lie-Poisson form in the Eulerian description, as we shall show.

\begin{remark}[Property of materially covariant Hamiltonian densities]\label{cancellation_H}\rm When a Hamiltonian density $ \mathfrak{H}$ is materially covariant, see \eqref{mat_cov_H}, then its partial derivatives satisfy: 
\[
\frac{\partial \mathfrak{H} }{\partial \nabla \varphi } + \Big(\frac{\partial \mathfrak{H} }{\partial \varrho   }  \varrho +\frac{\partial \mathfrak{H} }{\partial S  }S + \frac{\partial \mathfrak{H} }{\partial M  } \cdot M - \mathfrak{H} \Big) \nabla _x \varphi ^{-1} \circ \varphi =0.
\]
This property is similar with the one mentioned in Remark \ref{cancellation} for Lagrangian densities, and is responsible for the cancellation of \eqref{term_2} and \eqref{term_3} in the materially covariant case. 
\end{remark}

\paragraph{Bracket formulation in the Eulerian description.} As we have seen earlier for the Hamiltonian, see \eqref{mathfrak_h} and \eqref{relabelling_symm_H}, when the density $ \mathfrak{F}$ in \eqref{functional_F} is materially covariant, the functional $F$ is $\operatorname{Diff}( \mathbb{R} ^n )$-invariant, and induces a functional $f:\big(\Omega ^1 ( \mathbb{R} ^n ) \otimes \Omega ^n ( \mathbb{R} ^n ) \big) \times \Omega ^n( \mathbb{R} ^n) \times \Omega ^n( \mathbb{R} ^n) \rightarrow \mathbb{R}$ of the form
\[
f(m, \rho  , s)= \int_ \Omega \mathfrak{f}(m(x), \rho  (x), s(x)) {\rm d} x,
\]
see Remark \ref{remark_quotient_space_H}. By taking the time derivative of such functionals on the solutions of system \eqref{LP_open}, we get the following result.

\begin{proposition} The system \eqref{LP_open} can be equivalently written in the extended Lie-Poisson form
\begin{equation}\label{Poisson_formulation} \begin{aligned} 
\frac{d}{dt}   f
&=\{f,h\}_{\rm LP} + \underbrace{\int_ \Omega \Big( b \cdot \frac{ \partial  \mathfrak{f} }{ \partial  m} + \theta_ \rho   \frac{\partial \mathfrak{f}}{\partial \rho  }+ \theta_ s   \frac{\partial \mathfrak{f}}{\partial s  } \Big)  {\rm d}x}_{  \rm bulk\; contribution} + \underbrace{
\int_{ \partial \Omega } \Big( J \cdot \frac{\partial \mathfrak{f}}{\partial m} + j_ \rho   \frac{\partial \mathfrak{f}}{\partial \rho  }+j_s \frac{\partial \mathfrak{f}}{\partial s  } \Big) {\rm d} a}_{ \rm boundary\; contribution}, \quad \text{for all $f$},
\end{aligned}
\end{equation} 
where
\begin{align*} 
\{f,h\}_{\rm LP}&= \int_ \Omega m \cdot \Big[ \frac{\partial \mathfrak{f}}{\partial m} , \frac{\partial \mathfrak{h}}{\partial m} \Big] {\rm d} x + \int_ \Omega \rho  \Big(  \frac{\partial \mathfrak{h}}{\partial m} \cdot \nabla  \frac{\partial \mathfrak{f}}{\partial \rho  } -
\frac{\partial \mathfrak{f}}{\partial m} \cdot \nabla \frac{\partial \mathfrak{h}}{\partial \rho  }\Big) {\rm d} x \\
& \qquad \qquad + \int_ \Omega s  \Big(  \frac{\partial \mathfrak{h}}{\partial m} \cdot \nabla  \frac{\partial \mathfrak{f}}{\partial s  } -
\frac{\partial \mathfrak{f}}{\partial m} \cdot \nabla \frac{\partial \mathfrak{h}}{\partial s  }\Big) {\rm d} x 
\end{align*} 
is the Lie-Poisson bracket.
\end{proposition}

In particular, taking $f=h$ gives the energy balance of the open fluid system \eqref{LP_open} as
\[
\frac{d}{dt} h = \underbrace{\int_ \Omega \left( b \cdot \frac{ \partial  \mathfrak{h} }{ \partial  m} + \theta_ \rho   \frac{\partial \mathfrak{h}}{\partial \rho  }+ \theta _ s\frac{\partial \mathfrak{h}}{\partial s  } \right)  {\rm d}x}_{  \rm bulk\; source} + \underbrace{
\int_{ \partial \Omega } \left( J \cdot \frac{\partial \mathfrak{h}}{\partial m} + j_ \rho   \frac{\partial \mathfrak{h}}{\partial \rho  } +j_s \frac{\partial \mathfrak{h}}{\partial s  }\right) {\rm d} a}_{ \rm boundary\; source},
\]
while the balance of momentum, mass, and entropy follow by appropriate choice for $ f$, namely
\begin{align*}
\frac{d}{dt} m_i&= \int_ \Omega \partial _i \Big( \mathfrak{h} - m \cdot  \frac{\partial \mathfrak{h} }{\partial m} - \rho  \frac{\partial \mathfrak{h} }{\partial \rho  } - s  \frac{\partial \mathfrak{h} }{\partial s} \Big) {\rm d} x + \int_ \Omega b_i {\rm d} x + \int_{ \partial \Omega } J_i {\rm d} a\\
\frac{d}{dt} \rho  &=\int_ \Omega \theta_ \rho   {\rm d} x + \int_{ \partial \Omega }j_ \rho   {\rm d} a\\
\frac{d}{dt} s&=\int_ \Omega \theta_s {\rm d} x + \int_{ \partial \Omega }j_s {\rm d} a.
\end{align*} 

\begin{remark}[On the Lie-Poisson bracket]\rm The Lie-Poisson bracket appearing above can be seen as the Lie-Poisson associated to the Lie algebra $ \mathfrak{X} ( \Omega )$ of \textit{all} smooth vector fields on $ \Omega $, i.e., without any conditions on the boundary.
We have obtained it via reduction by Lie group symmetry applied to a specific class of functionals appearing for open fluid systems, namely, of the form \eqref{functional_F}. It should be noted that it \textit{does not} come from the group $ \operatorname{Diff}( \Omega )$ of diffeomorphism $ \Omega $. First $ \operatorname{Diff}( \Omega )$ is not the Lie group configuration of our system (it is the configuration Lie group of a closed fluid in $ \Omega $), and, second, its Lie algebra is the space $ \mathfrak{X}_{\|} ( \Omega )$ of vector fields parallel to the boundary, and not the Lie algebra $ \mathfrak{X} ( \Omega )$ we are considering here.
\end{remark}

\paragraph{Particular cases.} In order to present particular cases and compare our bracket formulation with earlier works, it is useful to split the Lie-Poisson bracket in two terms, obtained by integration by parts in a way to isolate the partial derivatives $ \frac{\partial \mathfrak{f}}{\partial m}$, $ \frac{\partial \mathfrak{f}}{\partial \rho    } $, and $ \frac{\partial \mathfrak{f}}{\partial s} $. This results in the expression
\begin{align*} 
\{f,h\}_{\rm LP} &= - \int_ \Omega \Big[ \Big( \nabla _{ \frac{\partial \mathfrak{h} }{\partial m} } m + \Big( \nabla \frac{\partial \mathfrak{h}}{\partial m} \Big) ^\mathsf{T} \cdot m + m \operatorname{div} \frac{\partial \mathfrak{h}}{\partial m} + \rho  \nabla \frac{\partial \mathfrak{h}}{\partial \rho  } +s  \nabla \frac{\partial \mathfrak{h}}{\partial s  }  \Big) \cdot \frac{\partial \mathfrak{f}}{\partial m} \\
&\qquad \qquad \qquad \qquad -  \operatorname{div} \Big( \rho  \frac{\partial \mathfrak{h}}{\partial m  } \Big) \frac{\partial \mathfrak{f}}{\partial \rho  }- \operatorname{div} \Big(s  \frac{\partial \mathfrak{h}}{\partial m  } \Big) \frac{\partial \mathfrak{f}}{\partial s  }\Big] {\rm d} x \\
& \qquad + \int_{ \partial \Omega }\Big[ \Big( \frac{\partial \mathfrak{h}}{\partial m} \cdot m \Big) \Big( m \cdot \frac{\partial \mathfrak{f}}{\partial m} \Big) + \rho  \Big( \frac{\partial \mathfrak{h}}{\partial m}\cdot n \Big)\frac{\partial \mathfrak{f}}{\partial \rho  }+ s  \Big( \frac{\partial \mathfrak{h}}{\partial m}\cdot n \Big)\frac{\partial \mathfrak{f}}{\partial s  } \Big] {\rm d} a\\
&= :\{f,h\}_{\rm bulk} + \{f, h\}_{\rm boundary },
\end{align*} 
with $\{f,h\}_{\rm bulk}$, resp., $\{f, h\}_{\rm boundary }$, defined as the integral over $ \Omega $, resp., $ \partial \Omega $. With this notation, our bracket formulation \eqref{Poisson_formulation} reads
\begin{equation}\label{PB_formualtion}
\begin{aligned}
\dot f 
&= \{f,h\}_{\rm bulk} \\
& \qquad +  \underbrace{\int_{ \partial \Omega } \Big[ m\Big( \frac{\partial \mathfrak{h}}{\partial m} \cdot n \Big)  + J \Big]  \cdot \frac{\partial \mathfrak{f}}{\partial m} + \Big[ \rho  \Big( \frac{\partial \mathfrak{h}}{\partial m}\cdot n\Big)  +j_ \rho   \Big] \frac{\partial \mathfrak{f}}{\partial \rho  }+ \Big[ s  \Big( \frac{\partial \mathfrak{h}}{\partial m}\cdot n\Big)  +j_ s   \Big] \frac{\partial \mathfrak{f}}{\partial s  } {\rm d} a}_{= \{f,h\}_{\rm boundary} + \text{boundary contribution}}\\
& \qquad  +\underbrace{\int_ \Omega \Big( b \cdot \frac{ \partial  \mathfrak{f}}{ \partial  m} + \theta_ \rho   \frac{\partial \mathfrak{f}}{\partial \rho  } + \theta_ s   \frac{\partial \mathfrak{f}}{\partial s  } \Big)  {\rm d}x}_{ \text{ interior contribution}}
\end{aligned}
\end{equation} 
which shows how the boundary term $\{ \cdot , \cdot \}_{\rm boundary}$ can be combined with the term including the boundary fluxes $J$, $j_ \rho  $, and $j_ \rho  $.
We now consider several particular cases of the extended Lie-Poisson formulation. For this discussion we set
\[
\
b=0 \quad\text{and}\quad  \theta _ \rho  = \theta _s=0.
\]
\begin{itemize}
\item[(A)] If the fluxes are chosen as $J=0$ and $j_ \rho  = j_s=0$, then  the bracket formulation is
\[
\dot  f=\{f,h\}_{\rm LP}= \{f,h\}_{\rm bulk} + \{f,h\}_{\rm boundary}
\]
which yields the boundary conditions:
\begin{equation}\label{BC_A} 
\left( \frac{\partial \mathfrak{h}}{\partial m} \cdot n \right)  m =0, \quad \rho   \frac{\partial \mathfrak{h}}{\partial m}   \cdot n=0 \quad\text{and}\quad s   \frac{\partial \mathfrak{h}}{\partial m}   \cdot n=0 \quad\text{on}\quad \partial \Omega, 
\end{equation} 

the fluid is closed. It formally recovers the usual Lie-Poisson bracket from semidirect product theory. There is however a major difference in the formulation, which is due to the fact that we are considering the Lie-Poisson structure on the whole Lie algebra $ \mathfrak{X} ( \Omega )$, not $ \mathfrak{X} _{\|}( \Omega )$. In our case the Lie-Poisson formulation $ \dot  f=\{f,h\}_{\rm LP}$ for all $f$, not only produces the fluid equations of motion, but also imposes the boundary conditions
\eqref{BC_A}.

\item[(B)] If the boundary fluxes $J$, $j_ \rho  $, and $j_s$ are chosen as
\[
J:= - m\left( \frac{\partial h}{\partial m} \cdot m \right), \quad j_ \rho  := \rho  \left( \frac{\partial h}{\partial m}\cdot n\right) \quad\text{and}\quad j_s:= s  \left( \frac{\partial h}{\partial m}\cdot n\right),
\]
(recall that our approach allows them to depend on $m$, $ \rho  $, $s$), then the bracket formulation \eqref{Poisson_formulation} reduces to
\[
\dot  f=\{f,h\}_{\rm bulk}
\]
and the boundary conditions just vanish, and the fluid is freely open. It is important to remember that $\{f,h\}_{\rm bulk}$ is not a Poisson bracket. We can write it in terms of the Lie-Poisson bracket as $\{f,h\}_{\rm bulk}= \{f,h\}_{\rm LP}- \{f,h\}_{\rm boundary}$. This setting recovers the bracket formulation for open fluids in \cite{Ot2006}, see also \cite{BaMaBeMe2018}.

\item[(C)] The other cases are when $J$ and $j$ are prescribed, i.e., the fluid is open, but its behavior through the boundary is prescribed.
\end{itemize}

It is important to note that the boundary can be divided into non-intersecting pieces on which (A), (B) or (C) hold. For example there could be inflow along one part of the domain, no-flux along another, and open flow along a third part.

\section{Extensions}
\label{extensions}

We extend the previously developed framework to address several key aspects, including the treatment of multicomponent fluids, the consideration of general advected quantities, the analysis of higher-order fluids, and the incorporation of boundary stresses. As we shall see, some of these extensions lead to nontrivial modifications of the boundary conditions for open fluids.

\color{black}

\subsection{Multicomponent fluids}

To treat multicomponent fluids, one considers Lagrangian densities of the form $ \mathfrak{l}( u, \rho  _k, s)$, i.e., depending on the mass density of each component $k=1,...,N$. These single velocity, single temperature, multiple component models are ubiquitous in geophysical fluid dynamics \cite{gay2019variational,makarieva2017equations,bott2008theoretical,catry2007flux,eldred2020single,wacker2003continuity,bannon2002theoretical} and porous media \cite{helfferich1981theory,quintard2006nonlinear,chella1998multiphase,FaGBPu2020,GBPu2022}. The extension of Proposition \ref{EP_open} to this case is straightforward. In the material description, we replace the terms $ \varrho  \dot  W $, $ \Theta _{ \rho  } \delta W $ and $  \mathfrak{j} _{ \rho  } \delta W$ by the terms $\sum_k \varrho _k \dot  W_k $, $\sum_k \Theta _{ \rho  _k} \delta W _k$ and $\sum_k \mathfrak{j} _{ \rho  _k} \delta W_k$ in the Lagrange-d'Alembert principle \eqref{extended_HP}. In the Eulerian version this gives the variational principle
\begin{equation}\label{extended_reduced_multi}
\begin{aligned} 
&\delta \int_0^T \!\!\int_ \Omega \Big[ \mathfrak{l}(u, \rho ,s ) + \sum_k\rho_k D_t w_k + s D_t \gamma \Big]  {\rm d} x {\rm d} t \\
& \qquad \quad + \int_0^T \!\!\int_ \Omega \Big[b \cdot \zeta + \sum_k\theta_ {\rho_k}   D_ \delta w_k  + \theta_s D_ \delta \gamma  \Big] {\rm d} x {\rm d} t\\
& \qquad \quad  + \int_0^T \!\!\int_{ \partial \Omega } \Big[J \cdot \zeta + \sum_k j _ {\rho_k}  D_ \delta w_k + j_ s D_ \delta \gamma \Big]   {\rm d} a {\rm d} t=0
\end{aligned} 
\end{equation}
with respect to variations $\delta u=\partial _t \zeta+ [ \zeta, u]$ and free variations $ \delta \rho_k  $, $ \delta s$, $ \delta w_k$, $ \delta \gamma $, with $ \delta w_k$ and $ \delta \gamma $ vanishing at $t=0,T$ and $\zeta $ an arbitrary time dependent vector field vanishing at $t=0,T$. From this, the following system is obtained
\begin{equation}\label{Open_fluid_Eulerian_multi} 
\left\{
\begin{array}{l}
\vspace{0.2cm}\displaystyle\partial _t  \frac{\partial \mathfrak{l} }{\partial u} + \pounds _u \frac{\partial \mathfrak{l} }{\partial u} - \sum_k\rho_k  \nabla  \frac{\partial \mathfrak{l} }{\partial \rho _k }- s\nabla  \frac{\partial \mathfrak{l} }{\partial s  } = b\\
\vspace{0.2cm}\displaystyle\partial _t \rho_k  + \operatorname{div}( \rho_k  u) = \theta_ {\rho_k},\;\; \text{for all $k$}, \qquad \partial _t s  + \operatorname{div}( s  u) = \theta_s\\
\vspace{0.2cm}\displaystyle (u \cdot n) \frac{\partial \mathfrak{l} }{\partial u} =-J\quad\text{on}\quad \partial \Omega \\
\vspace{0.2cm}\displaystyle\rho_k   u   \cdot n=- j _ {\rho_k},\;\; \text{for all $k$} , \qquad  s   u   \cdot n=- j _s\quad\text{on}\quad \partial \Omega .
\end{array}
\right. 
\end{equation}
From this point, the variational formulation on the Hamiltonian side, as well as the extension of the Lie-Poisson bracket can be derived in a straightforward way.

\subsection{General advected quantities}
In this section we consider the case in which the fluid depends or other quantities than the densities $ \rho  $ and $s$. We treat the general case of advected tensor fields and tensor field densities. These sorts of fluid models occur in the complex fluids literature \cite{gay2009geometric,holm2002euler} and the plasma physics literature \cite{holm1983poisson,holm1986hydrodynamics,holm1986hamiltonian,morrison2005hamiltonian}, amongst others. For example, magnetohydrodynamics \cite{HoKu1983} can be understood as a compressible fluid with an advected magnetic field $2$-form.

\paragraph{Case of tensor field densities.} Let us assume that the Lagrangian density depends on a $(p,q)$ tensor field density, denoted $ \Pi \,{\rm d} X$ in the material description. Following the approach described in \S\ref{sec_2_1} for the mass and entropy densities, we consider a ``distributed source" $ \Theta _ \Pi{\rm d} X = \Theta _ \Pi ( \varphi , \dot  \varphi , \nabla \varphi , \varrho , K){\rm d} X $ of the tensorial quantity, given as a map
\[
\Theta _ \Pi {\rm d} X: \varphi ^{-1} ( \Omega ) \rightarrow T^p_q( \varphi ^{-1} ( \Omega )) \otimes \Lambda ^n(  \varphi ^{-1} ( \Omega ) )
\]
with $T^p_q( \varphi ^{-1} ( \Omega ))$ the bundle of $(p,q)$ tensors. We also consider a ``boundary source" $\mathfrak{J}_ \Pi {\rm d} A=\mathfrak{J}_ \Pi ( \varphi , \dot  \varphi , \nabla \varphi , \varrho , \Pi ) {\rm d} A$ of this tensorial quantity, given as a map
\[
\mathfrak{J}_ \Pi {\rm d} A: \varphi ^{-1} ( \partial \Omega ) \rightarrow T ^p_q( \varphi ^{-1} ( \partial \Omega )) \otimes \Lambda ^{n-1}( \varphi ^{-1} ( \partial \Omega )).
\]
The corresponding virtual displacement to be considered in this case, is a $(q,p)$ tensor field $Z$. Denoting by $ \Pi \!:\! Z$ the full contraction of tensors, the Lagrange-d'Alembert principle \eqref{extended_HP} adapted to this case becomes 
\begin{equation}\label{extended_HP_tensor_density}
\begin{aligned} 
&\left. \frac{d}{d\varepsilon}\right|_{\varepsilon=0}\int_0^T \!\!\int_{ \varphi _ \varepsilon ^{-1} ( \Omega )}\left[  \mathfrak{L}( \varphi _ \varepsilon , \dot  \varphi _ \varepsilon , \nabla \varphi _ \varepsilon, \varrho  _ \varepsilon , \Pi _ \varepsilon ) + \varrho_ \varepsilon  \dot  W_ \varepsilon + \Pi _ \varepsilon \!:\!\dot  Z_ \varepsilon \right]  {\rm d} X {\rm d} t\\
& \qquad  +  \int_0^T\!\!\int_{ \varphi ^{-1} ( \Omega)} \Big[ \mathfrak{B} \cdot \delta \varphi +\Theta_ \varrho   \delta W +\Theta_ \Pi  \!:\!\delta Z  \Big]
{\rm d} X {\rm d} t\;\;\;\;\;\leftarrow \text{bulk contribution}\\
&  \qquad  + \int_0^T\!\!\int_{ \varphi ^{-1} (\partial\Omega)} \Big[ \mathfrak{J}  \cdot \delta \varphi +\mathfrak{j}_ \varrho  \delta W  +\mathfrak{j}_\Pi \!:\! \delta Z   \Big] {\rm d} A {\rm d} t=0,\;\;\leftarrow \text{boundary contribution}
\end{aligned}
\end{equation} 
where we assume $ \delta \varphi $, $ \delta W$, and $ \delta Z $ vanish at $t=0,T$. We leave the derivation of equations in the material description to the reader and focus below on the Eulerian description.

\medskip

In a similar way with \eqref{B_Theta} and \eqref{J_j}, the relabelling symmetry of the tensorial source terms are 
\begin{equation}\label{Theta_J_Pi}
\begin{aligned} 
&\Theta_ \Pi   \big( \varphi \circ \psi , \dot  \varphi \circ \psi , \nabla_X ( \varphi \circ \psi ) ,(\varrho \circ \psi )J \psi,\psi ^* \Pi J \psi \big)= \psi ^* \big(\Theta _ \Pi  (\varphi , \dot \varphi , \nabla \varphi , \varrho , \Pi ) \big) J \psi\\
&\mathfrak{j} _ \Pi   \big( \varphi \circ \psi , \dot  \varphi \circ \psi , \nabla_X ( \varphi \circ \psi ) ,(\varrho \circ \psi )J \psi,\psi ^* \Pi J \psi \big)= \psi _ \partial ^*\big(\mathfrak{j}  _ \varrho (\varphi , \dot \varphi , \nabla \varphi , \varrho , \Pi ) \big) J \psi_ \partial ,
\end{aligned} 
\end{equation}
for all $ \psi \in \operatorname{Diff}( \mathbb{R} ^n )$, with $ \psi ^*$ and $ \psi _ \partial ^* $ denoting the pull-back of tensor fields. From this, one gets the corresponding Eulerian source terms  
\begin{align*} 
&\theta _ \pi  ( u, \rho  , \pi ){\rm d} x:  \Omega  \rightarrow T^p_q( \Omega ) \otimes \Lambda ^n ( \Omega )\\
&j_ \pi  (u, \rho  , \pi ) {\rm d} a : \partial \Omega \rightarrow T^p_q( \Omega )\Lambda ^{n-1} (  \partial \Omega ),
\end{align*} 
related to $ \Theta _ \Pi $ and $ \mathfrak{J}_ \Pi $ via
\begin{equation}\label{Theta_J_Pi_red}
\begin{aligned}
& \varphi _* \big(\Theta _ \Pi  (\varphi,...)\big)  J \varphi ^{-1}= \theta _ \pi ( u, \rho  , \pi )\\
&( \varphi _ \partial )_*\big(\mathfrak{j}_ \pi  (\varphi,...)\big) J \varphi _ \partial ^{-1} = j_ \pi  (u, \rho  , \pi ),
\end{aligned}
\end{equation} 
with $u= \dot  \varphi \circ \varphi  ^{-1} $ and $ \rho  = (\varrho \circ \varphi ^{-1} )J \varphi ^{-1} $ as earlier, and with the Eulerian version of $ \Pi $ defined by $ \pi   = (\varphi _* \Pi )J \varphi ^{-1} $. Above, we have used the abbreviation $(\varphi,...)= ( \varphi , \dot  \varphi , \varrho , \Pi )$.

We extend to tensor fields the notations $ \bar D_t \rho  $ and $D_tw$ used earlier for the Eulerian time derivatives of densities and functions: given a tensor field density $ \pi $ and a tensor field $ \kappa $, we use
\[
\bar D_t \pi := \partial _t \pi + \pounds _u \pi \quad\text{and}\quad D_t \kappa := \partial _t \kappa + \pounds _u \kappa,
\]
with $ \pounds _u \kappa $ and $ \pounds _u\pi $ the Lie derivative of the tensor field $ \kappa $ and the tensor field density $ \pi $. We recall their coordinate expressions
\begin{align*}
(\pounds _u \kappa )^{a_1...a_p}_{b_1...b_q}&= u^ c \partial _ c \kappa ^{a_1...a_p}_{b_1...b_q} - \kappa _{b_1...b_q}^{a_1...a_{r-1}c a_{r+1}...a_p} \partial _c u^{a_r}+ \kappa ^{a_1...a_p}_{b_1...b_{r-1} c b_{r+1}...b_q} \partial _{b_r} u^c\\
&=u^ c \partial _ c \kappa ^{a_1...a_p}_{b_1...b_q}+ \widehat{ \kappa }^{a_1...a_p c}_{b_1...b_q d} \partial _c u^d\\
(\pounds _u \pi)^{a_1...a_p}_{b_1...b_q} &= \partial _ c(u^ c  \pi ^{a_1...a_p}_{b_1...b_q}) - \pi _{b_1...b_q}^{a_1...a_{r-1}c a_{r+1}...a_p} \partial _c u^{a_r}+ \pi ^{a_1...a_p}_{b_1...b_{r-1} c b_{r+1}...b_q} \partial _{b_r} u^c\\
&=\partial _ c(u^ c  \pi ^{a_1...a_p}_{b_1...b_q})+ \widehat{ \pi }^{a_1...a_p c}_{b_1...b_q d} \partial _c u^d,
\end{align*} 
where we introduced the $(p+1,q+1)$ tensor field defined by
\begin{equation}\label{kappa_hat} 
\widehat{\kappa }^{a_1...a_pc}_{b_1...b_qd}=\sum_r \big(  \kappa ^{a_1...a_p}_{b_1...b_{r-1} d b_{r+1}...b_q} \delta ^c_{b_r}- \kappa _{b_1...b_q}^{a_1...a_{r-1}c a_{r+1}...a_p} \delta ^{a_r}_d\big),
\end{equation} 
similarly for $\widehat{ \pi }$. We can hence write these Lie derivatives as
\[
\pounds _u \kappa = \nabla _u \kappa + \widehat{ \kappa }\!:\! \nabla u \quad\text{and}\quad \pounds _u \pi = \operatorname{div} (u \pi) + \widehat{ \pi }\!:\! \nabla u.
\]
These notations are convenient to write explicitly the equations of motion below. We shall also use the following technical result, see \cite{GB2024}.

\begin{lemma}\label{GB2024_lemma} Let $ \kappa $ be a $(p,q)$ tensor field, $ \pi  $ a $(q,p)$ tensor field density, and $u$ a vector field. Then
\[
\pounds _ u \kappa \!:\! \pi = (u \cdot \nabla  \kappa )\!:\! \pi  - \operatorname{div}( \pi  \therefore \widehat{ \kappa }) \cdot u+ \operatorname{div}(( \pi \therefore \widehat{ \kappa }) \cdot u )   
\]
and
\[
\pounds _ u \pi \!:\! \kappa   = - \pi \!:\! (u \cdot \nabla \kappa   ) - \operatorname{div}( \kappa  \therefore \widehat{ \pi }) \cdot u+ \operatorname{div}\big( ( \pi \!:\! \kappa ) u + ( \kappa  \therefore \widehat{ \pi }) \cdot u \big),
\]
where $ \pi  \therefore \widehat{ \kappa  }$ is the $(1,1)$ tensor field density obtained by contracting all the respective
indices of $\widehat{ \kappa  }$ and $ \pi $ except the last covariant and contravariant indices of $\widehat{ \kappa  }$, similarly for  $\kappa  \therefore \widehat{ \pi }$.
Also, we have $ \pi \therefore \widehat{ \kappa }+ \kappa \therefore \widehat{ \pi }=0$ and $\pounds _ u \kappa  \!:\! \pi+ \pounds _ u \pi \!:\! \kappa= \operatorname{div}((\pi \!: \kappa ) u)$.
\end{lemma} 

We can now state the extension of Proposition \ref{EP_open} to the tensorial case.

\begin{proposition}[Eulerian variational principle for open fluids with tensorial advection I] \label{EP_open_tensor_density} Assume that the Lagrangian density and the distributed and boundary sources satisfy the relabelling symmetries. Then the following hold:
\begin{itemize}
\item[\bf (i)]
The Lagrange-d'Alembert principle \eqref{extended_HP_tensor_density} yields the following variational formulation in the Eulerian frame:
\begin{equation}\label{extended_reduced_tensor_density}
\begin{aligned} 
&\delta \int_0^T \!\!\int_ \Omega \Big[ \mathfrak{l}(u, \rho , \pi  ) + \rho D_t w + \pi \!:\! D_t z \Big]  {\rm d} x {\rm d} t \\
& \qquad \quad + \int_0^T \!\!\int_ \Omega \Big[b \cdot \zeta + \theta_ \rho   D_ \delta w  + \theta_ \pi \!:\! D_ \delta z  \Big] {\rm d} x {\rm d} t\hspace{1.2cm}\leftarrow \text{bulk contribution}\\
& \qquad \quad  + \int_0^T \!\!\int_{ \partial \Omega } \Big[J \cdot \zeta + j _ \rho  D_ \delta w + j_ \pi \!:\! D_ \delta z \Big]   {\rm d} a {\rm d} t=0\;\;\leftarrow \text{boundary contribution}
\end{aligned} 
\end{equation}
with respect to variations $\delta u=\partial _t \zeta+ [ \zeta, u]$ and free variations $ \delta \rho  $, $ \delta \pi $, $ \delta w$, $ \delta z $, with $ \delta w$ and $ \delta z $ vanishing at $t=0,T$ and $\zeta $ an arbitrary time dependent vector field vanishing at $t=0,T$.

The relation with the material variables used in \eqref{extended_HP_tensor_density} is given by
\[
u= \dot  \varphi \circ  \varphi ^{-1} , \quad \rho  = (\varrho \circ \varphi ^{-1} )J \varphi ^{-1} , \quad \pi   =  \varphi _* \Pi J \varphi ^{-1} , \quad w = W \circ \varphi ^{-1}, \quad z= \varphi  _*Z  .
\]

\item[\bf (ii)] This principle yields the following equations and boundary conditions for open fluids:
\begin{equation}\label{Open_fluid_Eulerian_tensor_density} 
\left\{
\begin{array}{l}
\vspace{0.2cm}\displaystyle\partial _t  \frac{\partial \mathfrak{l} }{\partial u} + \pounds _u \frac{\partial \mathfrak{l} }{\partial u} - \rho  \nabla  \frac{\partial \mathfrak{l} }{\partial \rho  }- \pi \!:\!\nabla  \frac{\partial \mathfrak{l} }{\partial \pi   } + \operatorname{div}\Big(\pi\therefore \widehat{\frac{\partial \mathfrak{l}}{\partial \pi } } \Big) = b\\
\vspace{0.2cm}\displaystyle\partial _t \rho  + \operatorname{div}( \rho  u) = \theta_ \rho  , \qquad \partial _t \pi   + \pounds _u \pi  = \theta_\pi \\
\vspace{0.2cm}\displaystyle (u \cdot n) \frac{\partial \mathfrak{l} }{\partial u} + \Big(\pi\therefore \widehat{\frac{\partial \mathfrak{l}}{\partial \pi } } \Big) \cdot n=-J\quad\text{on}\quad \partial \Omega \\
\vspace{0.2cm}\displaystyle\rho   u   \cdot n=- j _ \rho , \qquad  \pi    u   \cdot n=- j _\pi \quad\text{on}\quad \partial \Omega .
\end{array}
\right. 
\end{equation}
\end{itemize}
\end{proposition}
\noindent\textbf{Proof.}  The proof of Part {\bf (i)} is similar to that in Proposition \ref{EP_open}.
For {\bf (ii)}, the treatment of $ \rho  $ has already been done earlier, so it is enough to present the computation for the variation of $\int_0^T \int_ \Omega \left[  \mathfrak{l} (u, \pi) + \pi\!:\! D_tz \right]   {\rm d} x {\rm d} t $. The variation with respect to $ \pi  $ gives the condition 
\begin{equation}\label{z_condition} 
D_tz = - \frac{\partial \mathfrak{l} }{\partial \pi}.
\end{equation} 
The variation with respect to $u$ and $z$ is computed as follows:
\begin{align*}
&= \int_0^T \int_ \Omega  \left[ \frac{\partial \mathfrak{l} }{\partial u} \cdot \left( \partial _t \zeta + \pounds_ u \zeta \right) + \pi\!:\!  D_t \delta z + \pi \!:\! \pounds _{ \partial _t \zeta + \pounds_ u \zeta} z\right]{\rm d} x {\rm d} t \\
&= \int_0^T \int_ \Omega  \Bigg[ \left(- \partial _t \frac{\partial \mathfrak{l} }{\partial u} - \pounds_ u \frac{\partial \mathfrak{l} }{\partial u} \right) \cdot \zeta + \operatorname{div} \left( \Big( \frac{\partial \mathfrak{l} }{\partial u} \cdot \zeta \Big)  u \right)\\
& \qquad \qquad \qquad \qquad \qquad + \pi \!:\!D_t D_\delta z - \pi \!:\! D_t \pounds _ \zeta z +  \pi \!:\! \pounds _{ \partial _t \zeta + \pounds_ u \zeta} z \Bigg]{\rm d} x {\rm d} t  \\
&= \int_0^T \int_ \Omega  \Bigg[ \left(- \partial _t \frac{\partial \mathfrak{l} }{\partial u} - \pounds_ u \frac{\partial \mathfrak{l} }{\partial u} \right) \cdot \zeta + \operatorname{div} \left( \Big( \frac{\partial \mathfrak{l} }{\partial u} \cdot \zeta \Big)  u \right)\\
& \qquad \qquad \qquad \qquad \qquad + \bar D_t( \pi \!:\! D_\delta z ) - \bar D_t \pi \!:\!D_\delta z- \pi \!:\!\pounds _ \zeta D_t z \Bigg]{\rm d} x {\rm d} t  \\
&= \int_0^T \int_ \Omega  \Bigg[ \left(- \partial _t \frac{\partial \mathfrak{l} }{\partial u} - \pounds_ u \frac{\partial \mathfrak{l} }{\partial u} \right) \cdot \zeta + \operatorname{div} \left( \Big( \frac{\partial \mathfrak{l} }{\partial u} \cdot \zeta + \pi \!:\!  D_ \delta z \Big)  u \right) - \bar D_t \pi \!:\!D_\delta z + \pi \!:\!\pounds _ \zeta \frac{\partial \mathfrak{l}}{\partial \pi }  \Bigg]{\rm d} x {\rm d} t  \\
&=  \int_0^T \int_ \Omega  \Bigg[ \left(- \partial _t \frac{\partial \mathfrak{l} }{\partial u} - \pounds_ u \frac{\partial \mathfrak{l} }{\partial u} + \pi \!:\! \nabla \frac{\partial \mathfrak{l}}{\partial \pi }- \operatorname{div}\Big(  \pi \therefore \frac{\partial \mathfrak{l}}{\partial \pi } \Big)\right) \cdot \zeta  -\bar D_t \pi \!:\! D_ \delta z \Bigg]{\rm d} x {\rm d} t \\
& \qquad \qquad \qquad +  \int_0^T \int_ {\partial \Omega} \left[ \left(\frac{\partial \mathfrak{l} }{\partial u} \cdot \zeta \right) (u \cdot n) +  \Big(\pi \therefore \frac{\partial \mathfrak{l}}{\partial \pi }\Big) \cdot n+ (u \cdot n)\pi \!:\!  D_ \delta z\right] {\rm d} a {\rm d} t,
\end{align*}
where we used \eqref{z_condition}, we used $ \zeta $, $ \delta w$, and $ \delta z$ vanishing at $t=0,T$, and, in the last equality we used Lemma \ref{GB2024_lemma}.

Hence, when \eqref{extended_reduced_tensor_density} is considered, and by collecting the terms proportional to $ \zeta $, $ D_ \delta w$, and $D_ \delta z$, both at the interior and at the boundary, and using the equalities $D_tw = - \frac{\partial \mathfrak{l} }{\partial \rho  } $, $D_t z  = - \frac{\partial \mathfrak{l} }{\partial \pi  } $ arising from the variations $ \delta \rho  $ and $\delta \pi  $, we get the six equations \eqref{Open_fluid_Eulerian_tensor_density}. $\qquad\blacksquare$ 

\medskip 

The local and global energy balance equations are found as
\[
\partial _t \mathfrak{e}= \operatorname{div}\Big( - \Big(\frac{\partial \mathfrak{l} }{\partial u} \cdot u \Big)u  + \rho \frac{\partial \mathfrak{l}}{\partial \rho} u  + \pi \!:\! \frac{\partial \mathfrak{l}}{\partial \pi }  u +  \Big(\frac{\partial \mathfrak{l}}{\partial \pi }\therefore\widehat{ \pi }\Big)\cdot u  \Big) + b \cdot u - \frac{\partial \mathfrak{l}}{\partial \rho  }  \theta _ \rho  - \frac{\partial \mathfrak{l}}{\partial \pi   } \!:\! \theta _ \pi 
\]
and
\[
\frac{d}{dt} \int_ \Omega \mathfrak{e} {\rm d} x = \int_ \Omega \left[  b \cdot u  - \theta_{\rho} \frac{\partial \mathfrak{l}}{\partial \rho} - \theta_\pi \!:\! \frac{\partial \mathfrak{l}}{\partial \pi } \right]  {\rm d} x + \int_{\partial \Omega} \left[ J\cdot u - j_ \rho   \frac{\partial \mathfrak{l}}{\partial \rho  }- j_ \pi \!:\!   \frac{\partial \mathfrak{l}}{\partial  \pi   }  \right] {\rm d} a.
\]

\begin{remark}[Additional stress and examples]\rm
The dependence of the fluid Lagrangian on the tensor density $ \pi $ induces the additional stress
\[
\sigma _ \pi := \pi \therefore \widehat{\frac{\partial \mathfrak{l}}{\partial \pi }}
\]
appearing both in the fluid momentum equation and its boundary condition in \eqref{Open_fluid_Eulerian_tensor_density}. It can be explicitly computed depending on the nature of the tensor $\pi $, by using the general formula \eqref{kappa_hat}. For example, we get the following cases
\begin{align*}
& \text{$ \pi $ vector field density} & &(\sigma _ \pi)^c_d = \pi ^c \frac{\partial \mathfrak{l} }{\partial \pi ^d}\\
& \text{$ \pi $ 1-form density} & &(\sigma _ \pi)^c_d = -\pi _d \frac{\partial \mathfrak{l} }{\partial \pi _c}\\
& \text{$ \pi $ 2-contravariant symmetric density} & &(\sigma _ \pi)^c_d = 2\pi ^{ac} \frac{\partial \mathfrak{l} }{\partial \pi ^{ad}}\\
& \text{$ \pi $ 2-covariant symmetric density} & &(\sigma _ \pi)^c_d = -2\pi _{ad} \frac{\partial \mathfrak{l} }{\partial \pi _{ac}}.
\end{align*} 
We note that $ \sigma _ \pi $ are  (1,1) tensor field densities, hence the question of its symmetry needs to explicitly use a metric (here the Euclidean metric) for rising or lowering one of the indices. From the Lagrangian variational point of view we have here, this symmetry is related to the spatial (as opposed to material) covariance of the Lagrangian density. We refer to \cite{GBMaRa2012,GB2024} for these questions.
\end{remark}

\begin{remark}[Formulation on manifolds and bundle valued forms]\rm Although we have assumed a bounded domain in $\mathbb{R}^n$, the developments in this section can be extended to arbitrary Riemannian manifolds. Our treatment can also naturally be extended to (vector) bundle-valued differential forms \cite{KaArToYaMaDe2007,rashad2021exterior,califano2021geometric,gilbert2023geometric,rashad2023intrinsic}, of which $\kappa$ and $\pi$ are special cases (tensor-valued $0$-forms and $n$-forms, respectively). This will be the subject of future work.
\end{remark}

\paragraph{Case of tensor fields.} Let us now assume that the Lagrangian density depends on a $(p,q)$ tensor field (not tensor field density), denoted $K$ in the material description. The sources are hence given as $ \Theta _K= \Theta _K( \varphi , \dot  \varphi , \nabla \varphi , \varrho , K) $ with
\[
\Theta _K: \varphi ^{-1} ( \Omega ) \rightarrow T ^p_q( \varphi ^{-1} ( \Omega ))
\]
and $\mathfrak{J}_K=\mathfrak{J}_K( \varphi , \dot  \varphi , \nabla \varphi , \varrho , K)$ with
\[
\mathfrak{J}_K: \varphi ^{-1} ( \partial \Omega ) \rightarrow T ^p_q( \varphi ^{-1} ( \partial \Omega )).
\]
The corresponding virtual displacement to be considered in this case, is a $(q,p)$ tensor field density $Z \,{\rm d} X$ and the Lagrange-d'Alembert principle \eqref{extended_HP} adapted to this case becomes 
\begin{equation}\label{extended_HP_tensor}
\begin{aligned} 
&\left. \frac{d}{d\varepsilon}\right|_{\varepsilon=0}\int_0^T \!\!\int_{ \varphi _ \varepsilon ^{-1} ( \Omega )}\left[  \mathfrak{L}( \varphi _ \varepsilon , \dot  \varphi _ \varepsilon , \nabla \varphi _ \varepsilon, \varrho  _ \varepsilon , K_ \varepsilon ) + \varrho_ \varepsilon  \dot  W_ \varepsilon + K_ \varepsilon \!:\!\dot  Z_ \varepsilon \right]  {\rm d} X {\rm d} t\\
& \qquad  +  \int_0^T\!\!\int_{ \varphi ^{-1} ( \Omega)} \Big[ \mathfrak{B} \cdot \delta \varphi +\Theta_ \varrho   \delta W +\Theta_ K  \!:\!\delta Z  \Big]
{\rm d} X {\rm d} t\;\;\;\;\;\;\leftarrow \text{bulk contribution}\\
&  \qquad  + \int_0^T\!\!\int_{ \varphi ^{-1} (\partial\Omega)} \Big[ \mathfrak{J}  \cdot \delta \varphi +\mathfrak{j}_ \varrho  \delta W  +\mathfrak{j}_K \!:\! \delta Z   \Big] {\rm d} A {\rm d} t=0,\;\;\;\leftarrow \text{boundary contribution}
\end{aligned}
\end{equation} 
where we assume $ \delta \varphi $, $ \delta W$, and $ \delta \Gamma $ vanish at $t=0,T$.

\begin{remark}\rm
We note that the geometric setting in this case is different from that in \eqref{extended_HP_tensor_density} regarding the construction of the boundary integral. In \eqref{extended_HP_tensor_density} the boundary source term $\mathfrak{J}_ \Pi {\rm d} A$ is given as a tensor field density on the boundary, which can be contracted, on the boundary $ \varphi ^{-1} ( \partial \Omega )$, with the tensor field $ \delta Z$  to give the density on the boundary to be integrated. In the present case $ \delta Z {\rm d} X$ is a tensor field density on the domain $ \varphi ^{-1} ( \Omega )$ and one needs to induce from it a tensor field density on the boundary, that we denoted as $ \delta Z {\rm d} A$ in \eqref{extended_HP_tensor_density}, to be contracted with the tensor field $ \mathfrak{J}_K$ on the boundary. This step of passing from $ \delta Z {\rm d} X$ to $ \delta Z {\rm d} A$ needs the explicit use of a metric, here the Euclidean one, which was not necessary in \eqref{extended_HP_tensor_density}. This issue appears more clearly when formulating the framework on Riemannian manifolds, which will be pursued elsewhere.
\end{remark}

\begin{proposition}[Eulerian variational principle for open fluids with tensorial advection II] \label{EP_open_tensor_density_2} Assume that the Lagrangian density and the distributed and boundary sources satisfy the relabelling symmetries. Then the following hold:
\begin{itemize}
\item[\bf (i)]
The Lagrange-d'Alembert principle \eqref{extended_HP_tensor} yields the following variational formulation in the Eulerian frame:
\begin{equation}\label{extended_reduced_tensor}
\begin{aligned} 
&\delta \int_0^T \!\!\int_ \Omega \Big[ \mathfrak{l}(u, \rho , \kappa   ) + \rho D_t w + \kappa  \!:\! \bar D_t z \Big]  {\rm d} x {\rm d} t \\
& \qquad \quad + \int_0^T \!\!\int_ \Omega \Big[b \cdot \zeta + \theta_ \rho   D_ \delta w  + \theta_ \kappa  \!:\! \bar D_ \delta z  \Big] {\rm d} x {\rm d} t\hspace{1.2cm}\leftarrow \text{bulk contribution}\\
& \qquad \quad  + \int_0^T \!\!\int_{ \partial \Omega } \Big[J \cdot \zeta + j _ \rho  D_ \delta w + j_ \kappa  \!:\! \bar D_ \delta z \Big]   {\rm d} a {\rm d} t=0\;\;\leftarrow \text{boundary contribution}
\end{aligned} 
\end{equation}
with respect to variations $\delta u=\partial _t \zeta+ [ \zeta, u]$ and free variations $ \delta \rho  $, $ \delta \kappa  $, $ \delta w$, $ \delta z $, with $ \delta w$ and $ \delta z $ vanishing at $t=0,T$ and $\zeta $ an arbitrary time dependent vector field vanishing at $t=0,T$.

The relation with the material variables used in \eqref{extended_HP_tensor_density} is given by
\[
u= \dot  \varphi \circ  \varphi ^{-1} , \quad \rho  = (\varrho \circ \varphi ^{-1} )J \varphi ^{-1} , \quad \kappa    =  \varphi _* K , \quad w = W \circ \varphi ^{-1}, \quad z= \psi _*Z J \varphi ^{-1}  .
\]

\item[\bf (ii)] This principle yields the following equations and boundary conditions for open fluids:
\begin{equation}\label{Open_fluid_Eulerian_tensor_density_2}
\left\{
\begin{array}{l}
\vspace{0.2cm}\displaystyle\partial _t  \frac{\partial \mathfrak{l} }{\partial u} + \pounds _u \frac{\partial \mathfrak{l} }{\partial u} - \rho  \nabla  \frac{\partial \mathfrak{l} }{\partial \rho  }  - \kappa  \!:\!\nabla  \frac{\partial \mathfrak{l} }{\partial \kappa    } + \operatorname{div}\Big(\Big( \frac{\partial \mathfrak{l}}{\partial \kappa }: \kappa \Big) \delta +  \kappa \therefore \widehat{\frac{\partial \mathfrak{l}}{\partial\kappa  } } \Big) = b\\
\vspace{0.2cm}\displaystyle\partial _t \rho  + \operatorname{div}( \rho  u) = \theta_ \rho  , \qquad \partial _t \kappa    + \pounds _u \kappa  = \theta_ \kappa  \\
\vspace{0.2cm}\displaystyle (u \cdot n) \frac{\partial \mathfrak{l} }{\partial u} + \Big(\Big( \frac{\partial \mathfrak{l}}{\partial \kappa }: \kappa \Big) \delta +  \kappa \therefore \widehat{\frac{\partial \mathfrak{l}}{\partial\kappa  } } \Big) \cdot n=-J\quad\text{on}\quad \partial \Omega \\
\vspace{0.2cm}\displaystyle\rho   u   \cdot n=- j _ \rho , \qquad   \kappa     u   \cdot n=- j _ \kappa \quad\text{on}\quad \partial \Omega .
\end{array}
\right. 
\end{equation}
\end{itemize}
\end{proposition}

The local and global energy balance equations are found as
\[
\partial _t \mathfrak{e}= \operatorname{div}\Big( - \Big(\frac{\partial \mathfrak{l} }{\partial u} \cdot u \Big)u  + \rho \frac{\partial \mathfrak{l}}{\partial \rho} u  +  \Big(\frac{\partial \mathfrak{l}}{\partial \kappa  }\therefore\widehat{ \kappa  } \Big)\cdot u \Big) + b \cdot u - \frac{\partial \mathfrak{l}}{\partial \rho  }  \theta _ \rho  - \frac{\partial \mathfrak{l}}{\partial \kappa    } \!:\! \theta _ \kappa  
\]
and
\[
\frac{d}{dt} \int_ \Omega \mathfrak{e} {\rm d} x = \int_ \Omega \left[  b \cdot u  - \theta_{\rho} \frac{\partial \mathfrak{l}}{\partial \rho} - \theta_ \kappa  \!:\! \frac{\partial \mathfrak{l}}{\partial \kappa  } \right]  {\rm d} x + \int_{\partial \Omega} \left[ J\cdot u - j_ \rho   \frac{\partial \mathfrak{l}}{\partial \rho  }- j_ \kappa  \!:\!   \frac{\partial \mathfrak{l}}{\partial \kappa   }  \right] {\rm d} a.
\]

\paragraph{Examples.} The additional stress due to the dependence on the tensor $ \kappa $ is now given by
\[
\sigma _ \kappa  := \Big( \frac{\partial \mathfrak{l}}{\partial \kappa }: \kappa  \Big) \delta + \kappa  \therefore \widehat{\frac{\partial \mathfrak{l}}{\partial \kappa  }},
\]
which also appears both in the fluid momentum equation and its boundary condition in \eqref{Open_fluid_Eulerian_tensor_density_2}. It can be explicitly computed depending on the nature of the tensor $\kappa  $, by using the general formula \eqref{kappa_hat}. For example, we get the following cases
\begin{align*}
& \text{$\kappa  $ vector field} & &(\sigma _ \kappa )^c_d = \frac{\partial \mathfrak{l} }{\partial \kappa ^a} \kappa ^a \delta ^c_d + \kappa  ^c \frac{\partial \mathfrak{l} }{\partial \kappa  ^d}\\
& \text{$\kappa  $ 1-form} & &(\sigma _ \kappa )^c_d = \frac{\partial \mathfrak{l} }{\partial \kappa _a} \kappa _a \delta ^c_d- \kappa _d \frac{\partial \mathfrak{l} }{\partial \kappa _c}\\
& \text{$ \kappa  $ 2-form} & &(\sigma _ \kappa )^c_d = \frac{\partial \mathfrak{l} }{\partial \kappa _{ac}} \kappa _{ac} \delta ^c_d - 2 \kappa _{ad} \frac{\partial \mathfrak{l} }{\partial \kappa _{ac}}.
\end{align*} 

\paragraph{Magnetohydrodynamics with open boundaries.} In particular, for $ \kappa $ a closed 2-form, which we identify with a divergence free vector field $B$, this is the setting of magnetohydrodynamics. In this case, the equations \eqref{Open_fluid_Eulerian_tensor_density_2} become

\begin{equation}\label{mhd_equations}
\left\{
\begin{array}{l}
\vspace{0.2cm}\displaystyle\partial _t  \frac{\partial \mathfrak{l} }{\partial u} + \pounds _u \frac{\partial \mathfrak{l} }{\partial u} - \rho  \nabla  \frac{\partial \mathfrak{l} }{\partial \rho  } - B \cdot \nabla   \frac{\partial \mathfrak{l} }{\partial B   }  + \operatorname{div}\Big(B \otimes  \frac{\partial \mathfrak{l}}{\partial B }  \Big) = b\\
\vspace{0.2cm}\displaystyle\partial _t \rho  + \operatorname{div}( \rho  u) = \theta_ \rho  , \qquad \partial _t B    + \operatorname{curl}( B \times u)   = \theta_ B  \\
\vspace{0.2cm}\displaystyle (u \cdot n) \frac{\partial \mathfrak{l} }{\partial u} + \Big(B \otimes  \frac{\partial \mathfrak{l}}{\partial B }\Big) \cdot n=-J\quad\text{on}\quad \partial \Omega \\
\vspace{0.2cm}\displaystyle\rho   u   \cdot n=- j _ \rho , \qquad   B   u   \cdot n=- j _ B\quad\text{on}\quad \partial \Omega .
\end{array}
\right. 
\end{equation}

The local and global energy balance equations are found as
\[
\partial _t \mathfrak{e}= \operatorname{div}\Big( - \Big(\frac{\partial \mathfrak{l} }{\partial u} \cdot u \Big)u  + \rho \frac{\partial \mathfrak{l}}{\partial \rho} u  +( B \times u) \times\frac{\partial \mathfrak{l}}{\partial B } \Big) + b \cdot u - \frac{\partial \mathfrak{l}}{\partial \rho  }  \theta _ \rho  - \frac{\partial \mathfrak{l}}{\partial B   }\cdot  \theta _ B 
\]
and
\[
\frac{d}{dt} \int_ \Omega \mathfrak{e} {\rm d} x = \int_ \Omega \left[  b \cdot u  -  \frac{\partial \mathfrak{l}}{\partial \rho}\theta_{\rho} - \frac{\partial \mathfrak{l}}{\partial B   }\cdot  \theta _ B\right]  {\rm d} x + \int_{\partial \Omega} \left[ J\cdot u - j_ \rho   \frac{\partial \mathfrak{l}}{\partial \rho  }- j_ B\cdot   \frac{\partial \mathfrak{l}}{\partial  B}  \right] {\rm d} a.
\]
In particular, it is interesting to note how the term associated to $B$ in the local energy balance, i.e., $( B \times u) \times\frac{\partial \mathfrak{l}}{\partial B } $, contributes to both the momentum, $J \cdot u$, and magnetic, $j_B \cdot \frac{\partial \mathfrak{l}}{\partial B } $, effect in the boundary term of the global energy balance.
Using
\begin{equation}\label{2_contrib}
(B \times u)\times \frac{\partial \mathfrak{l}}{\partial B  }= \Big(\frac{\partial \mathfrak{l}}{\partial B  }\cdot B\Big) u - \Big(\frac{\partial \mathfrak{l}}{\partial B  }\cdot u\Big) B
\end{equation}
we get
\begin{align*}
&\int_\Omega\operatorname{div}\Big(- \Big(\frac{\partial \mathfrak{l} }{\partial u} \cdot u \Big)u+ (B \times u)\times \frac{\partial \mathfrak{l}}{\partial B  }\Big) {\rm d}x \\
&= \int_{\partial\Omega} \left[- \Big(\frac{\partial \mathfrak{l} }{\partial u} \cdot u \Big)u\cdot n - \Big(\frac{\partial \mathfrak{l}}{\partial B  }\cdot u\Big) B\cdot n+ \Big(\frac{\partial \mathfrak{l}}{\partial B  }\cdot B\Big) u\cdot n   \right] {\rm d} a\\
&=\int_{\partial \Omega} \left[ J\cdot u - j_ B\cdot   \frac{\partial \mathfrak{l}}{\partial  B}  \right] {\rm d} a,
\end{align*}
showing the different role played by both terms in \eqref{2_contrib}. It also explains the occurrence of the magnetic field contribution $\big(B \otimes  \frac{\partial \mathfrak{l}}{\partial B }\big) \cdot n$ in the open boundary condition for the fluid momentum.

\color{black}

\subsection{Higher order fluids}\label{HOF}

We consider the case of a Lagrangian density depending also on the first order derivative of $\rho$. An example of a fluid model with this feature is the Euler-Kortweg equations for multiphase flow \cite{bresch2019navier}. As we shall see this induces the appearance of an additional boundary condition. Given $\mathfrak{l}(u, \rho  , \nabla \rho ) $, we consider the same Lagrange-d'Alembert principle as in \eqref{extended_HP}, in which we do not include the entropy for simplicity. We thus have
\begin{equation}\label{extended_reduced_nable}
\delta \int_0^T \!\!\int_ \Omega \Big[ \mathfrak{l}(u, \rho , \nabla \rho   ) + \rho D_t w \Big]  {\rm d} x {\rm d} t + \int_0^T \!\!\int_ \Omega \Big[b \cdot \zeta + \theta_ \rho   D_ \delta w   \Big] {\rm d} x {\rm d} t  + \int_0^T \!\!\int_{ \partial \Omega } \Big[J \cdot \zeta + j _ \rho  D_ \delta w  \Big]   {\rm d} a {\rm d} t=0
\end{equation}
with respect to variations $\delta u=\partial _t \zeta+ [ \zeta, u]$ and free variations $ \delta \rho  $, $ \delta w$, with $ \delta w$ vanishing at $t=0,T$ and $\zeta $ an arbitrary time dependent vector field vanishing at $t=0,T$.

This principle yields the following equations and boundary conditions for open fluids:
\begin{equation}\label{Open_fluid_Eulerian_nabla} 
\left\{
\begin{array}{l}
\vspace{0.2cm}\displaystyle\partial _t  \frac{\partial \mathfrak{l} }{\partial u} + \pounds _u \frac{\partial \mathfrak{l} }{\partial u} - \rho  \nabla \Big( \frac{\partial \mathfrak{l} }{\partial \rho  }- \operatorname{div} \frac{\partial \mathfrak{l} }{\partial \nabla \rho  } \Big)    = b\\
\vspace{0.2cm}\displaystyle\partial _t \rho  + \operatorname{div}( \rho  u) = \theta_ \rho \\
\vspace{0.2cm}\displaystyle (u \cdot n) \frac{\partial \mathfrak{l} }{\partial u} =-J\quad\text{on}\quad \partial \Omega \\
\vspace{0.2cm}\displaystyle\rho   u   \cdot n=- j _ \rho , \qquad  \frac{\partial \mathfrak{l} }{\partial \nabla \rho  } \cdot n   =0\quad\text{on}\quad \partial \Omega .
\end{array}
\right. 
\end{equation}
We note the occurrence of the last boundary condition on $ \frac{\partial \mathfrak{l}}{\partial \nabla \rho  }$. Its role in the energy balance is illustrated as follows. System \eqref{Open_fluid_Eulerian_nabla} implies the following local energy balance for the total energy density $ \mathfrak{e}( u, \rho  , \nabla \rho  )= \frac{\partial \mathfrak{l}}{\partial u} \cdot u - \mathfrak{l}(u, \rho  , \nabla \rho)$ defined from $ \mathfrak{l}$:
\begin{align*}
\partial _t \mathfrak{e} &=\operatorname{div}  \Big(- \Big(\frac{\partial \mathfrak{l} }{\partial u} \cdot u \Big)u  + \rho \Big(\frac{\partial \mathfrak{l}}{\partial \rho} - \operatorname{div}\frac{\partial \mathfrak{l}}{\partial \nabla \rho}\Big) u + \frac{\partial \mathfrak{l}}{\partial \nabla \rho}( \operatorname{div}( \rho  u)- \theta _ \rho  )  \Big)\\
& \qquad \qquad  + b \cdot u  - \theta_{\rho} \Big(\frac{\partial \mathfrak{l}}{\partial \rho} - \operatorname{div}\frac{\partial \mathfrak{l}}{\partial \nabla \rho}\Big).
\end{align*} 
This should be compared to \eqref{loc_energy_Eulerian}. 
Now, by integrating it over the domain, we use the four boundary conditions in \eqref{Open_fluid_Eulerian_nabla} to finally get 
\[
\frac{d}{dt} \int_ \Omega \mathfrak{e} \,{\rm d} x = \int_ \Omega \Big[ b \cdot u - \theta _ \rho  \Big(\frac{\partial \mathfrak{l}}{\partial \rho} - \operatorname{div}\frac{\partial \mathfrak{l}}{\partial \nabla \rho}\Big)\Big] {\rm d} x+ \int_{ \partial \Omega } \Big[ J \cdot u- j_ \rho  \Big(\frac{\partial \mathfrak{l}}{\partial \rho} - \operatorname{div}\frac{\partial \mathfrak{l}}{\partial \nabla \rho}\Big)  \Big] {\rm d} a.
\]
We note in particular, how the additional condition on $ \frac{\partial \mathfrak{l}}{\partial \nabla \rho  }$ is used so that the term $\frac{\partial \mathfrak{l}}{\partial \nabla \rho}( \operatorname{div}( \rho  u)- \theta _ \rho  )$ has no boundary flow contribution.

\subsection{Boundary stresses}

As a final extension, we add the contribution of a (bulk) external stress tensor $\sigma$. This is a common feature in geophysical fluid dynamics models, used to model surface stresses due to the wind in ocean models and the land or waves in atmospheric models \cite{kraus1994atmosphere}. By appropriately adding the stress contribution in the Euler-Poincar\'e-d'Alembert principle \eqref{extended_reduced}, we get:
\begin{equation}\label{extended_reduced_sigma}
\begin{aligned}
\delta \int_0^T \!\!\int_ \Omega \Big[ \mathfrak{l}(u, \rho  ) + \rho D_t w \Big]  {\rm d} x {\rm d} t &+ \int_0^T \!\!\int_ \Omega \Big[b \cdot \zeta - \sigma \!:\! \nabla \zeta + \theta_ \rho   D_ \delta w   \Big] {\rm d} x {\rm d} t  \\
&+ \int_0^T \!\!\int_{ \partial \Omega } \Big[J \cdot \zeta + j _ \rho  D_ \delta w  \Big]   {\rm d} a {\rm d} t=0.
\end{aligned}
\end{equation}
This gives
\begin{equation}\label{Open_fluid_Eulerian_sigma} 
\left\{
\begin{array}{l}
\vspace{0.2cm}\displaystyle\partial _t  \frac{\partial \mathfrak{l} }{\partial u} + \pounds _u \frac{\partial \mathfrak{l} }{\partial u} - \rho  \nabla  \frac{\partial \mathfrak{l} }{\partial \rho  }   = b +  \operatorname{div} \sigma\\
\vspace{0.2cm}\displaystyle\partial _t \rho  + \operatorname{div}( \rho  u) = \theta_ \rho \\
\vspace{0.2cm}\displaystyle (u \cdot n) \frac{\partial \mathfrak{l} }{\partial u} =-J + \sigma \cdot n\quad\text{on}\quad \partial \Omega \\
\vspace{0.2cm}\displaystyle\rho   u   \cdot n=- j _ \rho \quad\text{on}\quad \partial \Omega 
\end{array}
\right. 
\end{equation}
in which we note the appearance of the external stress tensor in the boundary condition for the fluid momentum.

The local and balance equations are found
\begin{align*}
\partial _t \mathfrak{e} &=\operatorname{div}  \Big(- \Big(\frac{\partial \mathfrak{l} }{\partial u} \cdot u \Big)u  + \rho \frac{\partial \mathfrak{l}}{\partial \rho}  u \Big)  + b \cdot u + \operatorname{div} \sigma \cdot u   - \theta_{\rho} \frac{\partial \mathfrak{l}}{\partial \rho}
\end{align*} 
and
\[
\frac{d}{dt} \int_ \Omega \mathfrak{e}\, {\rm d} x = \int_ \Omega \Big[ b \cdot u - \theta _ \rho  \frac{\partial \mathfrak{l}}{\partial \rho} - \sigma : \nabla u\Big] {\rm d} x+ \int_{ \partial \Omega } \Big[ J \cdot u- j_ \rho   \frac{\partial \mathfrak{l}}{\partial \rho} \Big] {\rm d} a.
\]


\section{Conclusions}
\label{conclusions}
We have developed a geometric framework based on Lie groups for the motion of fluids with fixed, permeable boundaries, that extends the usual geometric desciption in closed domains. Our setting is based on the classic approach with appropriate extensions to incorporate bulk and boundary forces via a Lagrange-d’Alembert approach, taking into account only the fluid parcels present in the domain. Using this, we were able to develop both a variational formulation in the Eulerian description that extends the Euler-Poincar\'e framework to open fluids; and a bracket-based formulations that extends the Lie-Poisson bracket. It is capable of reproducing, as particular cases, existing results from the GENERIC literature \cite{Ot2006,BaMaBeMe2018}. A comparison with related port-Hamiltonian results \cite{rashad2021portA,rashad2021portB,lohmayer2022exergetic,mora2023irreversible,cardoso2021dissipative,haine2021incompressible,rashad2021exterior,califano2021geometric} will be the subject of future work. This general framework was illustrated with several examples (shallow water and compressible Euler), and extensions to multicomponent fluids, general tensor advected quantities, higher-order fluids and boundary stresses were demonstrated. In future work we will consider the extension of these ideas to moving boundaries, following the approach in \cite{GBMaRa2012}; and also an extension to interactive coupling of multiple fluid domains.

\appendix

\section{Equivalence between Lagrangian and Eulerian open fluid equations}
\label{equivalance_lag_eul}

We have seen in Proposition \ref{EP_open} that the variational principles in regard to the action integrals in the Lagrangian perspective (\ref{extended_HP}) and the Eulerian perspective (\ref{extended_reduced}) are equivalent. We can also show directly that the Euler-Poincar\'e equations derived from the two variational principles are equivalent by change of perspective. To see this, we first calculate the relationship between the derivatives of reduced and non-reduced Lagrangians. 

In fact, from the relation 
\[
\mathfrak{L}(\varphi, \dot \varphi, \nabla \varphi, \varrho, S) (X)= \mathfrak{l} (u, \rho, s) \circ \varphi(X) J \varphi(X) 
\]
we have
\[
\mathfrak{L}( \varphi (X), \dot  \varphi (X), \nabla \varphi (X), \varrho (X))= \mathfrak{l} \left( \dot \varphi (X), \frac{\varrho (X)}{J \varphi (X)} ,\frac{S (X)}{J \varphi (X)} \right)  J \varphi (X)
\]
by definition of $u$, $\rho$ and $s$. From this equation, we have the following equalities:
\begin{align*}
\frac{\partial \mathfrak{L}}{\partial \varphi }&=0 \\
\frac{\partial \mathfrak{L}}{\partial \dot  \varphi } &= \frac{\partial \mathfrak{l}}{\partial u} \circ \varphi J \varphi \\
\frac{\partial \mathfrak{L}}{\partial \nabla_X \varphi }&= \left( \left( \mathfrak{l} - \rho  \frac{\partial \mathfrak{l}}{\partial \rho  }- s \frac{\partial \mathfrak{l}}{\partial s  }  \right) \circ \varphi\right) \left( \nabla_x \varphi ^{-1}  \circ \varphi \right) J \varphi \\
\frac{\partial \mathfrak{L}}{\partial \varrho } &= \frac{\partial \mathfrak{l}}{\partial \rho  }\circ \varphi \\
\frac{\partial \mathfrak{L}}{\partial S } &= \frac{\partial \mathfrak{l}}{\partial s  }\circ \varphi .
\end{align*}

\noindent{\bf (1)} We first prove the equivalence of the following two equations: 
\[D_t \frac{\partial \mathfrak{L}}{\partial \dot  \varphi }  +\operatorname{DIV} \frac{\partial \mathfrak{L}}{\partial \nabla \varphi } - \frac{\partial \mathfrak{L}}{\partial \varphi } = \mathfrak{B}\]
\[
\partial _t  \frac{\partial \mathfrak{l} }{\partial u} + \pounds _u \frac{\partial \mathfrak{l} }{\partial u} - \rho  \nabla  \frac{\partial \mathfrak{l} }{\partial \rho  }- s\nabla  \frac{\partial \mathfrak{l} }{\partial s  } = b.
\]

First of all, taking the material time derivative on both sides of the relation $\frac{\partial \mathfrak{L}}{\partial \dot  \varphi } = \frac{\partial \mathfrak{l}}{\partial u} \circ \varphi J \varphi$, we have that: 
\[
D_t \frac{\partial \mathfrak{L}}{\partial \dot  \varphi } = \partial _t  \frac{\partial \mathfrak{l} }{\partial u} \circ \varphi J \varphi + \nabla \frac{\partial \mathfrak{l} }{\partial u} \circ \varphi \cdot \dot \varphi J \varphi + \left( \frac{\partial \mathfrak{l} }{\partial u} \mathrm{div} u \right) \circ \varphi J \varphi.
\]
Then, we have that:
\[
\operatorname{DIV} \frac{\partial \mathfrak{L}}{\partial \nabla \varphi } = \nabla_X \left( \left( \ell- \rho  \frac{\partial \ell}{\partial \rho  }-s  \frac{\partial \ell}{\partial s  }  \right) \circ \varphi \right) \cdot  \left( \nabla_x \varphi ^{-1}  \circ \varphi \right) J \varphi.
\]
In local coordinates, it is:
\[
\left( \frac{\partial \mathfrak{L}}{\partial \nabla \varphi^i_{,A} } \right)_{,A}= \left( \left( \ell- \rho  \frac{\partial \ell}{\partial \rho  }-s  \frac{\partial \ell}{\partial s  }  \right) \circ \varphi \right)_{,A}  ( \varphi ^{-1} )^A_{,i} \circ \varphi J \varphi,
\]
where we have used Piola's identity:
\[
\operatorname{DIV} \left(\left( \nabla_x \varphi ^{-1}  \circ \varphi \right) J \varphi\right) = 0.
\]
Next we have that:
\begin{align*}
&\left( \left( \ell- \rho  \frac{\partial \ell}{\partial \rho  }-s  \frac{\partial \ell}{\partial s  }  \right) \circ \varphi \right)_{,A}  ( \varphi ^{-1} )^A_{,i} \circ \varphi J \varphi =  \left( \left( \ell- \rho  \frac{\partial \ell}{\partial \rho  }-s  \frac{\partial \ell}{\partial s  }  \right)_{,k}  \circ \varphi \right)  \varphi^k_{,A} ( \varphi ^{-1} )^A_{,i} \circ \varphi J \varphi \\
& \; = \left( \left( \ell- \rho  \frac{\partial \ell}{\partial \rho  }-s  \frac{\partial \ell}{\partial s  }  \right)_{k}  \circ \varphi \right) e_{ki} J \varphi  = \left( \left( \ell- \rho  \frac{\partial \ell}{\partial \rho  }-s  \frac{\partial \ell}{\partial s  }  \right)_{,i}  \circ \varphi \right) J \varphi \\
& \; =\left( \left( \frac{\partial \ell}{\partial u^j} u^j_{,i} + \frac{\partial \ell}{\partial \rho} \rho_{,i} + \frac{\partial \ell}{\partial s} s_{,i} - \frac{\partial \ell}{\partial \rho} \rho_{,i} - \rho \left( \frac{\partial \ell}{\partial \rho}\right)_{,i} - \frac{\partial \ell}{\partial s} s_{,i} - s \left( \frac{\partial \ell}{\partial s}\right)_{,i}\right) \circ \varphi \right) J \varphi \\
& \; =\left( \left( \frac{\partial \ell}{\partial u^j} u^j_{,i}  - \rho \left( \frac{\partial \ell}{\partial \rho}\right)_{,i} - s \left( \frac{\partial \ell}{\partial s}\right)_{,i}\right) \circ \varphi \right) J \varphi.
\end{align*}
Thus 
\[
\operatorname{DIV} \frac{\partial \mathfrak{L}}{\partial \nabla \varphi } = \left( \left( (\nabla u)^{\intercal} \cdot \frac{\partial \ell}{\partial u}  - \rho  \nabla \frac{\partial \ell}{\partial \rho} - s \nabla \frac{\partial \ell}{\partial s}\right) \circ \varphi \right) J \varphi.
\]
Hence 
\begin{align*}
\partial_t \frac{\partial \mathfrak{L}}{\partial \dot  \varphi }  +\operatorname{DIV} \frac{\partial \mathfrak{L}}{\partial \nabla \varphi } - \frac{\partial \mathfrak{L}}{\partial \varphi } &= \left( \partial _t  \frac{\partial \mathfrak{l} }{\partial u}  + \left( \frac{\partial \mathfrak{l} }{\partial u} \mathrm{div} u \right) + (\nabla u)^{\intercal} \cdot \frac{\partial \mathfrak{l}}{\partial u} + u \cdot \nabla \frac{\partial \mathfrak{l}}{\partial u}  - \rho  \nabla \frac{\partial \mathfrak{l}}{\partial \rho} - s \nabla \frac{\partial \mathfrak{l}}{\partial s} \right) \circ \varphi J \varphi \\
&= \left( \partial _t  \frac{\partial \mathfrak{l} }{\partial u} + \pounds _u \frac{\partial \mathfrak{l} }{\partial u} - \rho  \nabla  \frac{\partial \mathfrak{l} }{\partial \rho  }- s\nabla  \frac{\partial \mathfrak{l} }{\partial s  } \right) \circ \varphi J \varphi .
\end{align*}
On the other hand, by definition, 
\[
b \circ \varphi J \varphi  = \mathfrak{B}.
\]
We therefore have the equivalence.
\[
D_t \frac{\partial \mathfrak{L}}{\partial \dot  \varphi }  +\operatorname{DIV} \frac{\partial \mathfrak{L}}{\partial \nabla \varphi } - \frac{\partial \mathfrak{L}}{\partial \varphi } = \mathfrak{B} \quad \Leftrightarrow \quad 
\partial _t  \frac{\partial \mathfrak{l} }{\partial u} + \pounds _u \frac{\partial \mathfrak{l} }{\partial u} - \rho  \nabla  \frac{\partial \mathfrak{l} }{\partial \rho  }- s\nabla  \frac{\partial \mathfrak{l} }{\partial s  } = b.
\]

\noindent{\bf (2)} Regarding the balance of mass, we have 
\[
D_t \varrho = D_t (\rho \circ \varphi J \varphi) = \partial_t \rho \circ \varphi J \varphi + \operatorname{div}(\rho u) \circ \varphi J \varphi.
\]
On the other hand
\[
\theta_{\varrho} \circ \varphi J \varphi  = \Theta_{\rho}.
\]
Thus we have 
\[ D_t \varrho =  \Theta_{\varrho} \quad \Leftrightarrow \quad \partial_t \rho  + \operatorname{div}(\rho u) = \theta_{\rho}. \]
A similar proof shows the equivalence 
\[
D_t S =  \Theta_S \quad \Leftrightarrow \quad \partial_t s  + \operatorname{div}(s u) = \theta_{s}.
\]

\noindent{\bf (3)} 
On the boundary $ \partial (\varphi^{-1}(\Omega))$, we have that for any vector field $w(x)$:
\begin{equation} \label{boundary_piola}
    N \cdot ((\nabla_x \varphi ^{-1} \cdot w )\circ \varphi)  J \varphi   = (n \cdot w) \circ \varphi_{\partial} J \varphi _ \partial ,
\end{equation}
where $N(X)$ at $X \in \partial (\varphi^{-1}(\Omega))$ and $n(x)$ at $x = \varphi(X) \in \partial (\Omega)$ are the normal vectors on the boundary pointing outward. 
In fact, due to Piola's identity,
\[
\operatorname{DIV} \left( ((\nabla_x \varphi ^{-1} \cdot w )\circ \varphi)  J \varphi \right) = \operatorname{div} (w) \circ \varphi J \varphi .
\]
Applying the divergence theorem on both sides and a change of variables, we have: 
\[
\int_{ \partial (\varphi^{-1}(\Omega))}N \cdot ((\nabla_x \varphi ^{-1}  \cdot w )\circ \varphi)  J \varphi {\rm d} A  = \int_{ \partial \Omega } n \cdot w {\rm d} a.
\]
Using the change of variables $\varphi_{\partial} : X \mapsto x $ on the boundary integral, we also have:
\[
\int_{ \partial \Omega } n \cdot w {\rm d} a = \int_{ \partial (\varphi^{-1}(\Omega))} (n \cdot w) \circ \varphi_{\partial} J \varphi_{\partial} {\rm d} A
\]
so 
\[
\int_{ \partial (\varphi^{-1}(\Omega))}N \cdot ((\nabla_x \varphi ^{-1}  \cdot w )\circ \varphi)  J \varphi {\rm d} A  = \int_{ \partial (\varphi^{-1}(\Omega))} (n \cdot w) \circ \varphi_{\partial} J \varphi_{\partial} {\rm d} A.
\]
Since $w$ is an arbitrary vector field, we can multiply $w$ with test functions that converge to the Dirac delta function, and the identity remains true. We can thus remove the integral and get the identity \ref{boundary_piola}.

Taking $w = \rho u$, we have that $w \circ \varphi J \varphi = \varrho \dot \varphi$. Piola's identity gives that: 
\[
N \cdot \left( \varrho( \nabla_x  \varphi ^{-1} \circ \varphi \cdot \dot \varphi) \right)  = \left( n \cdot (\rho  u ) \right)  \circ \varphi_{\partial} J \varphi_{\partial}.
\]
On the other hand, by definition 
\[
\mathfrak{j}_ \varrho = j_\rho \circ \varphi_{\partial} J  \varphi_{\partial} .
\]

We conclude that on the boundary
\[
\left(( \nabla_x  \varphi ^{-1} \circ \varphi \cdot \dot \varphi) \varrho  \right)\cdot N = -\mathfrak{j} _ \varrho  \quad  \Longleftrightarrow \quad \rho  u \cdot n =- j_ \rho .
\]

The same proof shows
\[
\left(( \nabla  \varphi_x ^{-1} \circ \varphi \cdot \dot \varphi) S  \right)\cdot N = -\mathfrak{j} _ S \quad  \Longleftrightarrow \quad s  u \cdot n =- j_ s  .
\]

\noindent{\bf (4)} 
Finally, note that 
\[
\left( \varrho \frac{\partial \mathfrak{L}}{\partial \varrho} + S \frac{\partial \mathfrak{L}}{\partial S} - \mathfrak{L}\right)  \nabla  \varphi ^{-1} \circ \varphi = \left( \rho  \frac{\partial \mathfrak{l}}{\partial \rho  }+s  \frac{\partial \mathfrak{l}}{\partial s  }- \mathfrak{l} \right) \circ \varphi J \varphi  \nabla  \varphi ^{-1} \circ \varphi =  -\frac{\partial \mathfrak{L}}{\partial \nabla \varphi }
\]
and that 
\[
N \cdot  (\nabla  \varphi ^{-1} \circ \varphi \cdot \dot \varphi )\frac{\partial \mathfrak{L}}{\partial \dot  \varphi } = N \cdot  (\nabla  \varphi ^{-1} \circ \varphi \cdot \dot \varphi )\frac{\partial \mathfrak{l}}{\partial u} \circ \varphi J \varphi = (n \cdot u) \frac{\partial \mathfrak{l}}{\partial u} \circ \varphi_{\partial} J \varphi_{\partial}
\] 
which is equation \ref{boundary_piola} with $w = \frac{\partial \mathfrak{l}}{\partial u} u$.
By definition 
\[
\mathfrak{J} = J \circ \varphi_{\partial} J \varphi_{\partial},
\]
we conclude that on the boundary
\[
N \cdot \left( \frac{\partial \mathfrak{L}}{\partial \nabla \varphi }   + (\nabla  \varphi ^{-1} \circ \varphi \cdot \dot \varphi )\frac{\partial \mathfrak{L}}{\partial \dot  \varphi } + \left( \varrho \frac{\partial \mathfrak{L}}{\partial \varrho} + S \frac{\partial \mathfrak{L}}{\partial S} - \mathfrak{L}\right)  \nabla  \varphi ^{-1} \circ \varphi \right) =-\mathfrak{J}
\]
\[
\Leftrightarrow
\]
\[
(u \cdot n) \frac{\partial \mathfrak{l} }{\partial u} =-J.
\]

\addcontentsline{toc}{section}{References}

\bibliographystyle{unsrtnat}
\bibliography{main}

\begin{thebibliography}{65}
\providecommand{\natexlab}[1]{#1}
\providecommand{\url}[1]{\texttt{#1}}
\expandafter\ifx\csname urlstyle\endcsname\relax
  \providecommand{\doi}[1]{doi: #1}\else
  \providecommand{\doi}{doi: \begingroup \urlstyle{rm}\Url}\fi

\bibitem[Holm et~al.(1998)Holm, Marsden, and Ratiu]{HoMaRa1998}
D.~D. Holm, J.~E. Marsden, and T.~S. Ratiu.
\newblock The {E}uler--{P}oincar{\'e} equations and semidirect products with
  applications to continuum theories.
\newblock \emph{Advances in Mathematics}, 137\penalty0 (1):\penalty0 1--81,
  1998.

\bibitem[{Gay-Balmaz} and Yoshimura(2017{\natexlab{a}})]{GBYo2017a}
F.~{Gay-Balmaz} and H.~Yoshimura.
\newblock A {L}agrangian variational formalism for nonequilibrium
  thermodynamics. {P}art {I}: discrete systems.
\newblock \emph{J. Geom. Phys.}, 111:\penalty0 169--193, 2017{\natexlab{a}}.

\bibitem[{Gay-Balmaz} and Yoshimura(2017{\natexlab{b}})]{GBYo2017b}
F.~{Gay-Balmaz} and H.~Yoshimura.
\newblock A {L}agrangian variational formalism for nonequilibrium
  thermodynamics. {P}art {II}: continuum systems.
\newblock \emph{J. Geom. Phys.}, 111:\penalty0 194--212, 2017{\natexlab{b}}.

\bibitem[Morrison and Greene(1980)]{MoGr1980}
P.~Morrison and J.~Greene.
\newblock Noncanonical {H}amiltonian density formulation of hydrodynamics and
  ideal magnetohydrodynamics.
\newblock \emph{Phys. Rev. Letters}, 45:\penalty0 790--794, 1980.

\bibitem[Holm and Kupershmidt(1983{\natexlab{a}})]{HoKu1983}
D.~D. Holm and B.~A. Kupershmidt.
\newblock Noncanonical {H}amiltonian formulation of ideal magnetohydrodynamics.
\newblock \emph{Physica D: Nonlinear Phenomena}, 7\penalty0 (1-3):\penalty0
  330--333, 1983{\natexlab{a}}.

\bibitem[Dzyaloshinskii and Volovick(1980)]{DzVo1980}
I.~E. Dzyaloshinskii and G.~E. Volovick.
\newblock Poisson brackets in condensed matter physics.
\newblock \emph{Ann. Phys.}, 125:\penalty0 67--97, 1980.

\bibitem[Marsden and Weinstein(1983)]{MaWe1983}
J~E. Marsden and A.~Weinstein.
\newblock {C}oadjoint orbits, vortices, and {C}lebsch variables for
  incompressible fluids.
\newblock \emph{Physica D: Nonlinear Phenomena}, 7\penalty0 (1):\penalty0
  305--323, 1983.

\bibitem[Marsden et~al.(1984)Marsden, Ratiu, and Weinstein]{MaRaWe1984}
J.~E. Marsden, T.~S. Ratiu, and A.~Weinstein.
\newblock Semidirect product and reduction in mechanics.
\newblock \emph{Trans. Amer. Math. Soc.}, 281:\penalty0 147--177, 1984.

\bibitem[Grmela(1984)]{Gr1984}
M.~Grmela.
\newblock Bracket formulation of dissipative fluid mechanics equations.
\newblock \emph{Phys. Lett. A}, 102:\penalty0 355--358, 1984.

\bibitem[Kaufman(1984)]{Ka1984}
A.~Kaufman.
\newblock Dissipative {H}amiltonian systems: A unifying principle.
\newblock \emph{Phys. Lett. A}, 100:\penalty0 419--422, 1984.

\bibitem[Morrison(1984{\natexlab{a}})]{Mo1984a}
P.~Morrison.
\newblock Bracket formulation for irreversible classical fields.
\newblock \emph{Phys. Lett. A}, 100:\penalty0 423--427, 1984{\natexlab{a}}.

\bibitem[Morrison(1984{\natexlab{b}})]{Mo1984b}
P.~Morrison.
\newblock Some observations regarding brackets and dissipation.
\newblock Technical report, University of California, Berkeley,
  1984{\natexlab{b}}.

\bibitem[Grmela and \"Ottinger(1997)]{GrOtt1997}
M.~Grmela and H.-C. \"Ottinger.
\newblock Dynamics and thermodynamics of complex fluids. {I}. {D}evelopment of
  a general formalism.
\newblock \emph{Phys. Rev. E}, 56:\penalty0 6620--6632, 1997.

\bibitem[\"Ottinger and Grmela(1997)]{OtGr1997}
H.-C. \"Ottinger and M.~Grmela.
\newblock Dynamics and thermodynamics of complex fluids. {II.} {I}llustrations
  of a general formalism.
\newblock \emph{Phys. Rev. E}, 56:\penalty0 6633--6655, 1997.

\bibitem[Holm et~al.(1985)Holm, Marsden, Ratiu, and Weinstein]{HoMaRaWe1985}
D.~D. Holm, J.~E. Marsden, T.~S. Ratiu, and A.~Weinstein.
\newblock Nonlinear stability of fluid and plasma equilibria.
\newblock \emph{Physics Reports}, 123\penalty0 (1-2):\penalty0 1--116, 1985.

\bibitem[Webb and Mace(2014)]{webb2014noether}
G.~M. Webb and R.~L. Mace.
\newblock {N}oether’s theorems and fluid relabelling symmetries in
  magnetohydrodynamics and gas dynamics.
\newblock \emph{J. Phys. A Math. Theoret. article No. JPHYSA-101-057, available
  at http://arxiv. org/abs/1403.3133 (submitted on March 10)}, 2014.

\bibitem[Webb et~al.(2014{\natexlab{a}})Webb, Dasgupta, McKenzie, Hu, and
  Zank]{webb2014locala}
G.~M. Webb, B.~Dasgupta, J.~F. McKenzie, Q.~Hu, and G.~P. Zank.
\newblock Local and nonlocal advected invariants and helicities in
  magnetohydrodynamics and gas dynamics {I}: {L}ie dragging approach.
\newblock \emph{Journal of Physics A: Mathematical and Theoretical},
  47\penalty0 (9):\penalty0 095501, 2014{\natexlab{a}}.

\bibitem[Webb et~al.(2014{\natexlab{b}})Webb, Dasgupta, McKenzie, Hu, and
  Zank]{webb2014localb}
G.~M. Webb, B.~Dasgupta, J.~F. McKenzie, Q.~Hu, and G.~P. Zank.
\newblock Local and nonlocal advected invariants and helicities in
  magnetohydrodynamics and gas dynamics {II}: {N}oether's theorems and
  {C}asimirs.
\newblock \emph{Journal of Physics A: Mathematical and Theoretical},
  47\penalty0 (9):\penalty0 095502, 2014{\natexlab{b}}.

\bibitem[Gay-Balmaz et~al.(2012)Gay-Balmaz, Marsden, and Ratiu]{GBMaRa2012}
F.~Gay-Balmaz, J.~E. Marsden, and T.~S. Ratiu.
\newblock Reduced variational formulations in free boundary continuum
  mechanics.
\newblock \emph{J. Nonlin. Sci.}, 22(4):\penalty0 463--497, 2012.

\bibitem[Makaremi-Esfarjani and Najafi-Yazdi(2022)]{Makaremi2022}
P.~Makaremi-Esfarjani and A.~Najafi-Yazdi.
\newblock Characteristic boundary conditions for magnetohydrodynamic equations.
\newblock \emph{Computers \& Fluids}, 241:\penalty0 105461, 2022.
\newblock ISSN 0045-7930.
\newblock \doi{https://doi.org/10.1016/j.compfluid.2022.105461}.
\newblock URL
  \url{https://www.sciencedirect.com/science/article/pii/S0045793022001116}.

\bibitem[Meier et~al.(2012)Meier, Glasser, Lukin, and Shumlak]{Meier2012}
E.~T. Meier, A.~H. Glasser, V.~S. Lukin, and U.~Shumlak.
\newblock Modeling open boundaries in dissipative mhd simulation.
\newblock \emph{Journal of Computational Physics}, 231\penalty0 (7):\penalty0
  2963--2976, 2012.
\newblock ISSN 0021-9991.
\newblock \doi{https://doi.org/10.1016/j.jcp.2012.01.003}.
\newblock URL
  \url{https://www.sciencedirect.com/science/article/pii/S0021999112000137}.

\bibitem[{\"O}ttinger(2006)]{Ot2006}
H.~C. {\"O}ttinger.
\newblock Nonequilibrium thermodynamics for open systems.
\newblock \emph{Physical Review E—Statistical, Nonlinear, and Soft Matter
  Physics}, 73\penalty0 (3):\penalty0 036126, 2006.

\bibitem[Badlyan et~al.(2018)Badlyan, Maschke, Beattie, and
  Mehrmann]{BaMaBeMe2018}
A.~Moses Badlyan, B.~Maschke, C.~Beattie, and V.~Mehrmann.
\newblock Open physical systems: from {GENERIC} to port-{H}amiltonian systems,
  2018.
\newblock URL \url{https://arxiv.org/abs/1804.04064}.

\bibitem[Rashad et~al.(2021{\natexlab{a}})Rashad, Califano, Schuller, and
  Stramigioli]{rashad2021portA}
R.~Rashad, F.~Califano, F.~P. Schuller, and S.~Stramigioli.
\newblock Port-{H}amiltonian modeling of ideal fluid flow: Part {I}.
  {F}oundations and kinetic energy.
\newblock \emph{Journal of Geometry and Physics}, 164:\penalty0 104201,
  2021{\natexlab{a}}.

\bibitem[Rashad et~al.(2021{\natexlab{b}})Rashad, Califano, Schuller, and
  Stramigioli]{rashad2021portB}
R.~Rashad, F.~Califano, F.~P. Schuller, and S.~Stramigioli.
\newblock Port-{H}amiltonian modeling of ideal fluid flow: Part {II}.
  {C}ompressible and incompressible flow.
\newblock \emph{Journal of Geometry and Physics}, 164:\penalty0 104199,
  2021{\natexlab{b}}.

\bibitem[Lohmayer and Leyendecker(2022)]{lohmayer2022exergetic}
M.~Lohmayer and S.~Leyendecker.
\newblock Exergetic port-{H}amiltonian systems: {N}avier-{S}tokes-{F}ourier
  fluid.
\newblock \emph{IFAC-PapersOnLine}, 55\penalty0 (18):\penalty0 74--80, 2022.

\bibitem[Mora et~al.(2023)Mora, Le~Gorrec, Matignon, and
  Ramirez]{mora2023irreversible}
L.~A. Mora, Y.~Le~Gorrec, D.~Matignon, and H.~Ramirez.
\newblock Irreversible port-{H}amiltonian modelling of 3{D} compressible
  fluids.
\newblock \emph{IFAC-PapersOnLine}, 56\penalty0 (2):\penalty0 6394--6399, 2023.

\bibitem[Cardoso-Ribeiro et~al.(2021)Cardoso-Ribeiro, Matignon, and
  Lef{\`e}vre]{cardoso2021dissipative}
F.~L. Cardoso-Ribeiro, D.~Matignon, and L.~Lef{\`e}vre.
\newblock Dissipative shallow water equations: A port-{H}amiltonian
  formulation.
\newblock \emph{IFAC-PapersOnLine}, 54\penalty0 (19):\penalty0 167--172, 2021.

\bibitem[Haine and Matignon(2021)]{haine2021incompressible}
G.~Haine and D.~Matignon.
\newblock Incompressible {N}avier-{S}tokes equation as port-{H}amiltonian
  systems: velocity formulation versus vorticity formulation.
\newblock \emph{IFAC-PapersOnLine}, 54\penalty0 (19):\penalty0 161--166, 2021.

\bibitem[Rashad et~al.(2021{\natexlab{c}})Rashad, Califano, Brugnoli, Schuller,
  and Stramigioli]{rashad2021exterior}
R.~Rashad, F.~Califano, A.~Brugnoli, F.~P. Schuller, and S.~Stramigioli.
\newblock Exterior and vector calculus views of incompressible
  {N}avier-{S}tokes port-{H}amiltonian models.
\newblock \emph{IFAC-papersonline}, 54\penalty0 (19):\penalty0 173--179,
  2021{\natexlab{c}}.

\bibitem[Califano et~al.(2021)Califano, Rashad, Schuller, and
  Stramigioli]{califano2021geometric}
F.~Califano, R.~Rashad, F.~P. Schuller, and S.~Stramigioli.
\newblock Geometric and energy-aware decomposition of the {N}avier--{S}tokes
  equations: A port-{H}amiltonian approach.
\newblock \emph{Physics of fluids}, 33\penalty0 (4), 2021.

\bibitem[Arnold(1966)]{Ar1966}
V.~Arnold.
\newblock Sur la g{\'e}om{\'e}trie diff{\'e}rentielle des groupes de {L}ie de
  dimension infinie et ses applications {\`a} l'hydrodynamique des fluides
  parfaits.
\newblock In \emph{Annales de l'institut Fourier}, volume~16, pages 319--361,
  1966.

\bibitem[Marsden and Hughes(2012)]{MaHu1983}
J.~E. Marsden and T.~J.~R. Hughes.
\newblock \emph{Mathematical foundations of elasticity}.
\newblock Courier Corporation, 2012.

\bibitem[Gay-Balmaz(2024)]{GB2024}
F.~Gay-Balmaz.
\newblock General relativistic {L}agrangian continuum theories part {I}:
  reduced variational principles and junction conditions for hydrodynamics and
  elasticity.
\newblock \emph{Journal of Nonlinear Science}, 34\penalty0 (3):\penalty0 46,
  2024.

\bibitem[{Gay-Balmaz} and Yoshimura(2019)]{GBYo2019a}
F.~{Gay-Balmaz} and H.~Yoshimura.
\newblock From {L}agrangian mechanics to nonequilibrium thermodynamics: a
  variational perspective.
\newblock \emph{Entropy}, 21:\penalty0 8, 2019.

\bibitem[Seliger and Whitham(1968)]{SeWh1968}
R.~L. Seliger and F.~R.~S. Whitham.
\newblock {V}ariational principles in continuum mechanics.
\newblock \emph{Proc. Roy. Soc. A}, 305:\penalty0 1--25, 1968.

\bibitem[{Gay-Balmaz} and Yoshimura(2020)]{GBYo2019b}
F.~{Gay-Balmaz} and H.~Yoshimura.
\newblock From variational to bracket formulations in nonequilibrium
  thermodynamics of simple systems.
\newblock \emph{J. Geom. Phys.}, 158:\penalty0 103812, 2020.

\bibitem[Gay-Balmaz and Putkaradze(2019)]{GBPu2019}
F.~Gay-Balmaz and V.~Putkaradze.
\newblock Geometric theory of ﬂexible and expandable tubes conveying ﬂuid:
  equations, solutions, and shock waves.
\newblock \emph{Journal of Nonlinear Science}, 29\penalty0 (2):\penalty0
  377--414, 2019.

\bibitem[Lee(2012)]{Le2013}
J.~M. Lee.
\newblock \emph{Smooth manifolds}.
\newblock Springer, 2012.

\bibitem[Bouchut et~al.(2009)Bouchut, Lambaerts, Lapeyre, and
  Zeitlin]{bouchut2009fronts}
F.~Bouchut, J.~Lambaerts, G.~Lapeyre, and V.~Zeitlin.
\newblock Fronts and nonlinear waves in a simplified shallow-water model of the
  atmosphere with moisture and convection.
\newblock \emph{Physics of Fluids}, 21\penalty0 (11), 2009.

\bibitem[Serrin(1959)]{Se1959}
J.~Serrin.
\newblock On the uniqueness of compressible fluid motions.
\newblock \emph{Archive for Rational Mechanics and Analysis}, 3:\penalty0
  271--288, 1959.

\bibitem[Dubos and Tort(2014)]{dubos2014equations}
T.~Dubos and M.~Tort.
\newblock Equations of atmospheric motion in non-eulerian vertical coordinates:
  Vector-invariant form and quasi-hamiltonian formulation.
\newblock \emph{Monthly Weather Review}, 142\penalty0 (10):\penalty0
  3860--3880, 2014.

\bibitem[Eldred and Gay-Balmaz(2020)]{eldred2020single}
C.~Eldred and F.~Gay-Balmaz.
\newblock Single and double generator bracket formulations of multicomponent
  fluids with irreversible processes.
\newblock \emph{Journal of Physics A: Mathematical and Theoretical},
  53\penalty0 (39):\penalty0 395701, 2020.

\bibitem[Gay-Balmaz(2019)]{gay2019variational}
F.~Gay-Balmaz.
\newblock A variational derivation of the thermodynamics of a moist atmosphere
  with rain process and its pseudoincompressible approximation.
\newblock \emph{Geophysical \& Astrophysical Fluid Dynamics}, 113\penalty0
  (5-6):\penalty0 428--465, 2019.

\bibitem[Makarieva et~al.(2017)Makarieva, Gorshkov, Nefiodov, Sheil, Nobre,
  Bunyard, Nobre, and Li]{makarieva2017equations}
A.~M. Makarieva, V.~G. Gorshkov, A.~V. Nefiodov, D.~Sheil, A.~D. Nobre,
  P.~Bunyard, P.~Nobre, and B.~Li.
\newblock The equations of motion for moist atmospheric air.
\newblock \emph{Journal of Geophysical Research: Atmospheres}, 122\penalty0
  (14):\penalty0 7300--7307, 2017.

\bibitem[Bott(2008)]{bott2008theoretical}
A.~Bott.
\newblock Theoretical considerations on the mass and energy consistent
  treatment of precipitation in cloudy atmospheres.
\newblock \emph{Atmospheric Research}, 89\penalty0 (3):\penalty0 262--269,
  2008.

\bibitem[Catry et~al.(2007)Catry, Geleyn, Tudor, B{\'e}nard, and
  Troj{\'a}akov{\'a}a]{catry2007flux}
B.~Catry, J.-F. Geleyn, M.~Tudor, P.~B{\'e}nard, and A.~Troj{\'a}akov{\'a}a.
\newblock Flux-conservative thermodynamic equations in a mass-weighted
  framework.
\newblock \emph{Tellus A: Dynamic Meteorology and Oceanography}, 59\penalty0
  (1):\penalty0 71--79, 2007.

\bibitem[Wacker and Herbert(2003)]{wacker2003continuity}
U.~Wacker and F.~Herbert.
\newblock Continuity equations as expressions for local balances of masses in
  cloudy air.
\newblock \emph{Tellus A: Dynamic Meteorology and Oceanography}, 55\penalty0
  (3):\penalty0 247--254, 2003.

\bibitem[Bannon(2002)]{bannon2002theoretical}
P.~R. Bannon.
\newblock Theoretical foundations for models of moist convection.
\newblock \emph{Journal of the Atmospheric Sciences}, 59\penalty0
  (12):\penalty0 1967--1982, 2002.

\bibitem[Helfferich(1981)]{helfferich1981theory}
F.~G. Helfferich.
\newblock Theory of multicomponent, multiphase displacement in porous media.
\newblock \emph{Society of Petroleum Engineers Journal}, 21\penalty0
  (01):\penalty0 51--62, 1981.

\bibitem[Quintard et~al.(2006)Quintard, Bletzacker, Chenu, and
  Whitaker]{quintard2006nonlinear}
M.~Quintard, L.~Bletzacker, D.~Chenu, and S.~Whitaker.
\newblock Nonlinear, multicomponent, mass transport in porous media.
\newblock \emph{Chemical Engineering Science}, 61\penalty0 (8):\penalty0
  2643--2669, 2006.

\bibitem[Chella et~al.(1998)Chella, Lasseux, and
  Quintard]{chella1998multiphase}
R.~Chella, D.~Lasseux, and M.~Quintard.
\newblock Multiphase, multicomponent fluid flow in homogeneous and
  heterogeneous porous media.
\newblock \emph{Revue de l'Institut Fran{\c{c}}ais du P{\'e}trole}, 53\penalty0
  (3):\penalty0 335--346, 1998.

\bibitem[Farkhutdinov et~al.(2020)Farkhutdinov, Gay-Balmaz, and
  Putkaradze]{FaGBPu2020}
T.~Farkhutdinov, F.~Gay-Balmaz, and V.~Putkaradze.
\newblock Geometric variational approach to the dynamics of porous media filled
  with incompressible fluid.
\newblock \emph{Acta Mechanica}, 431\penalty0 (9):\penalty0 3897--3924, 2020.

\bibitem[Gay-Balmaz and Putkaradze(2022)]{GBPu2022}
F.~Gay-Balmaz and V.~Putkaradze.
\newblock Variational geometric approach to the thermodynamics of porous media.
\newblock \emph{Zeitschrift f\"ur Angewandte Mathematik und Mechanik},
  102\penalty0 (11), 2022.

\bibitem[Gay-Balmaz and Ratiu(2009)]{gay2009geometric}
F.~Gay-Balmaz and T.~S. Ratiu.
\newblock The geometric structure of complex fluids.
\newblock \emph{Advances in Applied Mathematics}, 42\penalty0 (2):\penalty0
  176--275, 2009.

\bibitem[Holm(2002)]{holm2002euler}
D.~D. Holm.
\newblock Euler-{P}oincar{\'e} dynamics of perfect complex fluids.
\newblock In \emph{Geometry, mechanics, and dynamics}, pages 169--180.
  Springer, 2002.

\bibitem[Holm and Kupershmidt(1983{\natexlab{b}})]{holm1983poisson}
D.~D. Holm and B.~A. Kupershmidt.
\newblock Poisson brackets and {C}lebsch representations for
  magnetohydrodynamics, multifluid plasmas, and elasticity.
\newblock \emph{Physica D: Nonlinear Phenomena}, 6\penalty0 (3):\penalty0
  347--363, 1983{\natexlab{b}}.

\bibitem[Holm and Kupershmidt(1986)]{holm1986hydrodynamics}
D.~D. Holm and B.~A. Kupershmidt.
\newblock Hydrodynamics and electrohydrodynamics of adiabatic multiphase fluids
  and plasmas.
\newblock \emph{International journal of multiphase flow}, 12\penalty0
  (4):\penalty0 667--680, 1986.

\bibitem[Holm(1986)]{holm1986hamiltonian}
D.~D. Holm.
\newblock Hamiltonian dynamics of a charged fluid, including electro-and
  magnetohydrodynamics.
\newblock \emph{Physics Letters A}, 114\penalty0 (3):\penalty0 137--141, 1986.

\bibitem[Morrison(2005)]{morrison2005hamiltonian}
P.~J. Morrison.
\newblock Hamiltonian and action principle formulations of plasma physics.
\newblock \emph{Physics of plasmas}, 12\penalty0 (5), 2005.

\bibitem[Kanso et~al.(2007)Kanso, Arroyo, Tong, Yavari, Marsden, and
  Desbrun]{KaArToYaMaDe2007}
E.~Kanso, M.~Arroyo, Y.~Tong, A.~Yavari, J.~E. Marsden, and M.~Desbrun.
\newblock On the geometric character of stress in continuum mechanics.
\newblock \emph{Zeitschrift f\"ur Angewandte Mathematik und Mechanik},
  58:\penalty0 1--14, 2007.

\bibitem[Gilbert and Vanneste(2023)]{gilbert2023geometric}
A.~D. Gilbert and J.~Vanneste.
\newblock A geometric look at momentum flux and stress in fluid mechanics.
\newblock \emph{Journal of Nonlinear Science}, 33\penalty0 (2):\penalty0 31,
  2023.

\bibitem[Rashad et~al.(2023)Rashad, Brugnoli, Califano, Luesink, and
  Stramigioli]{rashad2023intrinsic}
R.~Rashad, A.~Brugnoli, F.~Califano, E.~Luesink, and S.~Stramigioli.
\newblock Intrinsic nonlinear elasticity: An exterior calculus formulation.
\newblock \emph{Journal of Nonlinear Science}, 33\penalty0 (5):\penalty0 84,
  2023.

\bibitem[Bresch et~al.(2019)Bresch, Gisclon, and
  Lacroix-Violet]{bresch2019navier}
D.~Bresch, M.~Gisclon, and I.~Lacroix-Violet.
\newblock On {N}avier--{S}tokes--{K}orteweg and {E}uler--{K}orteweg systems:
  application to quantum fluids models.
\newblock \emph{Archive for Rational Mechanics and Analysis}, 233:\penalty0
  975--1025, 2019.

\bibitem[Kraus and Businger(1994)]{kraus1994atmosphere}
E.~B. Kraus and J.~A. Businger.
\newblock \emph{Atmosphere-{O}cean {I}nteraction}, volume~27.
\newblock Oxford University Press, 1994.

\end{thebibliography}

\end{document}